\documentclass[12pt]{article}
\title{Harmonic analysis on a finite homogeneous space}

\author{Fabio Scarabotti, Filippo Tolli}

\usepackage{amssymb}
\usepackage{amsmath}
\usepackage{latexsym}
\usepackage{amsthm}
\usepackage{latexsym,amsbsy,amsmath,amscd,amssymb}

\usepackage{amsmath}
\usepackage{amssymb}

 \newtheorem{definition}{Definition} [section]
 \newtheorem{remark}[definition]{Remark}
 \newtheorem{example}[definition]{Example}

 \newtheorem{proposition}[definition]{Proposition}
 \newtheorem{theorem}[definition]{Theorem}
 \newtheorem{corollary}[definition]{Corollary}
  \newtheorem{lemma}[definition]{Lemma}

\def\NN{{\mathbb{N}}}

\def\RR{{\mathbb{R}}}

\def\CC{{\mathbb{C}}}

\renewcommand{\dim}{\mbox{ dim }}

\begin{document}
\maketitle

\begin{abstract}  In this paper, we study harmonic analysis on finite homogeneous spaces
whose associated permutation representation decomposes with multiplicity. After a careful look at Frobenius
reciprocity and transitivity of induction, and the introduction of three types of spherical functions,
 we develop a theory of Gelfand Tsetlin bases for permutation representations. Then we study several concrete
 examples on the symmetric groups, generalizing the Gelfand pair of the Johnson scheme; we also consider statistical
  and probabilistic applications. After that, we consider the composition of two permutation representations, giving
  a non commutative generalization of the Gelfand pair associated to the ultrametric space; actually, we study the more general notion of crested product. Finally, we consider the exponentiation action, generalizing the decomposition of the Gelfand pair of the Hamming scheme; actually, we study a more general construction that we call wreath product of permutation representations, suggested by the study of finite lamplighter random walks. We give several examples of concrete decompositions of permutation representations and several explicit 'rules' of decomposition.
\footnote{{\it AMS 2002 Math. Subj. Class.Primary: 43A95. Secondary:20C15,20C30,20E22,43A90,}:  }
\end{abstract}

\tableofcontents

%%%%%%%%%%%%%%%%%%%%%%%%%%%%%%%%%%%%%%%%%%%%%%%%%%%%%%%%%%%%%%%%%%%%%%%%%%%%%%%%%%%%%%%%%%%%%%%%%%%
%%%%%%%%%%%%%%%%%%%%%%%%%%%%%%%%%%%%%%%%%%%%%%%%%%%%%%%%%%%%%%%%%%%%%%%%%%%%%%%%%%%%%%%%%%%%%%%%%%
\section{Introduction}
%%%%%%%%%%%%%%%%%%%%%%%%%%%%%%%%%%%%%%%%%%%%%%%%%%%%%%%%%%%%%%%%%%%%%%%%%%%%%%%%%%%%%%%%%%%%%%%%%%
%%%%%%%%%%%%%%%%%%%%%%%%%%%%%%%%%%%%%%%%%%%%%%%%%%%%%%%%%%%%%%%%%%%%%%%%%%%%%%%%%%%%%%%%%%%%%%%%%%%%%%%%
Let $G$ be a finite group acting transitively on a finite set $X$ and let $K$ be the stabilizer of a point $x_0 \in X$,
so that $X = G/K$. If $T$ is a  $G$-invariant operator on $L(X)$
then its eigenspaces are representations of the group $G$ and therefore one could expect that
 the representation theory of $G$ could be used to determine the spectral decomposition of the operator $T$. The case of a finite Gelfand pair, that is when $L(X)$ decomposes without multiplicity, has been extensively studied; \cite{CST3} is a forthcoming monograph which treats this case, both the theoretical aspects and the applications (another important book that covers this subject is A. Terras' book\cite{Terras}). The present paper is devoted to the study of general spaces that decomposes with multiplicity.

 If the invariant operator $T$ does not belong to the center of $\text{\rm Hom}_G(L(X),L(X))$,  then it  has a nontrivial action on some
 isotypic component and therefore, in order to diagonalize the operator, it is necessary to find an explicit
 orthogonal decomposition   of that isotypic component. Although there is not a canonical way to derive such a
   decomposition,  in several  cases one can find a {\it natural} one.
In this paper, by Harmonic Analysis we mean the search of such orthogonal decompositions, in an explicit form suitable for the spectral study of invariant operators. The main application is to the study of invariant random walks along the lines of P. Diaconis work \cite{Diaconis}. Another application is to the spectral analysis of statistical data, again originated from Diaconis work.

Now we give two examples.
Let $X$ be a finite $G$-homogeneous graph and put in each
vertex a lamp, which may be on or off. A lamplighter performs the
simple random walk on $X$ and when he moves from a vertex $x$ to a
vertex $y$ he change randomly the state the lamps in $x$ and $y$,
that is both the lamps may be turned on or off with equal
probability. This is an example of a random walk on a finite homogeneous space that decomposes with multiplicity. Its analysis is in \cite{ScaTol2}; in the last part of the present paper, we develop a general harmonic analysis on wreath products that has the representation theoretic results in \cite{ScaTol2} as a particular case.
In the second example, consider a model made up of $3$ urns and $3m$ numbered balls. Put in the first urn the balls numbered $1,2,\dotsc,m$, in the second urn the balls $m+1,m+2,\dotsc,2m$ and in the third urn the balls $2m+1,2m+2,\dotsc, 3m$. At each time, choose at random two balls that belongs to different urns and switch them. This is a Markov chain in a space with multiplicity, but its Markov operator belongs to the center of the commuting algebra \cite{ScarabottiLapBer}. Now consider the following variation: we can only switch balls between urns 1 and 3 and between urns 2 and 3. The resulting Markov operator is not in the center of the commutant algebra and its spectral analysis requires an explicit orthogonal decomposition of each isotypic component; the details are in subsection \ref{GZMabc} of the present paper.

The paper is divided into four sections. In the first section, we develop some general theory.  It is known that the multiplicity of an irreducible representation
 in $L(X)$ is equal to  the dimension of the $K$-invariant vectors. This is Frobenius reciprocity. We examine it from two different sides, and derive some general rules to decompose a permutation representation; in particular, we introduce three kinds of spherical functions. We also show that the transitivity of induction play a fundamental role in the search of ``general principles'' that could be used to find interesting intermediate decompositions. This part of the paper culminates with a subsection devoted to the theory of Gelfand-Tsetlin bases.
  If there is  a chain of subgroups
  $G\geq H_1\geq H_2\geq \cdots H_m = K$ satisfying the Gelfand-Tsetlin  condition  (see Definition \ref{GZ}),
then  there is (modulo scalars) a unique choice  for the basis  of the $K$-invariant vectors (adapted  to such a
 chain of subgroups): this basis  leads to a ``natural'' orthogonal decomposition of the isotypic component. Recent fundamental papers on Gelfand-Tsetlin bases for the group algebra of the symmetric group are \cite{Ok-Ver1,Ok-Ver2}.

In the other parts of the paper we treat some specific spaces. There
are three basics examples of finite Gelfand pairs: the Johnson
scheme, the finite ultrametric space and the Hamming scheme; a
general reference on these pairs is \cite{CST3} (other useful
references are \cite{BHG,CST,DS,DunklProcSymp,FT,Stanton}). We give
examples that generalizes these pairs and have multiplicities.

As a noncommutative analog of the Johnson scheme, in the second part
of the paper we study a family of Gelfand-Tsetlin bases on the
symmetric group. In particular, we determine an explicit basis for
the representation $S^{n-1,1}$ in the permutation module
$M^{a_1,a_2,\ldots,a_k}$ (Theorem \ref{vettori}). Then we give two
different decompositions of the isotypic component $2S^{n-1,1}$ of
$M^{n-2,1,1}$; this way we illustrate the general theory developed
in the previous section. Moreover, as concrete application, we give
an analysis of the results of an election for president and director
of an association. Since the outcome of the election is an element
of  $M^{2,1,1} = S^{(n)}\oplus 2S^{n-1,1}\oplus S^{n-2,2}\oplus
S^{n-2,1,1}$, Diaconis \cite{Diaconis,DiaconisW}  computed and
interpreted the projections onto the isotypic components . We refine
its analysis  giving two different orthogonal decompositions of
isotypic component $2S^{n-1,1}$ (adapted to 2 different
Gelfand-Tsetlin bases) and interpreting the corresponding
projections. We determine also a Gelfand-Tsetlin decompositions for
all irreducible representations that appear in the permutation
module corresponding to a three parts partition, i.e. $M^{a,b,c}$
(Theorem \ref{t;sabc}). This is obtained by the spectral
decomposition of some laplacians and is connected with the results
in \cite{DunklPacific,ScarabottiSabc}.

In the multiplicity free setting, a generalization of the ultrametric space is given by the composition of Gelfand pairs. It has been studied in \cite{CST2}; see also \cite{BPRS,HH}. Recently, R. Bailey and P.J. Cameron \cite{Bailey,Ba-Ca} have developed a theory of crested products that generalizes both the direct product and the composition of two permutation representations (actually, they work both on groups and on association schemes). They treat the case of a symmetric association scheme, that generalizes the case of a symmetric Gelfand pair. In the third part of the paper, we extend the Bailey-Cameron theory to the case of crested product of arbitrary (i.e. with multiplicity) permutation representations. Now we have to deal with a noncommutative algebra of bi-$K$-invariant functions, and we extend in this setting the notion of ideal partition from \cite{Bailey,Ba-Ca}. Then we give explicit decomposition rules for the composition action of a crested product and also ex!
 plicit formulas for the spherical functions. We also gives a series of examples that illustrate the theory.

The last part of the paper is devoted to the analysis of the exponentiation action. This may be considered a generalization of the usual Hamming scheme; see \cite{CST2,Mi}. Actually, we treat a more general case: we give the decomposition rules for a more general notion, that we call wreath product of permutation representations; it has as a particular case both the exponentiation action and the left regular action of a wreath product on itself. This generalization is suggested by our recent work on the harmonic analysis of finite lamplighter random walks; \cite{ScaTol2}. As basic references on wreath products and their representations, we refer to \cite{JK,Hu}.

The present paper may be considered as a companion of \cite{ScaTol}. In that paper, in order to give partial solution of an open problem in \cite{BDPX}, we studied finite Markov chains invariant with respect a nontransitive permutation group $G$. We made the hypothesis that each orbit was a space without multiplicity (a Gelfand pair), but multiplicities were originated by considering the same representation in different orbits. Now we have only one orbit but with multiplicity. It remains to study the general case, a nontransitive group with multiplicity allowed on each orbit. This would lead to a complete solution of the open problem in \cite{BDPX}; we plan to devote a future paper on this case.

The present paper may be also considered as a sequel of \cite{CST2} (a joint paper with T. Ceccherini-Silberstein) where, among other things, we studied the composition and the exponentiation of Gelfand pairs. Now we are considering more general constructions (crested and wreath products) for arbitrary (i.e. with multiplicity) spaces.

{\bf Notation} If $X$ is a finite set, we denote by $L(X)$ the
vector space of all complex valued functions defined on $X$. The
space $L(X)$ will be endowed with the scalar product $\langle
f_1,f_2 \rangle_{L(X)}=\sum\limits_{x\in X}f_1(x)\overline{f_2(x)}$,
$f_1,f_2\in L(X)$. If $x\in X$, $\delta_x$ is the Dirac function
centered at $x$; if $A\subseteq X$, $\mathbf{1}_X$ is the
characteristic function of $A$. Suppose that $Y$ is another finite
set. We will use the natural isomorphism $L(X\times Y)\cong
L(X)\otimes L(Y)$, where for $f_1\in L(X), f_2\in L(Y)$,
$(f_1\otimes f_2)(x,y)=f_1(x)f_2(y)$, for all $(x,y)\in X\times Y$.
We denote by $Y^X$ the set of all functions $\varphi:X\rightarrow
Y$. A partition $\mathcal{P}$ of $X$ will be denoted by writing
$X=\coprod\limits_{A\in\mathcal{P}}A$, that is the simbol $\coprod$
denotes a disjoint union. If $G$ is a group, $1_G$ denotes the
identity of $G$. If $(\rho,V)$ is a representation of $G$ and

\begin{equation}\label{defiso}
V=\bigoplus_{i\in\mathcal{I}}m_iW_i,
\end{equation}

where $W_i,i\in\mathcal{I}$, are irreducible, pairwise inequivalent
representation of $G$ and $m_i$ is the multiplicity of $W_i$ in $V$,
then we will say that \eqref{defiso} is the {\em isotypic}
decomposition of $V$, and that $m_iW_i$ are the isotypic components.
If $U$ is a vector space and $u_1,u_2,\dotsc,u_m\in U$, then
$\langle u_1,u_2,\dotsc,u_m\rangle$ is the subspace generated by
$u_1,u_2,\dotsc,u_m$. If $U,V$ are unitary spaces, then
$\text{Hom}(U,V)$ (the space of all linear operators from $U$ to
$V$) will be endowed with the Hilbert-Schmidt scalar product: if
$T,S\in \text{Hom}(U,V)$ and $u_1,u_2,\dotsc,u_n$ is an orthonormal
basis in $U$, then $\langle
S,T\rangle_{\text{HS}}=\sum\limits_{k=1}^n\langle
Su_k,Tu_k\rangle_U$.

%%%%%%%%%%%%%%%%%%%%%%%%%%%%%%%%%%%%%%%%%%%%%%%%%%%%%%%%%%%%%%%%%%%%%%%%%%%%%%%%%
%%%%%%%%%%%%%%%%%%%%%%%%%%%%%%%%%%%%%%%%%%%%%%%%%%%%%%%%%%%%%%%%%%%%%%%%%%%%%%%%%%
\section{General theory}
%%%%%%%%%%%%%%%%%%%%%%%%%%%%%%%%%%%%%%%%%%%%%%%%%%%%%%%%%%%%%%%%%%%%%%%%%%%%%
%%%%%%%%%%%%%%%%%%%%%%%%%%%%%%%%%%%%%%%%%%%%%%%%%%%%%%%%%%%%%%%%%%%%%%%%%%%%%%

As general references on representation theory of finite groups, we refer to \cite{NS,Serre,Simon,Sternberg}. All the representations we consider are unitary.

%%%%%%%%%%%%%%%%%%%%%%%%%%%%%%%%%%%%%%%%%%%%%%%%%%%%%%%%%%%%%%%%%%%%%%%%%%%%%%%%%
%%%%%%%%%%%%%%%%%%%%%%%%%%%%%%%%%%%%%%%%%%%%%%%%%%%%%%%%%%%%%%%%%%%%%%%%%%%%%%%%%%
\subsection{Frobenius reciprocity for a permutation representation}
%%%%%%%%%%%%%%%%%%%%%%%%%%%%%%%%%%%%%%%%%%%%%%%%%%%%%%%%%%%%%%%%%%%%%%%%%%%%%
%%%%%%%%%%%%%%%%%%%%%%%%%%%%%%%%%%%%%%%%%%%%%%%%%%%%%%%%%%%%%%%%%%%%%%%%%%%%%%

\begin{definition}{\rm
Let $G$ be a finite group and $(\rho,V)$ a representation of
$G$.  The {\it commutant of $V$} is the algebra $Hom_G(V,V)$ of all linear operators
intertwining $V$ with  itself.}
\end{definition}

\begin{theorem}
Suppose that $V =  \bigoplus_{\rho \in I}m_\rho W_\rho$ is the isotypic decomposition of $V$. Then

\begin{equation}\label{dec}
Hom_G(V,V) \cong \bigoplus_{\rho \in I} M_{m_\rho,m_\rho}(\CC)
\end{equation}

as algebras. In particular, $\dim \text{Hom}_G(V,V) =\sum_{\rho \in I} m_\rho^2$.

\end{theorem}

\begin{proof}(Sketch)
Suppose that $m_\rho W_\rho = W_\rho^1\oplus W_\rho^2\oplus \cdots \oplus W_\rho^{d_\rho}$ is an orthogonal decomposition
of the {\it isotypic component} $m_\rho W_\rho$ for any $\rho \in I$.
Then, using Schur's lemma,  one can choose a basis $\{T_{i,j}^\rho: \rho \in I, i,j = 1, 2, \ldots,
 d_\rho\}$ of $Hom_G(V,V)$  with the following properties
\begin{itemize}
\item{
$\text{Ker}T^\rho_{i,j}=(W_\rho^j)^\bot$}
\item{$\text{Ran}T_{i,j}^\rho=W_\rho^i$}
\item{$T_{i,j}^\rho T^\rho_{j,k}=T^\rho_{i,k}.$}
\end{itemize}

Therefore any $T\in \text{Hom}_G(V,V)$, can be written uniquely as
$T=\sum_{\rho\in I}\sum_{i,j=1}^{m_\rho}\alpha_{i,j}^\rho
T_{i,j}^\rho$, and the map

\[
T\mapsto \oplus_{\rho\in I}(a_{i,j}^\rho)_{i,j=1,\dotsc,m_\rho}
\]

yield the isomorphism (\ref{dec}).

\end{proof}

\begin{corollary}\label{centro}
 A basis for the center
of $\text{Hom}_G(V,V)$ is given  by the orthogonal projections onto the isotypic components.
In particular, an operator $T$ is in
the center of $\text{Hom}(V,V)$ if and only if  every isotypic component $m_\rho W_\rho$ is an eigenspace for $T$.
\end{corollary}

Suppose that $G$ acts transitively on $X$, that $K$ is the stabilizer of a fixed point $x_0 \in X$ and
denote by $(\lambda,L(X))$ the corresponding permutation representation.

\begin{definition}
If $V$ is any $G-$representation,  we denote by
$V^K$ the subspace of $K-$invariant vectors in $V$.
\end{definition}

Suppose now that $(\rho,V)$ is irreducible and that $V^K$ is non--trivial. Set $d_\rho=\dim V$.
For any $v \in V^K$ define the linear map
 $T_v : V \to L(X)$ by setting

\begin{equation}\label{proiezione}
(T_v w)(g x_0) = \sqrt{\frac{d_\rho}{|X|}} \langle w, \rho(g)v\rangle_V,
\end{equation}

for any $w \in V$ and $g \in G$. Clearly, $T_vw\in L(X)$ is well defined because $v$ is $K$-invariant. Moreover, $T_v\in\text{Hom}_G(V,L(X))$: it is easy to see that $T_v\rho(w)=\lambda(g)T_vw$.

\begin{proposition}\label{p;3maroni} With the notation above, one has:

\begin{equation}\label{e;moccio}
\langle T_u w, T_v z \rangle_{L(X)} = \langle  w,  z \rangle_{V} \langle  v,  u \rangle_{V},
\end{equation}

for all $v,u\in V^K$ and $w,z\in V$.
In particular
\begin{enumerate}
\item \label{3maroni1}{if $\|v\|_V = 1$ then $T_v$ is an isometric immersion of $V$ into
$L(X)$;}
\item\label{3maroni2}{i$\text{\rm Im}(T_u)$ is orthogonal to $\text{\rm Im}(T_v)$ if and only if $u$ is orthogonal to $v$.}
\end{enumerate}
\end{proposition}
\begin{proof}
Fix $u,v\in V^K$ and define a linear map $R: V \to V$ by setting:

\begin{equation}\label{Rw}
Rw = \frac{1}{|K|}\sum_{g \in G} \langle w,\rho(g)u\rangle_V \rho(g)v
\end{equation}

for all $w\in V$. It is clear that $R \in \text{Hom}_G(V)$ and since
$V$ is irreducible, $R = a I_V$ with $a \in \CC$ and $I_V$ the identity map. But $\text{Tr}(R)=|X|\langle v, u \rangle_V$, as it can be easily checked from \eqref{Rw}. Then $ad_\rho = |X|\langle v, u \rangle$ and
$R = \frac{|X|}{d_\rho} \langle v, u \rangle_V I_V$.

Therefore, if $w,z \in V$, we have
\[
\begin{split}
\langle T_u w, T_v z \rangle_{L(X)} & = \frac{1}{|K|}\sum_{g \in G}\langle  w,
\rho(g)u \rangle_{V} \overline{\langle  z,  \rho(g)v \rangle_{V}} \cdot \frac{d_\rho}{|X|}\\
& = \frac{d_\rho}{|X|} \langle Rw,z\rangle_V \\
& = \langle  w,  z \rangle_{V} \langle  v,  u \rangle_{V}.
\end{split}
\]
Then (\ref{3maroni1}) and (\ref{3maroni2}) immediately follows from
(\ref{e;moccio}).

\end{proof}

\begin{theorem}[Frobenius reciprocity for permutation representations]\label{3maroni3}
The map

\[
\begin{split}
V^K &\longrightarrow \text{\rm Hom}_G(V,L(X))\\
v &\longmapsto \frac{1}{\sqrt{d_\rho}}T_v
\end{split}
\]

is an antilinear, isometric isomorphism between the vector spaces
$V^K$ and $\text{\rm Hom}_G(V,L(X))$. In particular,
the multiplicity of $V$ in $L(X)$ is
equal to $\dim V^K$, the dimension of the subspace of $K$-invariant vectors in $V$.

\end{theorem}
\begin{proof}
If $\alpha, \beta \in \CC$, $u,v \in V^K$, then clearly $T_{\alpha u+ \beta v} =
\overline{\alpha} T_u + \overline{\beta} T_v$, that is the map $v \mapsto T_v$ is antilinear.
Now we show that this map is a bijection. If $T \in \text{\rm Hom}_G(V,L(X))$, then $V \ni w \mapsto (Tw)(x_0) \in \CC$ is
a linear map and therefore there exists $v \in V$ such that $(Tw)(x_0) =
\langle w, v \rangle_V$, for all $w \in V$. It follows that for all $w \in V$ and $g \in G$

\begin{equation}\label{e;mocciolo}
\begin{split} [Tw](g x_0) & = [\lambda(g^{-1})Tw](x_0) \\
\mbox{(because $T \in \text{\rm Hom}_G(V,L(X))$) \ \ } \ & =  [T\rho(g^{-1})w](x_0)\\
& = \langle \rho(g^{-1})w,v \rangle_V \\
& = \langle w,\rho(g)v \rangle_V.
\end{split}
\end{equation}

This shows that $T = \sqrt{\frac{\lvert X\rvert}{d_\rho}}T_v$. Moreover, taking $g = k \in K$ in (\ref{e;mocciolo}) we get that $v \in V^K$.
It is also clear that for $T \in \text{\rm Hom}_G(V,L(X))$  such a vector $v \in V^K$ is uniquely determined.

It remains to show that the map is isometric: if $w_1,w_2,\dotsc,w_{d_\rho}$ is an orthonormal basis in $V$, then by \eqref{e;moccio}

\[
\langle T_u,T_w \rangle_\text{HS}=\sum_{j=1}^{d_\rho}\langle T_u w_j,T_vw_j\rangle_V=d_\rho\langle v,w\rangle_V.
\]
\end{proof}

\begin{corollary}\label{c;mocciolo} The vectors $v_1, v_2,\ldots, v_m$ form an orthogonal basis for $V^K$ if and only if $$T_{v_1}V \oplus
T_{v_2}V \oplus \cdots \oplus T_{v_m}V$$ is an orthogonal decomposition of the $V$-isotypic component of
$L(X)$.
\end{corollary}

Orthogonality relations may be proved also for the general form of Frobenius reciprocity; see \cite{Nebe}.

%%%%%%%%%%%%%%%%%%%%%%%%%%%%%%%%%%%%%%%%%%%%%%%%%%%%%%%%%%%%%%%%%%%%
%%%%%%%%%%%%%%%%%%%%%%%%%%%%%%%%%%%%%%%%%%%%%%%%%%%%%%%%%%%%%%%%%%%%
\subsection{Spherical functions}\label{sphefunctsection}
%%%%%%%%%%%%%%%%%%%%%%%%%%%%%%%%%%%%%%%%%%%%%%%%%%%%%%%%%%%%%%%%%%%%
%%%%%%%%%%%%%%%%%%%%%%%%%%%%%%%%%%%%%%%%%%%%%%%%%%%%%%%%%%%%%%%%%%%%

We present now the description of the  commutant of the permutation representations in terms of bi-$K$-invariant functions.
For $f \in L(X)$  denote by  $\tilde{f}$ the function in $L(G)$ defined by setting

\begin{equation}\label{e;Kinv}
\tilde{f}(g) = f(gx_0)
\end{equation}

The map $f\rightarrow\widetilde{f}$ is a linear isomorphism between the vector spaces $L(X)$ and $L(G/K)$ (and between the vector space of $K$-invariant functions on $X$ and $L(K\backslash G/K)$).

\begin{theorem}\label{comm}
The commutant  $\text{Hom}_G(L(X),L(X))$   is isomorphic  to the algebra of bi-$K$-invariant
 functions $L(K\backslash G/K): \{f \in L(G): f(k_1gk_2) = f(g), \ \forall k_1,k_2\in K, \ g \in G\}$.
\end{theorem}

\begin{proof}(sketch)
Suppose that $T \in \text{Hom}_G(L(X),L(X))$ and that
$[Tf](x) = \sum_{y\in X}r(x,y)f(y)$.
In virtue of the $G$-invariance of $T$, we have  $r(gx,gy) = r(x,y)$ for all $g \in G$. Set $\psi(x) = r(x,x_0)$ and
for $\xi \in L(G)$ define $\check{\xi}(g) = \xi(g^{-1})$. Then it is easy to show  that
$T \mapsto (\tilde{\psi})^{\check{\ }}$ is the desired isomorphism.
\end{proof}

\begin{remark}
{\rm
Form the proof of Theorem \ref{comm}, one can also get the following formula

\begin{equation}\label{repker}
\widetilde{Tf}=\frac{1}{\lvert K\rvert}\widetilde{f}*\widetilde{\psi}.
\end{equation}

We will say that $\psi$ is the representing kernel of $T$.
}
\end{remark}

Suppose again that $L(X) =  \bigoplus_{\rho \in I}m_\rho V_\rho$.
We can summarize Theorems \ref{dec} and \ref{comm} by writing

\begin{equation}\label{e;antiistaminico}
\text{Hom}_G(L(X),L(X))\cong L((K\backslash G/K) \cong \bigoplus_{\rho \in I} M_{m_\rho,m_\rho}(\CC).
\end{equation}

 In the remaining part of this section, we construct an explicit isomorphism
between $L(K\backslash G/K)$  and  $\bigoplus_{\rho \in I} M_{m_\rho,m_\rho}(\CC)$.
We also introduce two other kinds of algebras that are worthwhile to study, with the relative spherical functions.\\

For any irreducible representation $\rho \in I$, select an orthonormal basis $\{v_1^\rho,v_2^\rho, \ldots, v_{m_\rho}^\rho\}$ in $V_\rho^K$.

\begin{definition}\label{sphericalfunctions}{\rm
\begin{enumerate}

\item
For $\rho \in I$, the matrix coefficients
\begin{equation}\label{omeostinks}
\phi_{i,j}^\rho(g) = \langle v_i^\rho, \rho(g)v_j^\rho\rangle_{V_\rho}
\end{equation}
with $i,j = 1,2,\ldots, m_\rho$, are the {\it spherical matrix coefficients} of $(\rho, V_\rho)$.

\item For $\rho\in I$, the coefficients $\phi_{i,i}^\rho$, $i=1,2,\dotsc,m_\rho$, are the {\em spherical functions} of $(\rho, V_\rho)$.

\item For $\rho\in I$, $\chi_\rho^K=\sum_{i=1}^{m_\rho}\phi_{i,i}^\rho$ is the {\em spherical character} of $(\rho, V_\rho)$.
\end{enumerate}
}
\end{definition}

The spherical matrix coefficients form an orthonormal basis for  $L(K\backslash G/K)$ (note that from (\ref{e;antiistaminico})  it follows
that $\dim L(K\backslash G/K)= \sum_{\rho \in I}m_\rho^2$). Using the orthogonality relation for matrix coefficients  we can immediately write the Fourier transform relative to the matrix coefficients
(\ref{omeostinks}): for $F \in  L(K\backslash G/K)$, $\rho \in I$ and $i,j = 1,2, \ldots, m_\rho$
\[
\widehat{f}_{i,j}(\rho) = \langle f, \phi_{i,j}^\rho\rangle_{L(G)} \ \ \ \ \ \mbox{spherical Fourier transform}
\]
and
\[
f(g) =\frac{1}{|G|}\sum_{\rho \in I}d_\rho\sum_{i,j = 1}^{m_\rho} \phi_{i,j}^\rho(g) \widehat{f}_{i,j}(\rho)
\ \ \ \ \mbox{Inversion  formula.}
\]
From the orthogonality relations of matrix coefficients, one can easily prove
the identity
\begin{equation}\label{2maroni}
 \phi_{i,j}^\rho* \phi_{h,k}^\sigma = \frac{|G|}{d_\rho}\delta_{j,h}\delta_{\rho,\sigma} \phi_{i,k}^\rho
 \end{equation}
 and this immediately yields the desired explicit isomorphism.

 \begin{theorem}
 The map
 \[
 \begin{array}{rcl}
  L(K\backslash G/K) & \longrightarrow & \bigoplus_{\rho \in  I}M_{m_\rho,m_\rho}(\CC)\\
  f & \mapsto & \oplus_{\rho \in  I}\left(\widehat{f}_{i,j}(\rho)\right)_{i,j = 1,2,\ldots, m_\rho}
  \end{array}
  \] is an isomorphism of algebras.
  \end{theorem}

  \begin{proof} Foe all $f, h \in L(K\backslash G/K)$, $\rho \in I$ and $i,j = 1,2, \ldots, m_\rho$,
  we have
  \[
  \widehat{(f*h)}_{i,j}(\rho) = \sum_{k = 1}^{m_\rho}\widehat{f}_{i,k}(\rho)\cdot \widehat{h}_{k,j}(\rho)
  \]
  as it easily follows from the inversion formula and (\ref{2maroni}).
  \end{proof}

We can define another algebra, namely $\mathcal{A} = \mbox{ span }\{\phi_{i,i}^\rho: \rho\in I, \ i = 1,2, \ldots, m_\rho\}$.
 Clearly, it depends on the choice of the bases $\{v_1^\rho, v_2^\rho, \ldots, v_{m_\rho}^\rho\}$, $\rho \in I$.

\begin{proposition}\label{p;algebraf11}

\begin{enumerate}

\item\label{f11punto1}{ $\mathcal{A}$ is a maximal Abelian subalgebra of $L(K\backslash G/K)$.}

\item\label{f11punto2}{The operator $E_i^\rho: L(X)\to L(X)$ defined by

\begin{equation}\label{proj1}
(E_i^\rho f)(gx_0) = \frac{d_\rho}{|G|}\langle \widetilde{f}, \lambda(g)\phi_{i,i}^\rho\rangle_{L(G)},
\end{equation}

$f\in L(X)$, $g\in G$ ($\widetilde{f}$ is as in (\ref{e;Kinv})), is the projection from $L(X)$ onto
$T_{v_i^\rho}V_\rho$, $\rho\in I, \ i = 1,2, \ldots, m_\rho$ ($T_{v_i^\rho}V_\rho$ as in Proposition \ref{p;3maroni}).}

\end{enumerate}

\end{proposition}

\begin{proof} (\ref{f11punto1}) The diagonal matrices form a maximal Abelian subalgebra in the
full matrix algebra $M_{n,n}(\CC)$.

(\ref{f11punto2}) Extend each basis of $V_\rho^K$ to an orthonormal basis  $\{v_1^\rho, v_2^\rho, \ldots, v_{d_\rho}^\rho\}$ of $V_\rho$, for all $\rho \in I$. Any function $f \in T_{v_k}^\sigma V_\sigma$ is the  $T_{v_k}-$image of a vector $\sum_{j = 1}^{d_\sigma}a_jv_j^\sigma\in V_\rho^K$, and then is of the form
\[
f(hx_0)= \sqrt{\frac{d_\sigma}{|X|}}\sum_{j = 1}^{d_\sigma}a_j\phi_{j,k}^{\sigma}(h).
\]

Therefore
\[
\begin{split}
(E_i^\rho f)(gx_0)&  = \frac{d_\rho}{|G|}\sum_{h \in G}f(hx_0)\overline{\phi_{i,i}^\rho(g^{-1}h)}\\
&=  \frac{d_\rho}{|G|}\sqrt{\frac{d_\sigma}{|X|}}\sum_{j = 1}^{d_\sigma}a_j[\phi_{j,k}^{\sigma}*\phi_{i,i}^{\rho}](g)\\
& = \sqrt{\frac{d_\sigma}{|X|}}\sum_{j = 1}^{d_\sigma}a_j\delta_{\sigma,\rho}\delta_{k,i}\phi_{j,i}^{\sigma}(g)\\
& = \left\{\begin{array}{ll}
f(gx_0)& \mbox{ if $\rho = \sigma$ and $i = k$}\\
0 & \mbox{ otherwise.}
\end{array}
\right.
\end{split}
\]
\end{proof}

\begin{remark}{\rm
Suppose that $T:L(X)\rightarrow L(X)$ is normal and $G$-invariant.
Then each eigenspace of $T$ is $G$-invariant, and therefore may be
decomposed into mutually orthogonal, irreducible subspaces. This way
one can get an orthogonal decomposition $L(X)=\bigoplus_{\rho\in
I}\bigoplus_{j=1}^{d_\rho}T_{v^\rho_j}V_\rho$ such that each
$T_{v^\rho_j}V_\rho$ is an eigenspace of $T$. From another point of
view, suppose that $\psi$ is the representing kernel of $T$. Set
$\psi^\diamond(g)=\overline{\psi(g^{-1})}$. Then $T$ is normal if
and only if $\psi*\psi^\diamond=    \psi^\diamond*\psi$. The last
identity is equivalent to say that the spherical Fourier transform
of $\psi$ is normal, for any $\rho\in I$. If this is the case, one
can diagonalize the spherical Fourier transform of $\psi$ by a
choice of a suitable orthonormal basis of $K$-invariant vectors in
each irreducible representation. Also this way one gets again a
diagonalization of $T$. }
\end{remark}

The last algebra that we introduce is $\mathcal{B} = \mbox{ span
}\{\chi_{\rho}^K: \rho\in I \}$.

\begin{proposition}
\begin{enumerate}
\item\label{centropunto1}{$\mathcal{B}$ is the center of  $L(K\backslash G/K)$.}
\item\label{centropunto2} {The linear operator $E^\rho: L(X)\to L(X)$ defined by
\[
(E^\rho f)(gx_0) = \frac{d_\rho}{|G|}\langle \widetilde{f}, \lambda(g)\chi_{\rho}^K\rangle_{L(G)},
\]

$f\in L(X)$, $g\in G$ ($\widetilde{f}$ is as in (\ref{e;Kinv})), is the projection from $L(X)$ onto
the isotypic component $m_\rho V_\rho$, $\rho\in I$.}
\item\label{centropunto3}{If $\chi_\rho$ is the character of $\rho$, then $\chi_\rho^K(g)
 = \frac{1}{|K|}\sum_{k\in K}\overline{\chi_\rho(kg)}$.}
\end{enumerate}
\end{proposition}
\begin{proof}  The first two points  follow immediately from the previous proposition and Remark \ref{centro}.
For the third point, observe that $\frac{1}{|K|}\sum_{k\in K}\rho(k)$ is the orthogonal projection
of $V_\rho$ onto the space $V^K$ of $K$-invariant vectors.
\end{proof}
\begin{corollary}\label{corolla}
Suppose that the multiplicity of $V_\rho$ in  $L(X)$ is equal to $d_\rho$. Then
$\chi_\rho^K \equiv \chi_\rho$ and
\[
E^\rho f = \frac{d_\rho}{|G|} \sum_{g\in G}\chi_\rho(g) \lambda (g) f
\]
with $f \in L(X)$.
\end{corollary}
\begin{proof}
Il is clear that $\chi_\rho^K \equiv \chi_\rho$.  Therefore,
\[
\begin{split}
E^\rho f(gx_0) &  = \frac{d_\rho}{|G|}\langle \widetilde{f}, \lambda(g)\chi_\rho\rangle_{L(G)}\\
& = \frac{d_\rho}{|G|}\sum_{h \in G}f(hx_0) \overline{\chi_\rho(g^{-1}h)}\\
\mbox{ as $\chi_\rho$ is central}& = \frac{d_\rho}{|G|}\sum_{h \in G}f(hx_0) \overline{\chi_\rho(hg^{-1})}\\
t^{-1} = hg^{-1} & = \frac{d_\rho}{|G|}\sum_{t \in G}f(t^{-1}gx_0) \overline{\chi_\rho(t^{-1})}\\
& =  \frac{d_\rho}{|G|}\sum_{t \in G}\lambda(t)f(gx_0) {\chi_\rho(t)}.
\end{split}
\]
\end{proof}

The computational aspects of Corollary \ref{corolla} were explored in \cite{DiaconisRockmore}.

\begin{remark}{\rm In the following, often we will write $gf$ instead of $\lambda(g) f$.}
\end{remark}

%%%%%%%%%%%%%%%%%%%%%%%%%%%%%%%%%%%%%%%%%%%%%%%%%%%%%%%%%%%%%%%%%%%%%%%%%%%%%%%%%%%%%%%%%%%%%%%%%%%%%%%%%%%%%%%%%

%%%%%%%%%%%%%%%%%%%%%%%%%%%%%%%%%%%%%%%%%%%%%%%%%%%%%%%
%%%%%%%%%%%%%%%%%%%%%%%%%%%%%%%%%%%%%%%%%
\subsection{The other side of Frobenius reciprocity}
%%%%%%%%%%%%%%%%%%%%%%%%%%%%%%%%%%%%
%%%%%%%%%%%%%%%%%%%%%%%%%%%%%%%%%%%%%%%%%%%%%%%%%%%%%

Usually, Frobenius reciprocity is stated as an explicit isomorphism $\text{Hom}_K\left(\text{Res}_K^GV,U\right)\cong\text{Hom}_G\left(V,\text{Ind}_K^GU\right)$, where $K\leq G$, $V$ is a $G$-representation and $U$ is a $K$ representation. In Theorem \ref{3maroni3} we have examined the special case in which $U$ is the trivial representation. But Frobenius reciprocity may be also stated as an explicit isomorphism
$\text{Hom}_K\left( U,\text{Res}_K^GV\right)\cong\text{Hom}_G\left(\text{Ind}_K^GU,V\right)$; see \cite{Bump}, Corollary 34.1.
This formulation of Frobenius reciprocity is particular useful when the irreducible representations of $G$ may be obtained as induced representations; this is the case of a wreath product (see section \ref{repwreath}).
We examine this side of Frobenius reciprocity in the particular case in which $V$ is a permutation representation.

Let $G$ be again a finite group acting transitively on $X$ and suppose that $H$ is a subgroup of $G$. Let $(\rho,W)$ be an $H$-representation. Set $\tau=\text{Ind}_H^G\rho$ and denote by $\lambda$ the permutation representation of $G$ on $X$. Let $S$ be a set of representatives for the right cosets of $H$ in $G$, that is $G=\coprod\limits_{s\in S}sH$; we suppose that $1_G\in S$.
We recall that $\text{Ind}_H^GV$ has two equivalent descriptions \cite{NS,Serre,Simon,Sternberg} (we will use both of them). In the first description it is formed by all vectors $\sum\limits_{s\in S}\tau(s)w_s$ such that $w_s\in W$ for every $s\in S$. In the second description, it is formed by all functions $F:G\rightarrow W$ such that: $F(gh)=\rho(h^{-1})F(g)$, for all $g\in G$ and $h\in H$.

\begin{proposition}[Frobenius reciprocity II]\label{Frobenius}
For $L\in\text{\rm Hom}_H(W,\text{\rm Res}_H^GL(X))$, define $\stackrel{\diamond}{L}$ by setting

\[
\left(\stackrel{\diamond}{L}\sum\limits_{s\in S}\tau(s)w_s\right)(x)=\frac{1}{\sqrt{\lvert S\rvert}}\sum_{s\in S}(Lw_s)(s^{-1}x),
\]

for every $\sum\limits_{s\in S}\tau(s)w_s\in\text{\rm Ind}_H^GW$ and $x\in X$.
Then $\stackrel{\diamond}{L}\in\text{\rm Hom}_G(\text{\rm Ind}_H^GW,L(X))$ and the map

\[
\begin{split}
\text{\rm Hom}_H(W,\text{\rm Res}_H^GL(X))&\longrightarrow\text{\rm Hom}_G(\text{\rm Ind}_H^GW,L(X))\\
L&\longmapsto \stackrel{\diamond}{L}
\end{split}
\]

is an isometric isomorphism.
\end{proposition}

\begin{proof}
Clearly, $L$ verifies the identity $L\rho(h)w=(Lw)(h^{-1}x)$, for all $x\in X,h\in H,w\in V$. Take any $\sum_{t\in S}\tau(t)w_t\in\text{Ind}_H^GW$ and fix $g\in G$. For every $t\in T$, let $s_t\in S$ and $h_t\in H$ be the unique elements such that $gt=s_th_t$. Then we have:

\[
\begin{split}
\left[\stackrel{\diamond}{L}\tau(g)\sum_{t\in S}\tau(t)w_t\right](x)=&\left[\stackrel{\diamond}{L}\sum_{t\in S}\tau(s_t)\rho(h_t)w_t\right](x)\\
=&\frac{1}{\sqrt{\lvert S\rvert}}\sum_{t\in S}[L\rho(h_t)w_t](s_t^{-1}x)\\
=&\frac{1}{\sqrt{\lvert S\rvert}}\sum_{t\in S}(Lw_t)(h_t^{-1}s^{-1}_tx)\\
=&\frac{1}{\sqrt{\lvert S\rvert}}\sum_{t\in S}(Lw_t)(t^{-1}g^{-1}x)\\
=&\left[\lambda(g)\stackrel{\diamond}{L}\sum_{t\in S}\tau(t)w_t\right](x).
\end{split}
\]

This shows that $\stackrel{\diamond}{L}\in\text{Hom}_G(\text{Ind}_H^GV,L(X))$.

For $L\in\text{Hom}_G(\text{Ind}_H^GW,L(X))$, define $\widetilde{L}$ in the following way.
If $w\in W$, set $w_s=0$ when $s\neq 1_G$ and $w_{1_G}=w$ and then set $\widetilde{L}w=L\sum\limits_{s\in S}\tau(s)w_s$. It is easy to check that $\widetilde{L}\in\text{Hom}_H(W,\text{Res}_H^GL(X))$. Moreover, $\widetilde{L}$ is just the restriction of $L$ to the subspace $\tau(1_G)W$ in $\text{Ind}_H^GW$, and therefore $L\sum\limits_{s\in S}\tau(s)w_s=\sum\limits_{s\in S}\lambda(s)Lw_s=\sum\limits_{s\in S}\lambda(s)\widetilde{L}w_s$ for any choice of $\sum\limits_{s\in S}\tau(s)w_s\in\text{Ind}_H^GW$; from this fact it follows that the map $L\mapsto \sqrt{\lvert S\rvert}\widetilde{L}$ is the inverse of $L\mapsto \stackrel{\diamond}{L}$.

It is also easy to check that $\langle \stackrel{\diamond}{L}_1,\stackrel{\diamond}{L}_2 \rangle_{HS}=\langle L_1,L_2\rangle_{HS}$.

\end{proof}

%%%%%%%%%%%%%%%%%%%%%%%%%%%%%%%%%%%%%%%%%%%%%%%%%%%%%%%%%%%%%%%%%%%%%%%%%%%%
%%%%%%%%%%%%%%%%%%%%%%%%%%%%%%%%%%%%%%%%%%%%%%%%%%%%%%%%%%%%%%%%%%%%%%%%%%%%
\subsection{Decomposition by transitivity of induction}\label{trainddec}
%%%%%%%%%%%%%%%%%%%%%%%%%%%%%%%%%%%%%%%%%%%%%%%%%%%%%%%%%%%%%%%%%%%%%%%%%%%
%%%%%%%%%%%%%%%%%%%%%%%%%%%%%%%%%%%%%%%%%%%%%%%%%%%%%%%%%%%%%%%%%%%%%%%%%%%

Now we show how transitivity of induction may be used to decompose $L(X)$ in terms of induced representations. Suppose that $K\leq H\leq G$, set $Y = G/H$, $Z= H/K$ and let $x_0 \in X$, $z_0\in Z$ be the points stabilized by $K$ and $y_0 \in Y$ the point stabilized
by $H$. Let $S$ be a set of transversal of $G/H$ and $T$ a set of transversal  of  $H/K$. Then
\[
G = \coprod_{s \in S}sH \ \ \ \ \ H = \coprod_{t \in T}tK \ \ \ \ \ G = \coprod_{s \in S}\coprod_{t \in T}stK.
\]
Clearly $X = \{stx_0: s \in S, t \in T\}$ and $Y = \{sy_0: s \in S\}$.
\begin{proposition}\label{pisurj} There exists a unique map $\pi: X \to Y$ with the following properties:
\begin{itemize}
\item{$\pi$ is surjective;}
\item{$\pi(x_0) = y_0$}
\item{$\pi$ is $G$-equivariant: $\pi(gx) = g\pi(x)$ for all $x \in X$ and $g \in G$.}
\end{itemize}
\end{proposition}
\begin{proof}
It is immediate to check that such a $\pi$ is given by: $\pi(stx_0) = sy_0$.
\end{proof}

 Let $V_0$ the trivial representation of $K$. Then, the transitivity of induction
\begin{equation}\label{sandiego}
L(X) = \text{\rm Ind}_K^G V_0 = \text{\rm Ind}_H^GL(Z)
\end{equation}

could be made explcit by saying that

\begin{equation}\label{sandiego2}
L(X) = \bigoplus_{y \in Y}L(\pi^{-1}(y)),
\end{equation}

where each $\pi^{-1}(y)$ is a copy of $Z$: we identify $Z$ with $\pi^{-1}(y_0)$ and if $y = sy_0$, then $\pi^{-1}(y) = sZ$.

If $L(Z) = \bigoplus_{\rho \in I}m_\rho W_\rho$ is the decomposition of $L(Z)$ in irreducible representations of $H$,
then $L(sZ) = \bigoplus_{\rho \in I}m_\rho s W_\rho$, where  for $f \in L(X)$, $(sf)(x) = f(s^{-1}x)$, for any $x \in X$
and we have

\begin{equation}\label{sandiego3}
\text{\rm Ind}_H^GW_\rho = \bigoplus_{s\in S}sW_\rho.
\end{equation}

In (\ref{sandiego3}) if $sy_0 = y$ then $sW_\rho$ is the subspace of $L(\pi^{-1}(y))$ built using the isomorphism
between $L(\pi^{-1}(y))$ and $L(Z)$ given by $L(Z)\ni f\rightarrow sf\in L(\pi^{-1}(y))$. In the following Proposition, we summarize these facts in an orthogonal form.\\

\begin{proposition}\label{genprinc}
For every $\rho\in I$, let $\bigoplus_{j =1}^{m_\rho}W_\rho^j$ be an explicit orthogonal decomposition  of the isotypic component $m_\rho W_\rho$ (i.e. $W_\rho^1\cong
W_\rho^2 \cong \cdots \cong W_\rho^{m_\rho}$). Set $W_\rho^i(\pi^{-1}(y)) = sW_\rho^i$ if $sy_0 = y$.
Then we have the orthogonal decomposition

\begin{equation}\label{sandiego4}
L(X) = \bigoplus_{\rho\in I}\bigoplus_{j = 1}^{m_\rho}\left[\bigoplus_{y \in Y}W_\rho^j(\pi^{-1}(y))\right]
\end{equation}

where
\[
\bigoplus_{y \in Y}W_\rho^j(\pi^{-1}(y)) =\text{\rm Ind}_H^GW_\rho^j.
\]
\end{proposition}

\begin{remark}{\rm
In Proposition \ref{genprinc} we have replaced $\stackrel{\diamond}{L}$ of Theorem \ref{Frobenius} with the map $L^\sharp$ given by

\[
\left(L^\sharp\sum\limits_{s\in S}\tau(s)w_s\right)(x)=(Lw_{s_0})(s_0^{-1}x)
\]

if $x\in s_0Z$. Now the map $L\mapsto L^\sharp$ is an injective homomorphism $\text{\rm Hom}_H(W,L(Z))\longrightarrow\text{\rm Hom}_G(\text{\rm Ind}_H^GW,L(X))$, which in general is {\em not} surjective (for instance, take $W=$ the trivial representation of $H$).

}
\end{remark}

Denote by $E^\rho_{j}:L(Z) \to W_\rho^j$ the projection operator and by $(e^\rho_j(z,z'))_{z,z' \in Z}$ the  corresponding representing matrix, that is
$E^\rho_jf(z) = \sum_{z'\in Z}e^\rho_{j}(z,z')f(z')$. Clearly, $(e^\rho_j(z,z'))_{z,z' \in Z}$  is $H$-invariant:  $e^\rho_{j}(hz,hz') = e^\rho_{j}(z,z')$
for all $z,z' \in Z$ and $h \in H$.

Therefore, the projection operator $sE^\rho_{j}s^{-1}: L(\pi^{-1}(y))\to W_{\rho,j}(\pi^{-1}(y))$,
where $s\in S$ and $sy_0=y$, is represented by the matrix
$(e^\rho_{j}(s^{-1}z,s^{-1}z'))_{z,z' \in Z}$. Indeed, if $x\in L(\pi^{-1}(y))$, $f \in L(\pi^{-1}(y))$, then

\[
\begin{split}
\left(sE^\rho_{j}s^{-1} f\right)(x) & = \sum_{z'\in Z}e^\rho_{j}(s^{-1}x,z')f(sz')\\
& =  \sum_{x'\in \pi^{-1}(y)}e^\rho_{j}(s^{-1}x,s^{-1}x')f(x').
\end{split}
\]

The projector operator onto $\text{Ind}_H^G W_{\rho,j}$ is given by $\bigoplus_{s\in S} sE^\rho_{j}s^{-1}$.

Suppose now that $W_\rho$ is contained in $L(Z)$ with multiplicity $d_\rho\equiv\dim W_\rho$. Then, by Corollary \ref{corolla}

\begin{equation}\label{gilo}
E^\rho f = \frac{d_\rho}{H}\sum_{h \in H}\chi_\rho(h) hf
\end{equation}

projects $L(X)$ onto the isotypic component $d_\rho W_\rho$.
Then the projection
$L(\pi^{-1}(y)) \to d_\rho s W_\rho$ is given by

\begin{equation}\label{projectioniso}
sE^\rho s^{-1} f =  \frac{d_\rho}{|H|} \sum_{h \in H}\chi_\rho(h)\cdot shs^{-1}f =
\frac{d_\rho}{|H|} \sum_{t \in sHs^{-1}}\chi_\rho(s^{-1}ts)\cdot tf.
\end{equation}

\begin{example}\label{exasymm}{\rm
We apply now the above theory to the particular case $G = S_n$, $K= S_{n-k}$ and $H = S_{n-k} \times S_k$.
 We have

\[
 X = S_n/S_{n-k} = \{(i_1,i_2, \ldots,i_k):  i_1, i_2, \ldots,i_k  \in \{1,2, \ldots n\}, i_j\neq i_l
\mbox{ for } j\neq l\},
\]

$Y$ is the family of all $k$-subsets of $\{1,2,\dotsc,n\}$ and

\[
Z = (S_{n-k} \times S_k)/S_{n-k} = \{(i_1,i_2, \ldots, i_k) : i_1, i_2, \ldots,i_k  \in \{1,2, \ldots k\}, i_j\neq i_l
\mbox{ for } j\neq l\}.
\]

Alternatively, we may see $X$ as the set of all injective functions $\xi:\{1,2,\dotsc,k\}\rightarrow\{1,2,\dotsc,n\}$: identify $(i_1,i_2,\dotsc,i_k)$ with $\xi$ given by $\xi(j)=i_j$. This way, $Z$ coincides with $S_k$ and if $A\in Y$ ($A$ is a $k$-subset) then  $\pi^{-1}(A)$ is the set of all $\xi:\{1,2,\dotsc,k\}\rightarrow A$.

Using the standard terminology for the symmetric group \cite{sagan}, we write $M^{n-k, 1^k} = L(S_n/S_{n-k})$ and
$M^{1^k} = L(Z)$.
We clearly have $M^{n-k, 1^k} = \text{Ind}_{S_{n-k}\times S_k}^{S_n}\left[S^{(n-k)}\times M^{1^k}\right]$.

We have the natural action of $S_n$ on $X$: if $\theta\in S_n$ and $(i_1,i_2,\dotsc,i_k)\in X$ then $\theta(i_1,i_2,\dotsc,i_k)=(\theta(i_1),\theta(i_2),\dotsc,\theta(i_k))$.
If $(i_1,i_2,\dotsc,i_k)$ is identified with a function $\xi$, this is equivalent to multiplication on the left: $\xi\mapsto \theta\xi$.
On $X$ there is also a natural action of $S_k$, that we denote by $r$: if $\sigma\in\ S_k$ then $r(\sigma) (i_1,i_2, \ldots, i_k)  = (i_{\sigma^{-1}(1)},
i_{\sigma^{-1}(2)}, \ldots, i_{\sigma^{-1}(k)})$.
This is  equivalent to a multiplication on the right: $\xi\mapsto \xi\sigma^{-1}$.
For $\lambda \vdash k$, denote by $S^\lambda$ the corresponding irreducible representation of $S_n$ and by $\chi_\lambda$ its character.
By (\ref{gilo}), $E^\lambda =\frac{d_\lambda}{k!}\sum_{\sigma \in S_k}\chi_\lambda(\sigma) \sigma$
is the projector  of $M^{1^k}$ onto $d_\lambda S^\lambda$.

\begin{proposition}\label{projsymm}
The orthogonal projection

\begin{equation*}
\begin{array}{lccc}
\widetilde{E}^\lambda:& M^{n-k,1^k}& \to & \text{\rm Ind}_{S_{n-k}\times S_k}^{S_n}\left[S^{(n-k)}\otimes d_\lambda S^\lambda\right]\\
& f& \mapsto &\widetilde{E}^\lambda f
\end{array}
\end{equation*}

is given by

\begin{equation}\label{symm}
(\widetilde{E}^\lambda f)(i_1,i_2, \ldots, i_k) = \frac{d_\lambda}{k!}
\sum_{\sigma\in S_k}\chi_\lambda(\sigma)f(i_{\sigma^{-1}(1)},
i_{\sigma^{-1}(2)}, \ldots, i_{\sigma^{-1}(k)}).
\end{equation}

\end{proposition}

\begin{proof}
Let $\Theta$ be a set of representatives of the cosets of $S_n/S_{n-k}$, that is $S_n = \coprod_{\theta \in \Theta}
\theta(S_{n-k}\times S_k)$. Then $\widetilde{E}^\lambda=\bigoplus_{\theta\in \Theta}\theta\widetilde{E}\theta^{-1}$, where each $\theta\widetilde{E}\theta^{-1}$ acts on $L(\theta Z)$. Therefore, we may fix $\theta\in\Theta$ and assume that $f\in L(\theta Z)$. If $(i_1,i_2,\dotsc,i_k)\in\theta Z$ and $\xi$ is the corresponding function $\xi:\{1,2,\dotsc,k\}\rightarrow \{\theta(1),\theta(2),\dotsc,\theta(k)\}$ then $\theta^{-1}\xi\in S_k$, and by the conjugacy invariance of $\chi_\lambda$, \eqref{projectioniso} becomes

\[
\begin{split}
\sum_{\sigma\in S_k}\chi_\lambda(\sigma)f(\theta\sigma^{-1}\theta^{-1}(i_1,i_2,\dotsc,i_k))=&\sum_{\sigma\in S_k}\chi_\lambda(\sigma)f(\theta\sigma^{-1}\theta^{-1}\xi)\\
=&
\sum_{\sigma\in S_k}\chi_\lambda(\sigma)f(\xi(\theta^{-1}\xi)^{-1}\sigma^{-1}(\theta^{-1}\xi))=\\
=&\sum_{\sigma\in S_k}\chi_\lambda(\sigma)f(\xi\sigma^{-1})\\
=&\sum_{\sigma\in S_k}\chi_\lambda(\sigma)f((i_{\sigma^{-1}(1)},i_{\sigma^{-1}(2)},\dotsc,i_{\sigma^{-1}(k)})).
\end{split}
\]

\end{proof}

For instance, when $k=2$ Proposition \ref{projsymm} yields the decomposition of $M^{n-2,1,1}$ into symmetric (i.e., $f(i,j)=f(j,i)$) and antisymmetric (i.e., $f(i,j)=-f(j,i)$) functions. Compare with \cite{MetropolisRota, MRS}; see also subsection \ref{decMn211}.
}
\end{example}

We end this section by giving an explicit rule to decompose the
induced representation in \eqref{sandiego3}. It will be fundamental
in the theory of Gelfand-Tsetlin bases. Assume all the notation in
Proposition \ref{pisurj}. Let $(\sigma,V)$ be an irreducible
representation of $G$ and $W$ an $H$-invariant, irreducible subspace
of $V$. That is, $(\rho,W)$ is an irreducible representation in
$\text{Res}^G_HV$ and $\rho(h)$ is the restriction of $\sigma(h)$ to
$W$, $h\in H$. If $w_0\in W$ is $K$-invariant, then, by mean of
\eqref{proiezione}, we can form two distinct intertwining operators:
$S_{w_0}:W\rightarrow L(Z)$ ($H$-invariant) and
$T_{w_0}:V\rightarrow L(X)$ ($G$-invariant). By \eqref{sandiego3},
$\text{Ind}_H^GS_{w_0}W$ is a well defined subspace of $L(X)$.

\begin{lemma}\label{traind}
The operator $T_{w_0}$ intertwines $V$ with the subspace $\text{\rm Ind}_H^GS_{w_0}W$.
\end{lemma}

\begin{proof}
Denote by $P_W:V\rightarrow W$ the orthogonal projection onto $W$.
Suppose that $x\in X$ and $x=shz_0$, with $s\in S$ and $h\in H$. By \eqref{proiezione}, for any $v\in V$ we have

\[
\begin{split}
(T_{w_0}v)(x)=(T_{w_0}v)(shx_0)=&\sqrt{\frac{d_\sigma}{\lvert X\rvert}}\langle v,\sigma(s)\rho(h)w_0\rangle_V\\
=&\sqrt{\frac{d_\sigma}{\lvert X\rvert}}\langle \sigma(s^{-1})v,\rho(h)w_0\rangle_V\\
=&\sqrt{\frac{d_\sigma}{\lvert X\rvert}}\langle P_W\sigma(s^{-1}) v,\rho(h)w_0\rangle_V\\
=&\sqrt{\frac{\lvert Z\rvert d_\sigma}{\lvert X\rvert d_\rho}}\left[S_{w_0}(P_W\sigma(s^{-1})v\right](hz_0)\\
(x=shz_0)\qquad=&\sqrt{\frac{\lvert Z\rvert d_\sigma}{\lvert X\rvert d_\rho}}\cdot\left\{s\left[S_{w_0}(P_W\sigma(s^{-1})v\right]\right\}(x),
\end{split}
\]

that is, $T_{w_0}v\in sS_{w_0}W$, and by \eqref{sandiego3} this shows that $T_{w_0}v\in\text{Ind}_H^GS_{w_0}W$.
\end{proof}

%%%%%%%%%%%%%%%%%%%%%%%%%%%%%%%%%%%%%%%%%%%%%%%%%%%%%%%%%%%%%%%%%%%%%
%%%%%%%%%%%%%%%%%%%%%%%%%%%%%%%%%%%%%%%%%%%%%%%%%%%%%%%%%%%%%%%%%%%%%
\subsection{Gelfand-Tsetlin bases}
%%%%%%%%%%%%%%%%%%%%%%%%%%%%%%%%%%%%%%%%%%%%%%%%%%%%%%%%%%%%%%%%%%%%%%%%%%%%
%%%%%%%%%%%%%%%%%%%%%%%%%%%%%%%%%%%%%%%%%%%%%%%%%%%%%%%%%%%%%%%%%%%%%%%%%

Let $G$ be a finite group, $K\leq G$ a subgroup and $X= G/K$ the corresponding homogeneous space.
Let $L(X) = \sum_{\rho \in \widehat{G}}m_\rho V_\rho$ be the decomposition of $L(X)$ into
irreducible representations. The main result of this section is an explicit orthogonal decomposition
of the isotypic component $m_\rho V_\rho$ under a particular condition that we now describe.

\begin{definition}[Gelfand-Tsetlin condition]\label{GZ}{\rm A chain of subgroups  $G = H_m\geq H_{m-1} \geq \cdots H_2 \geq H_1 = K$
 from $K$ to $G$ satisfies the {\it Gelfand-Tsetlin condition} if for every irreducible $H_j$-representation  $V$
  that contains nontrivial $K$-invariant vectors (and therefore $V$  is contained in $L(H_j/K)$),
the decomposition of $\text{Ind}_{H_j}^{H_{j+1}}V$ into irreducible $H_{j+1}-$representations is multiplicity free,
$j= 1,2, \ldots m-1$.}
\end{definition}

\begin{remark} {\rm From Frobenius reciprocity theorem, it follows that the above condition is equivalent to say that
if $W,V$ are irreducible representations of $H_{j+1}$ and $H_{j}$
respectively, and they both contain nontrivial $K$-invariant vectors, then
the multiplicity of $V$ in $\text{Res}_{H_j}^{H_{j+1}}W$ is $\leq 1$.}
\end{remark}

Let $H_1,H_2, \ldots H_m$ be as before. We denote by $\mathcal{B}_j$ the set of the irreducible inequivalent representations
of $H_j$ containing nontrivial $K-$invariant vectors. Clearly $\mathcal{B}_1 = \{\iota_K\}$, with $\iota_K$ the trivial representations
of $K$.

\begin{definition} {\rm The {\it Brattelli diagram} is the graph whose
vertex set is $\mathcal{B} =\cup_{j = 1}^k\mathcal{B}_j$; the edge
set is constituted by all  the pairs $\{\rho,\sigma\}$, with
 $\sigma \in \mathcal{B}_j$ and $\rho  \in \mathcal{B}_{j+1}$, such that
$\text{Res}_{H_j}^{H_{j+1}}\rho$ contains $\sigma$ (and this is equivalent to  say that $\text{Ind}_{H_j}^{H_{j+1}}\sigma$ contains
$\rho$.)}
\end{definition}

A {\it path} in the Brattelli diagram is a sequence $\rho_1, \rho_2, \ldots, \rho_m$ with $\rho_j \in \mathcal{B}_j$
for $j= 1,2, \ldots m$ and $\{\rho_j, \rho_{j+1}\}$ an edge for $j= 1,2, \ldots m-1$.

Suppose that the Gelfand-Tsetlin condition  holds.
Let $\rho$ be an irreducible representation of $G$ contained in $L(X)$ and consider  a path ${C}: \rho_1, \rho_2,
\ldots, \rho_m = \rho$  ending in $\rho$. Denote by $V= V_m$ the space on which $\rho$ acts.
Then $\text{Res}^{H_m}_{H_{m-1}}V_m$ contains a unique subspace
 $V_{m-1}$ such that $\text{Res}^{H_m}_{H_{m-1}}\rho_m$ restricted to $V_{m-1}$ coincides with $\rho_{m-1}$. Iterating this
argument we obtain a sequence $V= V_{m}\supseteq V_{m_1}\supseteq \cdots \supseteq V_2\supseteq V_1$ such that
$\text{Res}_{H_j}^{H_m}\rho_m$ restricted to $V_j$ coincides with $\rho_j$ (more generally, for $m\geq k\geq j$,
$\text{Res}_{H_j}^{H_k}\rho_k$ restricted to $V_j$ is $\rho_j$).
As $\text{Res}_{H_1}^{H_m}\rho_m$ restricted to $V_1$ must coincide with $\rho_1 = \iota_K$, we have that $\dim V_1=1$ and $V_1$
is spanned by  a $K$-invariant vector in $V$. Let $\mathcal{C}(\rho)$ be the set of all paths ending in $\rho$.
For $C \in \mathcal{C}(\rho)$ we denote by $V_1(C)$ the corresponding one dimensional  space  as above, and we select
a vector $v_C \in V_1(C)$ such that $\|v_C\| =1$ (such a vector is defined modulo a scalar multiple $\alpha\in \CC$ with $|\alpha| = 1$). By induction on $m$, it is easy to prove that the set $\{v_C: C \in \mathcal{C}(\rho)\}$ is an orthonormal basis for $V^K$.

\begin{definition}{\rm The set $\{v_C: C \in \mathcal{C}(\rho)\}$ is called {\it the Gelfand-Tsetlin basis} for the $K$-invariant vectors in $(\rho,V)$, associated to the chain $G=H_m\geq H_{m-1}\geq H_{m-2}\geq\dotsb\geq H_1=K$.}
\end{definition}

Suppose again that $C : \rho_1, \rho_2 \ldots, \rho_m \in
\mathcal{C}(\rho)$ and denote by $V_1 \equiv \CC$ the space on which
acts $\rho_1 = \iota_K$.  Then $L(H_2/H_1) = \text{Ind}_{H_1}^{H_2}V_1$
contains a unique irreducible representation $V_2$  isomorphic to $\rho_2$;
iterating, $\text{Ind}_{H_{m-1}}^{H_m}V_{m-1}$ contains a unique irreducible
representations $V_C$ isomorphic to $\rho_m$. Then

\begin{theorem}
\[
m_\rho V_\rho = \bigoplus_{C\in \mathcal{C}(\rho)}V_C
\]

is an orthogonal decomposition of the isotypic component $\rho$.  Moreover, the intertwiner operator associated to
$v_C$ by \eqref{proiezione} maps $V_\rho$ into $V_C$.
\end{theorem}
\begin{proof}
The only part to be proved is the statement about the image of the intertwining operator associate with $v_C$. But this follows by iterated use of Lemma \ref{traind}.
\end{proof}

 We now present an explicit example of Gelfand-Tsetlin basis in the setting of the representation theory of the symmetric group.

 Let $n$ be a positive integer. A {\it composition} of $n$ is an
ordered sequence  ${a} = (a_1,a_2,\ldots,a_h)$ of positive integers
such that $a_1 + a_2+\cdots + a_h = n$.
A {\it partition} of $n$ is a composition $\lambda =
(\lambda_1,\lambda_2\ldots,\lambda_h)$ such that $\lambda_1 \geq
\lambda_2 \geq\cdots\geq \lambda_h$.

Set  $I(n) = \{1,2,\ldots,n\}$  and let   $a$ be a composition of  $n$. A {\it composition} (or {\it ordered partition}) of $I(n)$
of type $a$  is an ordered sequence
$A = (A_1, A_2, \ldots, A_h)$ of subsets  $A_i \subseteq I(n)$, called the {\it $A-$elements}, such that
$|A_i| = a_i, \mbox{ and } A_i\cap A_j = \emptyset$ for   all $1 \leq i \neq j \leq n$.
Note that under such conditions  $\coprod_{i=1}^{n}A_i = I(n)$. In other words, the $A-$elements partition the set $I(n)$.
Denote by $\Omega_a$ the set of all compositions of $I(n)$ of type $a$.

For $A \in \Omega_a$, denote by $S_A = \{\sigma\in S_n: \sigma(A_i) \subseteq A_i, \ \forall i = 1,2, \dots, h\}$
the {\it stabilizer} of $A$  and by $S_{A_i} = \{\sigma\in S_n:\sigma(j) = j, \ \forall j \in I(n)\setminus A_i\}$
its subgroup  which acts non-trivially only on $A_i$: $S_{A_i}$ is clearly isomorphic to $S_{a_i}$, the symmetric group
on $a_i$ elements. This way, $S_A = S_{A_1}\times S_{A_2}\times \cdots \times S_{A_h} \cong
S_{a_1}\times S_{a_2}\times \cdots \times S_{a_h}$, for all $A \in \Omega_a$.

Fix,  once  for all, an element $A^* = (A_1^*, A_2^*, \ldots, A_h^*) \in \Omega_a$  and denote, by abuse of language,
its stabilizer  $S_{A^*}$ by $S_a = S_{a_1}\times S_{a_2}\times \cdots \times S_{a_h}$, which therefore will be
identified with a specific subgroup of $S_n$. Noting that $S_n$ acts transitively on $\Omega_a$, we  can  regard $\Omega_a$ as the homogeneous
space
$\Omega_a \equiv S_n/(S_{a_1}\times S_{a_2}, \times \cdots \times S_{a_h})$; in particular
$|\Omega_a| = \frac{n!}{a_1!a_2!\cdots a_h!}$.

Let $a$ and $b$ be two positive integers and consider $\lambda = (\lambda_1, \lambda_2, \ldots, \lambda_k)  \vdash a+b$
 and $\mu = (\mu_1, \mu_2, \ldots, \mu_h) \vdash a$ two partitions
of $a+b$ and $a$, respectively.   We write $\lambda \succeq \mu$ to
denote that  $h \leq k$ and
\begin{equation}
\lambda_1 \geq \mu_1 \geq \lambda _2 \geq \mu_2 \geq  \cdots \geq  \lambda_h \geq \mu_h.
 \end{equation}

\begin{proposition}\label{littlewood}
The multiplicity  of $S^\mu\otimes S^{(b)}$ in $Res^{S_{a+b}}_{S_a \times S_b}S^\lambda$ ($\equiv$ the multiplicity
of $S^\lambda$ in $Ind_{S_a\times S_b}^{S_{a+b}}(S^\mu \otimes S^{(b)})$ is equal to
\[
\left\{\begin{array}{ll}1  & \mbox{ if $\lambda \succeq \mu$}\\
 0 &\mbox{ otherwise.}
\end{array}
\right.
\]
\end{proposition}
\begin{proof}
It is a particular case of the {\it Littlewood - Richardson rule} (see \cite{James,sagan}).
\end{proof}

In virtue of Proposition \ref{littlewood}, if $(a_1, a_2, \ldots, a_m)$ is a {\it composition} of $n$
then

\begin{equation}\label{chain}
\begin{array}{lcl}
H_m & =& S_n = S_{a_1+a_2+\cdots+ a_m}\\
H_{m-1}& =& S_{a_1+a_2+\cdots+ a_{m-1}}\times S_{a_m}\\
\ \ \vdots & \vdots & \ \ \ \ \ \ \ \  \vdots\\
H_2 &= & S_{a_1+a_2} \times S_{a_3} \times \cdots \times  S_{a_m}\\
H_1 &= & S_{a_1} \times S_{a_2} \times \cdots \times  S_{a_m}
\end{array}
\end{equation}

satisfies the Gelfand-Tsetlin condition.

From the representation theory of $S_n$ \cite{James,sagan}, we know that if $\mu = (\mu_1, \mu_2, \ldots, \mu_k) \vdash a_1 + a_2+
\cdots + a_j$  then $S^\mu$ contains nontrivial $S_{a_1}\times S_{a_2} \times \cdots S_{a_j}$-invariant vectors (and this is
equivalent to say that
 $S^\mu$ is  contained in  $M^{a_1, a_2, \ldots, a_j}$) if and only if  $\mu \trianglerighteq (a_1, a_2, \ldots, a_j)$,
that is $\mu_1 \geq \max\{a_i\}$, $\mu_1+\mu_2\geq \max\{a_i +a_j\}$, $\ldots$, $\mu_1+ \mu_2+ \cdots+ \mu_{k-1}
\geq  \max\{a_{i_1}+ a_{i_2} + \cdots + a_{i_k}\}$ (Young's rule).

Then the levels of the Brattelli diagram of $S_{a_1+\cdots+ a_m}\geq \cdots \geq S_{a_1}\times \dots \times S_{a_m}$
are
\[
\mathcal{B}_j = \{\mu \vdash (a_1+a_2+ \cdots +a_j) : \mu \trianglerighteq (a_1, a_2, \ldots, a_j)\}
\]
and $\mu \in \mathcal{B}_j$ is connected to $\nu \in \mathcal{B}_{j+1}$ when $\nu \succeq \mu$, that is when
$S^\mu\otimes S^{(a_{j+1})}$ is contained in
$Res_{S_{a_1+a_2+ \cdots +a_j}\times S_{a_{j+1}}}^{S_{a_1+a_2+ \cdots +a_{j+1}}}S^\nu$.

%%%%%%%%%%%%%%%%%%%%%%%%%%%%%%%%%%%%%%%%%%%%%%%%%%%%%%%%%%%%%%%%%%%%%%%%%%%
%%%%%%%%%%%%%%%%%%%%%%%%%%%%%%%%%%%%%%%%%%%%%%%%%%%%%%%%%%%%%%%%%%%%%%%%%%%
\section{Explicit computations on the symmetric group}
%%%%%%%%%%%%%%%%%%%%%%%%%%%%%%%%%%%%%%%%%%%%%%%%%%%%%%%%%%%%%%%%%%%%%%%%%%%
%%%%%%%%%%%%%%%%%%%%%%%%%%%%%%%%%%%%%%%%%%%%%%%%%%%%%%%%%%%%%%%%%%%%%%%%%%%

%%%%%%%%%%%%%%%%%%%%%%%%%%%%%%%%%%%%%%%%%%%%%%%%%%%%%%%%%%%%%%%%%%%%%%%%%%%
%%%%%%%%%%%%%%%%%%%%%%%%%%%%%%%%%%%%%%%%%%%%%%%%%%%%%%%%%%%%%%%%%%%%%%%%%%%
\subsection{Gelfand-Tsetlin bases in $S^{n-1,1}$}
%%%%%%%%%%%%%%%%%%%%%%%%%%%%%%%%%%%%%%%%%%%%%%%%%%%%%%%%%%%%%%%%%%%%%%%%%%%
%%%%%%%%%%%%%%%%%%%%%%%%%%%%%%%%%%%%%%%%%%%%%%%%%%%%%%%%%%%%%%%%%%%%%%%%%%%

 We now compute  a Gelfand-Tsetlin basis for the representation $S^{n-1,1}$.

 \begin{theorem}\label{vettori}Let $(a_1,a_2, \ldots, a_m)$ be a composition of $n$.
 Set $A_1 = \{1,2,\ldots, a_1\}$,
 $A_2 = \{a_1+ 1, a_1+2, \ldots,a_1+ a_2\}$, $\ldots$, $A_m  = \{a_{m-1}+ 1, a_{m-1}+2, \ldots,a_{m-1}+ a_m\}$. Then the Gelfand-Tsetlin basis of $S^{n-1,1}$, with respect  the chain of subgroups in (\ref{chain}), is given by

\[
v_j = \sqrt{\frac{a_j}{(a_1+a_2+ \cdots +a_{j-1})(a_1+a_2+\cdots +a_j)}}{\bf{1}}_{A_1\cup A_2\cup \cdots A_{j-1}}-
\sqrt{\frac{a_1+a_2+ \cdots +a_{j-1}}{a_j(a_1+a_2+\cdots +a_j)}}{\bf{1}}_{A_{j}}
\]

for $j = m, m-1, \ldots ,2$.
\end{theorem}

\begin{proof}
First we make a preliminary computation.
For $a, b \in \NN$ set $A = \{1,2,\ldots ,a\}$ and $B = \{a+1, a+2, \ldots, a+b\}$. Then

\[
\begin{split}
M^{a+b-1,1} = &L(A\coprod B) \equiv L(A)\oplus L(B) = \left[M^{a-1,1}\otimes S^{(b)}\right]\oplus
\left[S^{(a)}\otimes M^{b-1,1}\right]\\
=&2\left( S^{(a)}\otimes S^{(b)} \right)\oplus\left(S^{a-1,1}\otimes S^{(b)}\right) \oplus \left(S^{(a)}\otimes S^{b-1,1}\right)
\end{split}
\]

Recalling that $S^{a+b-1,1}= \{f \in L(A\cup B): \sum_{j \in A\cup B}f(j) = 0\}$, we deduce that

\begin{equation}\label{resdec}
\text{Res}_{S_a\times S_b}^{S_{a+b}}S^{a+b-1,1} = \left(S^{(a)}\otimes S^{(b)}\right)\oplus
\left(S^{a-1,1}\otimes S^{(b)}\right) \oplus \left(S^{(a)}\otimes S^{b-1,1}\right)
\end{equation}

where

\[
\begin{split}
S^{(a)}\otimes S^{(b)} \equiv \left\langle\frac{1}{a}{\bf{1}}_A-\frac{1}{b}{\bf{1}}_B\right\rangle, \ \ \ \ \ \ \ \ \ \ \\
S^{a-1,1}\otimes S^{(b)} = \{f \in L(A\cup B): \sum_{i \in A}f(i) = 0,f|_B\equiv 0\},\\
S^{(a)} \otimes S^{b-1,1} = \{f \in L(A\cup B): \sum_{i \in B}f(i) = 0,f|_A\equiv 0\}.
\end{split}
\]

Then
simple computations show that a normalized  $S_a \times S_b$-invariant vector in $S^{(a)}\otimes S^{(b)}$ is

\begin{equation}\label{invvec}
\sqrt{\frac{b}{a(a+b)}}{\bf{1}}_A
-\sqrt{\frac{a}{b(a+b)}}{\bf{1}}_B
\end{equation}

and, clearly,  the other $S_a \times S_b$-invariant vectors in $S^{(a)}\otimes S^{(b)}$ are multiple
of such a vector.\\

Now we show how to obtain the Gelfand-Tseitlin basis for $S^{n-1,1}$
associated to the chain \eqref{chain}.
First of all, by \eqref{resdec} we can write

\begin{multline}\label{resiter}
\text{Res}_{S_{a_1+ a_2 + \cdots a_{m-1}}\times S_{a_m}}^{S_{a_1+ a_2 + \cdots a_{m}}}S^{n-1,1} = (S^{(a_1+a_2+ \cdots+ a_{m-1})}\otimes S^{(a_m)})\\ \oplus
 (S^{a_1+a_2+ \cdots a_{m-1}-1,1} \otimes S^{(a_m)})\oplus (S^{(a_1+a_2+ \cdots a_{m-1})} \otimes S^{a_m-1,1}).
\end{multline}

The first representation (that is the trivial representations) obviously contains a nontrivial
$S_{a_1}\times S_{a_2}\times \cdots S_{a_m}$-invariant vector, given by \eqref{invvec}

\begin{multline}
v_m = \sqrt{\frac{a_m}{(a_1+a_2+ \cdots +a_{m-1})(a_1+a_2+\cdots +a_m)}}{\bf{1}}_{A_1\cup A_2\cup \cdots A_{m-1}}\\
-\sqrt{\frac{a_1+a_2+ \cdots +a_{m-1}}{a_m(a_1+a_2+\cdots +a_m)}}{\bf{1}}_{A_{m}}.
\end{multline}

On the other hand, it is clear that the third representation in \eqref{resiter}
does not contain any nontrivial such a vector, as it is $S_{a_m}$-irreducible. Therefore we analyze the second block or, more precisely, its
restriction to $S_{a_1+ a_2 + \cdots a_{m-2}}\times S_{a_{m-1}}\times S_{a_m}$, that is

\begin{multline*}
\text{Res}_{S_{a_1+ a_2 + \cdots a_{m-2}}\times S_{a_{m-1}}\times S_{a_m}}^{S_{a_1+ a_2 + \cdots a_{m-1}}\times S_{a_m}}S^{a_1+a_2+ \cdots a_{m-1}-1,1}\otimes S^{(a_m)}  = (S^{(a_1+a_2+ \cdots+ a_{m-2})} \otimes S^{(a_{m-1})} \otimes S^{(a_m)})
\oplus \\
\oplus (S^{a_1+a_2+ \cdots+ a_{m-2}-1,1} \otimes S^{(a_{m-1})} \otimes S^{(a_{m})}) \oplus
(S^{(a_1+a_2+ \cdots+ a_{m-2})} \otimes S^{a_{m-1}-1,1} \otimes S^{(a_{m})}).
\end{multline*}

Again the first block is the trivial representation and therefore contains an invariant vector, namely
\begin{multline}
v_{m-1} =   \sqrt{\frac{a_{m-1}}{(a_1+a_2+ \cdots +a_{m-2})(a_1+a_2+\cdots +a_{m-1})}}{\bf{1}}_{A_1\cup A_2\cup \cdots A_{m-2}} \\ -
\sqrt{\frac{a_1+a_2+ \cdots +a_{m-2}}{a_{m-1}(a_1+a_2+\cdots +a_{m-1})}}{\bf{1}}_{A_{m-1}},
\end{multline}
the last representation does not contain any nontrivial invariant
vector and the intermediate term needs to be analyzed. Iterating
this argument, we can end the proof.
\end{proof}

We now write the spherical functions corresponding to the vectors $v_j$ (see Definition \ref{sphericalfunctions}) and we deduce in some cases an explicit
expression of the projector $E_i^\rho$ defined in (\ref{proj1}).

Let $(B_1, B_2, \ldots, B_m)\in\Omega_a$ and choose $\sigma \in S_n$ such that $\sigma (A_i) = B_i$ for $i = 1,2, \ldots,m$, i.e. $\sigma$ is a representative
of the coset of $S_n/(S_{a_1} \times S_{a_2}\times \cdots \times S_{a_m})$ identified by $(B_1,B_2, \ldots, B_m)$.
Set ${\bf{A}}^j = A_1\cup A_2\cup \cdots A_{j-1}$, ${\bf{B}}^j = B_1\cup B_2\cup \cdots B_{j-1}$,
${\bf{a}}_j = \sqrt{\frac{a_j}{(a_1+a_2+ \cdots +a_{j-1})(a_1+a_2+\cdots +a_j)}}$,
${\bf{b}}_j =\sqrt{\frac{a_1+a_2+ \cdots +a_{j-1}}{a_j(a_1+a_2+\cdots +a_j)}}$  and $h(j) = |{\bf{A}}^j| = |{\bf{B}}^j|\equiv a_1+a_2+\dotsb a_{j-1}$.

\begin{proposition}\label{p;projection}
An explicit expression for the spherical functions $\phi_j=\langle v_j,\sigma v_j\rangle$, seen as functions of $(B_1,B_2,\dotsc,B_m)$ (cf. \eqref{e;Kinv}), is given by

\begin{multline}
\phi_j = \sum_{s = \max\{0, 2h(j)-n\}}^{h(j)} \sum_{u = \max \{0, a_j+2h(j)-s-n\}}^{\min\{a_j, h(j)-s\}}
\sum_{v = \max \{0, a_j+2h(j)-s-n\}}^{\min\{a_j, h(j)-s\}}\\
\sum_{z = \max \{0, 2a_j-u-v-n+2h(j)-s\}}^{\min\{a_j-u, a_j-v\}}
\frac{sa_j^2-a_jh(j)(u+v)+zh(j)^2}{a_jh(j)h(j+1)}  {\bf{1}}_{\mathcal{C}_{s,u,v,z}}
\end{multline}

where

\[\mathcal{C}_{s,u,v,z} = \{(B_1,B_2,\ldots,B_m):  |{\bf{A}}^j\cap {\bf{B}}^j| = s, \ |{\bf{A}}^j\cap B_j| = u, \
|{\bf{B}}^j\cap A_j| = v \mbox{ and } |A_j \cap B_j| = z\}.
\]

\end{proposition}

\begin{proof}
The proposition follows from the explicit expression of $v_j$ in Theorem \ref{vettori},
applying the following observations ($C$ and $D$ are two finite subsets of a finite set $F$)

\begin{itemize}
\item{$\alpha{\bf{1}}_C \cdot \beta{\bf{1}}_D = \alpha \beta {\bf{1}}_{C\cap D}$},
\item $({\bf a}_j)^2s-{\bf a}_j{\bf b}_j(u+v)+({\bf b}_j)^2z=\frac{sa_j^2-a_jh(j)(u+v)+zh(j)^2}{a_jh(j)h(j+1)}$,
\item{$\max\{0, -|F|+|C|+|D|\} \leq |C \cap D| \leq \min \{|C|,|D|\}$,}
\item ${\bf A}^j\cap B_j=\left[{\bf A}^j\setminus ({\bf A}^j\cap{\bf B}^j)\right]\cap B_j\subseteq I(n)\setminus {\bf B}^j$,
\item ${\bf B}^j\cap A_j=\left[{\bf B}^j\setminus ({\bf B}^j\cap{\bf A}^j)\right]\cap A_j\subseteq I(n)\setminus {\bf A}^j$,
\item $A_j\cap B_j=\left[A_j\setminus (A_j\cap{\bf B}^j)\right]\cap\left[B_j\setminus(B_j\cap{\bf A}^j)\right]\subseteq I(n)\setminus ({\bf A}^j\cup {\bf B}^j)$.
\end{itemize}

\end{proof}

Proposition \ref{p;algebraf11}, writing $E^{1}_i$  instead of  $E^{S^{n-1,1}}_i$, becomes

\begin{equation}\label{ppp}
\begin{split}
(E^{1}_i f)(B_1, B_2, \ldots, B_m)= &(E^{1}_i f)(\sigma(A_1), \sigma(A_2), \ldots, \sigma(A_m)) =\frac{(n-1)a_1!a_2! \cdots  a_m!}{n!}\langle {f}, \sigma \phi_i\rangle_{L(X)}\\
=&\frac{(n-1)a_1!a_2! \cdots  a_m!}{n!}\sum_{s,u,v,z}\left\langle f,{\bf{1}}_{\sigma\mathcal{C}_{s,u,v,z}}
\right\rangle_{L(X)}.
\end{split}
\end{equation}

Note that $(B_1, B_2, \ldots, B_m)$ is determined by $(B_2, \ldots , B_m)$, as $B_1 = (B_2\cup \cdots \cup B_m)^c$.

\begin{example}\label{proexa}{\rm
 We  now apply the above  formulas to the particular case

\[
H_1 = S_{n-2}\times S_1 \times S_1\leq H_2 = S_{n-1}\times S_1\leq    H_3 = S_n.
\]

Therefore, $m= 3$, $a_1 = n-2$
and  $a_2 = a_3 = 1$. We identify $S_n/(S_{n-2}\times S_1\times S-1)$ with $\{(i,j):i,j\in I(n),i\neq j\}$. Then $s=\lvert\{1,2,\dotsc,n-2\}\cap\{i,j\}^C\rvert$, $u=\lvert\{1,2,\dotsc,n-2\}\cap\{i\}\rvert$, $v=\lvert\{i,j\}^C\cap\{n-1\}\rvert$ and $z=\lvert\{i\}\cap\{n-1\}\rvert$.

We determine the spherical function associated with
\[
v_2 = \sqrt{\frac{1}{(n-1)(n-2)}}{\bf{1}}_{\{1,2, \ldots, n-2\}}-\sqrt{\frac{n-2}{n-1}}{\bf{1}}_{\{n-1\}}.
\]

Clearly, $n-4 \leq s \leq n-2$, $h(2)=n-2$ and $h(3)=n-1$.
\begin{itemize}
\item{If $s = n-2$ then $u = v= 0$ and $0\leq z \leq 1$, and we have two possible cases:
\begin{itemize}
\item{$\mathcal{C}_{n-2,0,0,1} = \{(n-1,n)\}$}
\item{$\mathcal{C}_{n-2,0,0,0} = \{(n,n-1)\}$}
\end{itemize}}
\item{If $s = n-3$ then $0\leq u\leq 1$, $0\leq v \leq 1$, $ \max\{0,1-u-v\} \leq z \leq \min\{1-u, 1-v\}$,
and we have 4 possible cases
\begin{itemize}
\item{$\mathcal{C}_{n-3,1,0,0} = \{(i,n-1): i \in \{1,2, , \ldots,  n-2\}\}$}
\item{$\mathcal{C}_{n-3,0,0,1} = \{(n-1,j): j \in  \{1,2, , \ldots,  n-2\}\}$}
\item{$\mathcal{C}_{n-3,0,1,0} = \{(n,j): j \in \{1,2, , \ldots,  n-2\}\}$}
\item{$\mathcal{C}_{n-3,1,1,0} = \{(i,n): i \in \{1,2, \ldots, n-2\}\}$}.
\end{itemize}}
\item{For $s = n-4$ we have $u= v = 1$ and $z = 0$ and the only possible case is
\begin{itemize}
\item{$\mathcal{C}_{n-4,1,1,0} = \{(i,j): i,j \in \{1,2, , \ldots,  n-2\}\}$}
\end{itemize}}
\end{itemize}
Then by Proposition \ref{p;projection} we have
\begin{multline}
\phi_2 = {\bf{1}}_{\mathcal{C}_{n-2,0,0,1}} +\frac{1}{n-1}{\bf{1}}_{\mathcal{C}_{n-2,0,0,0}}
-\frac{1}{(n-1)(n-2)}{\bf{1}}_{\mathcal{C}_{n-3,1,0,0}} +\frac{n^2-3n+1}{(n-1)(n-2)}{\bf{1}}_{\mathcal{C}_{n-3,0,0,1}}\\
 -\frac{1}{(n-1)(n-2)}{\bf{1}}_{\mathcal{C}_{n-3,0,1,0}}
-\frac{1}{n-2}{\bf{1}}_{\mathcal{C}_{n-3,1,1,0}}
-\frac{n}{(n-1)(n-2)}{\bf{1}}_{\mathcal{C}_{n-4,1,1,0}},
\end{multline}
and therefore, by (\ref{ppp}) the corresponding projection is given by
\begin{multline}
E^1_2f(i,j) =\frac{1}{n}f(i,j)+\frac{1}{n(n-1)}f(j,i) -\frac{1}{n(n-1)(n-2)}\sum_{h\neq i,j}f(h,i)
+\frac{n^2-3n+1}{n(n-1)(n-2)}\sum_{k\neq i,j}f(i,k)
\\-\frac{1}{n(n-1)(n-2)}\sum_{k \neq i,j} f(j,k)-
\frac{1}{n(n-2)}\sum_{h\neq i,j}f(h,j)
  -\frac{1}{(n-1)(n-2)}\sum_{\substack{h\neq k\\h,k\neq i,j}}f(h,k).
\end{multline}

Simple computations show that

\begin{equation}\label{E12f}
(E^1_2 f)(i,j) = \frac{n-1}{n-2}\overline{f}(i)+\frac{1}{n-2}\overline{f}(j),
\end{equation}

where $$\overline{f}(t) = \frac{1}{n}\sum_{k\neq t}f(t,k) + \frac{1}{n(n-1)}\sum_{h \neq t}f(h,t)
-\frac{1}{n(n-1)}\sum_{h \neq k}f(h,k).$$
}\end{example}

%%%%%%%%%%%%%%%%%%%%%%%%%%%%%%%%%%%%%%%%%%%%%%%%%%%%%%
%%%%%%%%%%%%%%%%%%%%%%%%%%%%%%%%%%%%%%%%%%%%%%%%%%
\subsection{Two decompositions of $M^{n-2,1,1}$}\label{decMn211}
%%%%%%%%%%%%%%%%%%%%%%%%%%%%%%%%%%%%%%%%%%%%%%%%%%%%%%
%%%%%%%%%%%%%%%%%%%%%%%%%%%%%%%%%%%%%%%%%%%%%%%%%%
In this section we present two different decompositions of the module $M^{n-2,1,1}$. Since

\begin{equation}\label{Mn211}
M^{n-2,1,1} =
 L\left(S_n/(S_{n-2} \times S_1\times S_1)\right) =
 S^{(n)}\oplus 2S^{n-1,1} \oplus S^{n-2,2}\oplus S^{n-2,1,1},
\end{equation}

this means that we give two different ways to decompose $2S^{n-1,1}$. \\

The first decomposition is based on the chain of subgroups  $S_n \geq S_{n-2} \times S_{2} \geq S_{n-2}\times S_1 \times S_1$.
Using the notation of subsection \ref{trainddec}, we have $L(X)= M^{n-2,1,1}$, $L(Y) = M^{n-2,2}$ and $Z = \Omega_{1^2} =  S_2$.
Clearly $L(Z) = M^{1,1} = S^{(2)}\oplus S^{1,1}$ and therefore

\[
\begin{split}
M^{n-2,1,1} & = \text{Ind}_{S_{n-2}\times S_2}^{S_n}\left[S^{(n-2)} \otimes M^{1,1}\right]\\
& =\text{Ind}_{S_{n-2}\times S_2}^{S_n}\left[S^{(n-2)} \otimes S^{(2)}\right] \oplus
 \text{Ind}_{S_{n-2}\times S_2}^{S_n}\left[S^{(n-2)} \otimes S^{1,1}\right]  \\
& = M_S^{n-2,1,1} \oplus M_A^{n-2,1,1}.
\end{split}
\]

where

\[
M_S^{n-2,1,1} = \{f \in M^{n-2,1,1}: f(i,j) = f(j,i), \ \forall i,j \in \{1,2, \ldots n\}, \ i \neq j\}
\]

and

\[
M_A^{n-2,1,1} = \{f \in M^{n-2,1,1}: f(i,j) = -f(j,i), \ \forall i,j \in \{1,2, \ldots n\}, \ i \neq j\};
\]

this is also a particular case of Example \ref{exasymm}.
Clearly

\[
M_S^{n-2,1,1} \cong M^{n-2,2} = S^{(n)} \oplus S^{n-1,1}_S \oplus S^{n-2,2}.
\]

and therefore from \eqref{Mn211} we deduce that

\[
M_A^{n-2,1,1} = S^{n-1,1}_A \oplus S^{n-2,2}.
\]

The projections onto $M_S^{n-2,1,1}$ and $M_A^{n-2,1,1}$ are given by
\[
\begin{array}{lccc}
P_S :& M^{n-2,1,1} & \to & M^{n-2,1,1}_S\\
& f & \mapsto & f_S
\end{array}
\qquad\text{where} \quad f_S(i,j) = \frac{f(i,j)+f(j,i)}{2}
\]

and

\[
\begin{array}{lccc}
P_A :& M^{n-2,1,1} & \to & M^{n-2,1,1}_A\\
& f & \mapsto & f_A
\end{array}
\qquad\text{where}\quad f_A(i,j) = \frac{f(i,j)-f(j,i)}{2}.
\]

We now determine the projections onto each irreducible component.
If we apply Theorem \ref{vettori} with $a_1=a_2=1$, $a_3=n-2$ and $A_1=\{n\},A_2=\{n-1\}, A_3=\{1,2,\dotsc,n-2\}$, we get an orthonormal basis for the
$S_{n-2} \times S_1 \times S_1$-invariant vectors in the representation $S^{n-1,1}$, given by

\[
 v_S \equiv v_3= \sqrt{\frac{n-2}{2n}}{\bf{1}}_{\{n-1,n\}}-\sqrt{\frac{2}{n(n-2)}}{\bf{1}}_{\{1,2,\ldots, n-2\}}
\]

and

\[
v_A \equiv v_2= \sqrt{\frac{1}{2}}{\bf{1}}_{\{n\}}-\sqrt{\frac{1}{2}}{\bf{1}}_{\{n-1\}}.
\]

The corresponding matrix coefficients are
\[
\phi_S = {\bf{1}}_{\mathcal{C}_{n-2,0,0,2}} +\frac{n-4}{2(n-2)}{\bf{1}}_{\mathcal{C}_{n-3,1,1,1}}
-\frac{2}{n-2}{\bf{1}}_{\mathcal{C}_{n-4,2,2,0}}
\]
where
\begin{itemize}
\item{$\mathcal{C}_{n-2,0,0,2} = \{(i,j): i,j \in \{n-1,n\}, \ i \neq j\}$}
\item{$\mathcal{C}_{n-3,1,1,1} = \{(i,j):|\{i,j\} \cap \{n-1,n\}| = 1, i\neq j\}$}
\item{$\mathcal{C}_{n-4,2,2,0} = \{(i,j): \{i,j\}\cap\{n-1,n\}=\emptyset,i\neq j\}$}
\end{itemize}
and
\begin{multline}
\phi_A = {\bf{1}}_{\mathcal{C}_{1,0,0,1}} +\frac{1}{2}{\bf{1}}_{\mathcal{C}_{1,0,0,0}} + \frac{1}{2}{\bf{1}}_{\mathcal{C}_{0,0,0,1}}
-\frac{1}{2}{\bf{1}}_{\mathcal{C}_{0,1,0,0}} -\frac{1}{2}{\bf{1}}_{\mathcal{C}_{0,0,1,0}} - {\bf{1}}_{\mathcal{C}_{0,1,1,0}}
\end{multline}
where
\begin{itemize}
\item{$\mathcal{C}_{1,0,0,1} = \{(n-1,n)\}$}
\item{$\mathcal{C}_{1,0,0,0} = \{(n-1,j):j \neq n-1,n\}$}
\item{$\mathcal{C}_{0,0,0,1} = \{(i,n):i \neq n-1,n\}$}
\item{$\mathcal{C}_{0,1,0,0} = \{(i,n-1):i \neq n-1,n\}$}
\item{$\mathcal{C}_{0,0,1,0} = \{(n,j):j \neq n-1,n\}$}
\item{$\mathcal{C}_{0,1,1,0} = \{(n,n-1)\}$.}
\end{itemize}

The corresponding projectors are given by
\begin{multline}
E^1_Sf(i,j) = \frac{1}{n}\left[f(i,j)+ f(j,i) + \frac{n-4}{2(n-2)}\left(\sum_{\substack{k\neq i,j}}f(i,k)+
\sum_{\substack{k\neq i,j}}f(j,k) \right.\right.
\\
\left.\left. + \sum_{\substack{h\neq i,j}}f(h,i)+\sum_{\substack{h\neq i,j}}f(h,j)\right)-\frac{2}{n-2}
\sum_{\substack{k\neq h\\h,k\neq i,j}}f(h,k)\right]
\end{multline}
and
\begin{multline}
E^1_Af(i,j) = \frac{1}{n}\left[f(i,j)+ \frac{1}{2}\sum_{\substack{k\neq i,j}}f(i,k)+
\frac{1}{2}\sum_{\substack{h\neq i,j}}f(h,j)\right.
\\
 \left. -\frac{1}{2} \sum_{\substack{h\neq i,j}}f(h,i)-\frac{1}{2}\sum_{\substack{k\neq i,j}}f(j,k)
-f(j,i)\right].
\end{multline}

We can summarize the above analysis by saying that
a function $f \in M^{n-2,1,1}$ can be decomposed as
$f = f_S ( = P_S f \in M^{n-2,1,1}_S)  + f_A (= P_A f \in M^{n-2,1,1}_A)$ and these, in turn,  can be decomposed as
\[
f_S = f_{S,0}(i,j) ( \in S^{(n)}) + f_{S,1}(\in S_S^{n-1,1}) + f_{S,2}(\in S^{n-2,2})
\]
 and
 \[
 f_A = + f_{A,1}(\in S_A^{n-1,1}) + f_{1,1}(\in S^{n-2,1,1}).
\]

Note that the other projectors can be easily computed as difference: for instance the projector onto $S^{n-2,1,1}$
is equal to $P_A- E^1_A$ (and the projector onto $S^{(n)}$ is just the average).

We now present a statistical interpretation of this decomposition.
Suppose that we have an election for the president and the director
of an association with the following rule: every elector chooses a
pair $(i,j)$ where $i$ is its favorite  candidate  as president and
$j$ is its favorite candidate as director. Then the data of the
election is an element $f \in M^{n-2,1,1}$ and the
projections have the following interpretations: (compare with \cite{Diaconis,DiaconisW})
\begin{itemize}
\item{$f_S(i,j)$ represents the vote for the unordered pair  $\{i,j\}$: there is not preference for $i$ as president or
$j$ as director and viceversa;}
\item{$f_{S,0}(i,j)$ is just the average vote;}
\item{$f_{S,1}(i,j)$ is the effect of the popularity of the  single inside   an unordered pair;}
\item{$f_{S,2}(i,j)$ is the vote of the unordered pair without the effect of the singles;}
\item{$f_A(i,j)$ is the vote to the ordered pair without the effect of the  unordered pairs;}
\item{$f_{A,1}$ is the effect of the single in the ordered pair without the effect of the  unordered pairs;}
\item{$f_{1,1}$ is the vote to the ordered pair without all the other effects.}
\end{itemize}

We  now describe  an alternative decomposition of $M^{n-2,1,1}$. The only difference with respect the previous case it
is the different decomposition of the isotypic component $2S^{n-1,1}$.  This time we  use the chain of subgroups in Example \ref{proexa}.
In the notations of Example \ref{exasymm}, we have the following positions:
$H= S_{n-1}\times S_1$, $K= S_{n-2}$, $X = S_n/S_{n-2} = \{(i,j): i,j \in \{1,2,\ldots,n\}\}$, $Y = S_n/S_{n-1} = \{1,2,\ldots, n\}$ and $Z = S_{n-1}/S_{n-2} = \{1,2,\ldots, n-1\}$.
The map $\pi: X\to Y$ is given by $(i,j) \mapsto j$ and therefore $\pi^{-1}(j) = \{(i,j): i \in  \{1,2,\ldots,n\}\setminus\{j\}$.
Observe that $L(\pi^{-1}(j)) \cong M^{n-2,1} = S^{(n-1)}(j)\oplus S^{n-2,1}(j)$
where $S^{(n-1)}(j)$ (resp. $S^{n-2,1}(j)$) denotes the space of constant  functions (resp. of mean value zero)  on $\pi^{-1}(j)$.

Therefore, applying the theory developed in subsection \ref{trainddec}, one gets
\[
\begin{split}
M^{n-2,1,1}
& =  \bigoplus_{j = 1}^n L(\pi^{-1}(j))  \\
& = \left(\bigoplus_{j = 1}^nS^{(n-1)}(j)\right) \bigoplus \left(\bigoplus_{j = 1}^nS^{n-2,1}(j)\right).
\end{split}
\]

In other words

\[
M^{n-2,1,1}_1 = \text{Ind}_{S_{n-1}\times S_1}^{S_n}\left[S^{(n-1)}\otimes S^{(1)}\right] =  \bigoplus_{j = 1}^nS^{(n-1)}(j)
\]

represents the subspace of $M^{n-2,1,1}$ of all functions that depend only on the $j$ coordinate,
while

\[
\begin{split}
M^{n-2,1,1}_2 &  = \text{Ind}_{S_{n-1}\times S_1}^{S_n}\left[S^{n-2,1}\otimes S^{(1)}\right] =
 \bigoplus_{j = i}^nS^{n-2,1}(j) \\
 & = \{f\in M^{n-2,1,1}:\sum_{\substack{i \\i \neq j}}f(i,j) = 0 , \ \forall j \in \{1,2,\ldots,n\}\}.
 \end{split}
\]

Clearly $M^{n-2,1,1}_1 \cong M^{n-1,1} = S^{(n)}\oplus S^{n-1,1}_1$ and therefore
$M^{n-2,1,1}_2=  S^{n-1,1}_2\oplus S^{n-2,2}\oplus S^{n-2,1,1}$.
It remains to compute only the projection on $S^{n-1,1}_1$ (the projection onto $S^{n-1,1}_2$ has been computed in
Example \ref{proexa}).
It is easy to see that the projection $P_1$ onto $M^{n-2,1,1}_1$ is given, for  $f \in M^{n-2,1,1}$, by

\begin{equation}\label{P1f}
P_1f(i,j) = \frac{1}{n-1}\sum_{\substack{h \\h \neq j}}f(h,j).
\end{equation}

and therefore the projection $E^1_1$ onto $S^{n-1,1}_1$ is given by

\[
 E^1_1f(i,j)  =  \frac{1}{n-1}\sum_{\substack{h \\h \neq j}}f(h,j)- \frac{1}{n(n-1)}\sum_{\substack{h,k\\h \neq k}}f(h,k)
 \]

as the last term is nothing but the projection onto  $S^{(n)}$.

 In order to clarify the statistical meaning of this decomposition, we first observe that

\begin{equation}\label{Sn11}
S^{n-1,1}_1 = \{F \in M^{n-2,1,1}: F(i,j) = f(j), \sum_{h}f(h) = 0\},
\end{equation}

while from \eqref{E12f} one easily gets

\begin{equation}\label{Sn12}
S^{n-1,1}_2 = \{F \in M^{n-2,1,1}: F(i,j) = (n-1)f(i) + f(j), \sum_{h}f(h) = 0\}.
\end{equation}

Therefore the isotypic component is given by

\[
2S^{n-1,1} = \{F \in M^{n-2,1,1}: F(i,j) = f_1(i) + f_2(j), \sum_{h}f_1(h) = \sum_{h}f_2(h)= 0\}.
\]

The global projection onto the isotypic component is nothing but the
best approximation of the data  with a function of the form $f_1(i)
+ f_2(j)$, where $\sum_{h}f_1(h) = \sum_{h}f_2(h)= 0$. Therefore,
$f_1(i)$ is the effect that candidate $i$ is chosen as president and
$f_2(j)$ that $j$ is chosen as director (once the average vote has
been removed). On the other hand, $P_1f(i,j)$ is the average vote of
$j$ as director (see \eqref{P1f}). Therefore $E_1^1f$ equals $P_1f$,
once the average vote has been removed, while $E_2^1f$ is the best
approximation of the form $f_1(i)+f_j(j)$ (as above), once $E_1^1f$
has been subtracted out. For the interpretations of the projections
onto $S^{n-2,2}$ and $S^{n-2,1,1}$, we refer to
\cite{Diaconis,DiaconisW}.

 %%%%%%%%%%%%%%%%%%%%%%%%%%%%%%%%%%%%%%%%%%%%%%%%%%%%%%%%%%%%%%%%%%%%%%%%%%%%%%%%%%%%%%%%%%%%%%%%%%%
%%%%%%%%%%%%%%%%%%%%%%%%%%%%%%%%%%%%%%%%%%%%%%%%%%%%%%%%%%%%%%%%%%%%%
\subsection{The semistandard basis for the space $\text{Hom}_{S_n}(S^a, M^b)$}

The present subsection is a translation, into our framework, of the results in sections 2.9 and 2.10 of \cite{sagan}.  See also \cite{Sternberg}, Appendix C (and the original source \cite{James}).
Suppose that $a = (a_1, a_2,\ldots,  a_h)$ and $b = (b_1,b_2, \ldots, b_k)$ are  two compositions of $n$.
  We begin with a description of the orbits of $S_n$ on $\Omega_a \times \Omega_b$.

For $A \in \Omega_a$ and $B \in \Omega_b$  denote by $C = A\wedge B$ their {\it intersection}, namely
the composition of $I(n)$ whose elements are the non-empty  $A_i\cap B_j$'s ordered lexicographically (i.e.
$(i,j) \leq (i',j')$ if $i<i'$ or $i=i'$ and $j\leq j'$); denote by $c$ the corresponding type of $C$. The following
fact is obvious:

\begin{lemma}
Let $A\in \Omega_a$ and $B \in \Omega_b$ and  $C = A\wedge B \in \Omega_c$ and consider the   actions
of $S_n$ on $\Omega_a\times \Omega_b$ and on  $\Omega_c$. Then  $\sigma \in S_n$ fixes $(A,B) \in \Omega_a\times  \Omega_b$ if and only if
it fixes $C$. In other words, $S_A\cap S_B = S_{A\wedge B}$.
\end{lemma}

We replace the notion of intersection of two partition with the following notion:

 For $a \in C(n,h)$ and $b \in C(n,k)$ denote by
 $\mathfrak{M}_{a,b}$ the set  of all   matrices $(m_{i,j})\in M_{h\times k}(\NN)$ (with non-negative integer entries) such that
$\sum_{i = 1}^h m_{i,j} = b_j$ for all $j = 1,2, \ldots,k$ and
$\sum_{j = 1}^k m_{i,j} = a_i$ for all $i = 1,2, \ldots,h$.
We may also say that the column sums (resp. the row sums) of  the matrix $(m_{i,j})$ equal
the $b_j$'s (resp. $a_i$'s).

With $A \in \Omega_a$ and $B \in \Omega_b$, we associate the matrix $m = m(A,B) \in \mathfrak{M}_{a,b}$ defined by setting
$m_{ij} = |A_i\cap B_j|$, for all $j = 1,2, \ldots,k$ and $i = 1,2, \ldots,h$. We then have

\begin{lemma}\label{l;indsym}
$(A,B)$ and $(A',B') \in \Omega_a\times \Omega_b$
  belong to the same $S_n-$orbit  if and only if $m(A,B) = m(A',B')$. In particular
  the set of $S_n$ orbits on $\Omega_a \times \Omega_b$  is in one-to-one correspondence with $\mathfrak{M}_{a,b}$.
\end{lemma}
\begin{proof}
The ``only if'' part is obvious: $|A_i\cap B_j|  = |\sigma(A_i\cap B_j)|= |\sigma A_i\cap \sigma B_j|$ for all $\sigma \in S_n$. Conversely, suppose that $m(A,B) = m(A',B')$.  Then the intersection compositions $A\wedge B$ and $A'\wedge B'$
are of the same type, say $c$. As $S_n$ acts transitively on $\Omega_c$ there exists $\sigma \in S_n$ such that
$\sigma (A\wedge B) = A'\wedge B'$ and thus $\sigma (A, B) = (A',B')$.
\end{proof}

For $a \in C(n,h), \ b \in C(n,k)$ and $M \in \mathfrak{M}_{a,b}$ we define
$\mathcal{T}_M: M^a \to M^b$ by setting

\[
(\mathcal{T}_M f)(B) = \sum_{\substack{A \in \Omega_a:\\ m(A,B) = M}}f(A) \qquad\quad  \forall  B \in \Omega_b
\]

or equivalently

\[
\mathcal{T}_M  \delta_A = \sum_{\substack{B \in \Omega_b:\\ m(A,B) = M}}\delta_B \quad\qquad   \forall A \in \Omega_a.
\]

Then

\begin{theorem} $\{\mathcal{T}_M :M \in \mathfrak{M}_{a,b}\}$ is a basis for  $\text{\rm Hom}_{S_n}(M^a, M^b)$.
\end{theorem}
\begin{proof}(Sketch)
Theorem \ref{comm} may be generalized in the following way: if $X$ and $Y$ are $G$-space then $\text{Hom}_G(L(X),L(Y))$ is isomorphic to the vector space of all $G$-invariant functions on $X\times Y$ (see also \cite{ScarabottiRadon}). The correspondence

\[
\text{orbit associated to} M\longmapsto \mathcal{T}_M
\]

is just this isomorphism in the present setting.

\end{proof}

It is clear that the restriction of $\mathcal{T}_M $ to $S^a$ belongs to    $\text{Hom}_{S_n}(S^a, M^b)$. In \cite{sagan}, section 2.10, it is presented a basis of $\text{Hom}_{S_n}(S^a, M^b)$ in terms of
generalized tableaux. In the following we describe this basis in terms of the operators $\mathcal{T}_M$.
Keeping the notation of \cite{sagan} p.79, we first observe that the correspondence $\theta^{-1}$

\[
T \mapsto \{s\}
\]

could be expressed in the following way:
$\{s\}$ is the tabloid of shape $b$, obtained inserting $l$ into the row $T(l)$, $l = 1,2,\ldots,n$.
Therefore,
if $T$ is a generalized tableau of shape $a$ and content $b$,
and $\theta_T$ is the isomorphism associated to $T$ as in Definition 2.3.9 of \cite{sagan}, then
$\theta^{-1}\circ\theta_T: M^a \to M^b$, could be described by

\[
\theta^{-1}\circ\theta_T(\{t\}) = \sum_{s  \in \theta^{-1}(\{T\})} \{s\},
\]

where the sum is over all tabloids $\{s\}$ of shape $b$ such that the number of elements
in common between the $i$-th row of $\{s\}$
and the $j$-th row of $\{t\}$ is equal to the number of $i$ in the $j-$th row of $T$.
Indeed the elements  $S \in \{T\}$ are just the generalized tabloids obtained by permuting, in all possible ways, the rows
 of $T$ and therefore must go into the $S(l)$-th row. We deduce that the elements in the $j$-th rows of $\{t\}$  must go
in the rows of $\{s\} \in \theta^{-1}(\{T\})$ corresponding to the numbers in the $j$-th row of $T$.

In other words, if $T$ is a generalized tableau with the rows weakly
increasing (so that there is only one such a  tableau for each class
$\{T\}$) we obtain a  bijective correspondence
  $\theta^{-1}\circ \theta_T \mapsto  \mathcal{T}_M$ where  $M =\{m_{i,j}\}$ is the matrix given by

\[
m_{i,j} = \mbox{number of $i$'s in the $j$-th row of } T.
\]

We could also say that $\theta^{-1}\circ \theta_T$ is an
equivalent description of the operator $\mathcal{T}_M$. In virtue of this correspondence we say that a matrix
$M \in  \mathfrak{M}_{a,b}$ is semistandard if the corresponding generalized tableau is semistandard.

Finally, Theorem 2.10.1 can be translated in the following way:

\begin{theorem}
Suppose that $\lambda,\mu$ are partitions of $n$. Then the set $\{\mathcal{T}_M\vert_{S^\lambda}: M\in \mathfrak{M}_{\lambda,\mu}, \mbox{\rm $M$ is semistandard}\}$
is a basis for $\text{\rm Hom}_{S^n}(S^\lambda,M^\mu)$.
\end{theorem}
 In particular, $S^\lambda$ is contained in $M^\mu$ if and only if $\lambda\trianglerighteq \mu$ (Young's rule).

%%%%%%%%%%%%%%%%%%%%%%%%%%%%%%%%%%%%%%%%%%%%%%%%%%%%%%%%%%%%%%%%%%%%%%%%%%%%%%%%%%
\subsection{A Gelfand-Tsetlin decomposition of $M^{a,b,c}$}\label{GZMabc}
%%%%%%%%%%%%%%%%%%%%%%%%%%%%%%%%%%%%%%%%%%%%%%%%%%%%%%%%%%%%%%%%%%%%%%%%%%%%%%%
In \cite{ScarabottiSabc}, the first named author computed an
orthogonal basis for the $S_a\times S_b\times S_c$-invariant vectors
in the irreducible representation $S^{\alpha,\beta,\gamma}$. Indeed,
it was computed the Gelfand-Tsetlin basis associated to the chain of
subgroups $S_{a+b+c}\geq S_{a+b}\times S_c\geq S_a\times S_b \times
S_c$. However, the results in \cite{ScarabottiSabc} are expressed in
terms of a complicated family of orthogonal polynomials in four
variables. In the present subsection, we want to find the
corresponding decomposition of $M^{a,b,c}$ in a completely different
way: we use the semistandard basis in the preceding subsection. We
begin with some general notions. On the space $M^{a_1, a_2, \ldots,
a_k}$ we define the following operators

\[
d_{i,j}\delta_{(A_1,A_2, \ldots ,A_k)} = \sum_{x \in A_j}\delta_{(A_1, \ldots, A_i \cup\{x\}, \ldots , A_j \setminus \{x\},
\ldots, A_k)}
\]

while

\[
\Delta_{i,j}\delta_{(A_1,A_2, \ldots ,A_k)} = \sum_{\substack{x \in A_j\\ y \in A_i}}
\delta_{(A_1, \ldots,\left( A_i \setminus \{y\}\right)\cup\{x\}, \ldots , \left( A_i \setminus \{x\}\right)\cup\{y\},
\ldots, A_k)}.
\]

The case $k= 3$ has been studied in \cite{ScarabottiSabc} and we have the following results.
We recall the decomposition of $M^{a,b}$ into $S_{a+b}$-irreducible subspaces:

\[
M^{a,b}=\bigoplus\limits_{k=\max\{0,b-a\}}^b(d_{2,1}^k)\left[M^{a+k,b-k}\cap \ker d_{1,2}\right]
\]

where $M^{a+k,b-k}\cap \ker d_{1,2}$

is isomorphic to the irreducible representation $S^{a+k,b-k}$.

\begin{lemma}\cite{ScarabottiSabc}
Suppose that $(a,b,c)$ is a composition of $n$. Then we have:

\begin{enumerate}
\item{$(d_{2,1})^{k}\left[M^{a+k,b-k,c}\cap \ker d_{1,2}\right]$ is an eigenspace of $\Delta_{1,2}$ and the corresponding
eigenvalue is: $ab-(b-k)(a+k+1)$.}
\item{ The following is an orthogonal decomposition into $S_n$-invariant subspaces:
\[
M^{a,b,c} = \bigoplus_{k=\max\{0,b-a\}}^b(d_{2,1})^k\left[M^{a+k,b-k,c}\cap \ker d_{1,2}\right].
\]
}
\item{In the permutation module $$M^{a,b,c} =   \text{\rm Ind}_{S_{a+b} \times S_c}^{S_{a+b+c}}\left[M^{a,b}\otimes S^{(c)}\right]$$
the subspace $(d_{2,1})^k\left[M^{a+k,b-k,c}\cap \ker d_{1,2}\right]$ corresponds to
$$\text{\rm Ind}_{S_{a+b} \times S_c}^{S_{a+b+c}}\left[S^{a+k,b-k} \otimes S^{(c)}\right].$$}
\end{enumerate}
\end{lemma}

\begin{corollary}\label{charS}
The subspace $\text{\rm Ind}_{S_{a+b} \times S_c}^{S_{a+b+c}}\left[S^{a+k,b-k}\otimes S^{(c)}\right]$ of $M^{a,b,c}$ may be characterized as the eigenspace of $\Delta_{1,2}$ corresponding to the eigenvalue $ab-(b-k)(a+k+1)$.
\end{corollary}
\begin{proof}
The function $k\mapsto ab-(b-k)(a+k+1)$ is decreasing for $\max\{0,b-a\}\leq k\leq b$.
\end{proof}

Now suppose that $(a,b,c)$ is a partition, that is $a\leq b\leq c$. Let $(\alpha,\beta,\gamma)$ be another partition of $n$ such that $(\alpha,\beta, \gamma)\trianglerighteq (a,b,c)$ (that is $\alpha\geq a$ and $\alpha+\beta\geq a+b$), so that $S^{\alpha,\beta,\gamma}$ is contained in $M^{a,b,c}$. Consider the matrices

\[ M_l =
\begin{pmatrix}
a & l & \alpha -a -l\\
0 & b-l & \beta -b +l\\
0 & 0 & \gamma
\end{pmatrix}\]
\[
 M_l' =
\begin{pmatrix}
a-1 & l & \alpha -a -l+1\\
1& b-l & \beta -b+l -1\\
0 & 0 & \gamma
\end{pmatrix}\]
\[
 M_l'' =
\begin{pmatrix}
a & l & \alpha -a -l+1\\
0 & b-l & \beta -b +l-1\\
0 & 0 & \gamma
\end{pmatrix}.
\]

and let

\[
\begin{array}{cccc}
 \mathcal{T}_l: & M^{\alpha,\beta, \gamma} & \to & M^{a,b,c}\\
 \mathcal{T}_l': & M^{\alpha,\beta, \gamma} & \to & M^{a,b,c}\\
 \mathcal{T}_l'': & M^{\alpha+1,\beta-1, \gamma} & \to & M^{a,b,c}
\end{array}
\]

be the corresponding operators. For instance, if  for  $(D,E,F) \in
\Omega_{\alpha, \beta, \gamma}$ and $(A,B,C) \in  \Omega_{a, b, c}$,
  we set

\[
m\left((D,E,F), (A,B,C)\right) =
\begin{pmatrix}
|D\cap A| & |D\cap B| & |D\cap C|\\
|E\cap A| & |E\cap B| & |E\cap C|\\
|F\cap A| & |F\cap B| & |F\cap C|
\end{pmatrix},
\]

then

\[
\mathcal{T}_l\delta_{(D,E,F)} = \sum_{\substack{(A,B,C) \in \Omega_{a,b,c}:\\ m\left((D,E,F), (A,B,C)\right) = M_l}}
\delta_{(A,B,C)}.
\]

Note that $\mathcal{T}_l$ is semistandard when
 \begin{center}
\begin{picture}(400,60)
 \put(100,15){\line(1,0){75}}
\put(100,30){\line(1,0){150}}
\put(100,45){\line(1,0){150}}
\put(175,15){\line(0,1){15}}
\put(100,0){\line(0,1){45}}
\put(137,15){\line(0,1){15}}
\put(145,30){\line(0,1){15}}
\put(190,30){\line(0,1){15}}
\put(250,30){\line(0,1){15}}
\put(110,34){$1 \cdots 1$}
\put(154,34){$2\cdots 2$}
\put(210,34){$3 \cdots 3$}
\put(105,19){$2 \cdots 2$}
\put(140,19){$3 \cdots 3$}
\put(102,4){$3 \cdots 3$}
\put(100,0){\line(1,0){33}}
\put(133,0){\line(0,1){15}}

\put(123,50){$a$}
\put(170,50){$l$}
\put(200,50){$\alpha-a-l$}
\put(75,19){$b-l$}
\put(178,19){$\beta-b+l$}
\put(140,3){$ \gamma$}
\put(280,34){$\alpha$}
\put(280,19){$\beta$}
\put(280,4){$\gamma$}
\end{picture}
\end{center}
\[
\max\{0,b-a, \beta-a,  b-\beta\}\leq l\leq \min\{b-\gamma,\alpha-a\}.
\]

Since any semistandard tableaux of shape $(\alpha,\beta,\gamma)$ and content $(a,b,c)$ is as above, we can say that $\bigl\{\mathcal{T}_l|_{S^{\alpha,\beta,\gamma}}:\max\{0,b-a, \beta-a,  b-\beta\}\leq l\leq \min\{b-\gamma,\alpha-a\}\bigr\}$ is a basis for $\text{Hom}_{S_n}(S^{\alpha,\beta,\gamma},M^{a,b,c})$.

\begin{lemma}\label{zeil1}
\begin{enumerate}
\item{
\[\mathcal{T}_l''d_{1,2} =  \mathcal{T}_l' + (b-l+1) \mathcal{T}_{l-1}+ (\beta-b+l) \mathcal{T}_l.
\]}
\item{
\[
\Delta_{1,2} \mathcal{T}_l = al  \mathcal{T}_l + (l+1) \mathcal{T}_{l+1}'.
\]}
\end{enumerate}
\end{lemma}
\begin{proof}

First of all we observe that if $(D,E,F) \in \Omega_{\alpha,\beta,\gamma}$ then

\[
\begin{split}
\mathcal{T}_l''d_{1,2}\delta_{(D,E,F)} & = \sum_{\substack{(D',E',F') \in \Omega_{\alpha+1, \beta-1,\gamma}:\\
D \subseteq D', E \subseteq E', F = F'}}\mathcal{T}_l''\delta_{(D',E',F')}\\
& = \sum_{\substack{(D',E',F') \in \Omega_{\alpha+1, \beta-1,\gamma}:\\
D \subseteq D', E \subseteq E', F = F'}}
\sum_{\substack{(A,B,C) \in \Omega_{a, b,c}:\\
m\left( (D',E', F'), (A,B,C)\right) = M_l''}}\delta_{(A,B,C)}\\
= & \sum_{(A,B,C) \in \Omega_{a, b,c}}\xi(A,B,C)\delta_{(A,B,C)}
\end{split}
\]

where
\begin{multline}\label{forsabc}
\xi(A,B,C) = |\{(D',E',F') \in \Omega_{\alpha+1, \beta-1,\gamma}: D \subseteq D', E \subseteq E', F = F' \mbox{ and } \\
m\left( (D',E', F'), (A,B,C)\right) = M_l''\}|.
\end{multline}

Suppose that $(D,E,F) \in \Omega_{\alpha,\beta,\gamma}$ and $(A,B,C) \in \Omega_{a,b,c}$ and that  there exists
 $(D',E', F') \in \Omega_{\alpha+1, \beta-1,\gamma}$ such that  $D \subseteq D', E \subseteq E', F = F'$ and
$m\left( (D',E', F'), (A,B,C)\right) = M_l''$. This means that

\[
\begin{pmatrix}
|D'\cap A| & |D'\cap B| & |D'\cap C|\\
|E'\cap A| & |E'\cap B| & |E'\cap C|\\
|F'\cap A| & |F'\cap B| & |F'\cap C|
\end{pmatrix} =
\begin{pmatrix}
a & l & \alpha -a -l+1\\
0 & b-l & \beta -b +l-1\\
0 & 0 & \gamma
\end{pmatrix}
\]

and that there exists $x \in D'$ such that

\[
D = D'\setminus \{x\} \ \ \ \mbox{ and } \ \ \  E = E' \cup \{x\}.
\]

Therefore
\begin{itemize}
\item{if  $x \in D'\cap A$ then $m\left((D,E,F), (A,B,C)\right) = M_l'$}
\item{if  $x \in D'\cap B$ then $m\left((D,E,F), (A,B,C)\right) = M_{l-1}$}
\item{if  $x \in D'\cap C$ then $m\left((D,E,F), (A,B,C)\right) = M_l$.}
\end{itemize}
In other words, $m\left((D,E,F), (A,B,C)\right) \in \{M_l', M_{l-1}, M_l\}$.
On the other hand, the number $\xi(A,B,C)$ in (\ref{forsabc}) is equal to
\begin{itemize}
\item{$1$ if $m\left((D,E,F), (A,B,C)\right) = M_l'$}
\item{$b-l+1$ if $m\left((D,E,F), (A,B,C)\right) = M_{l-1}$}
\item{$\beta-b+l$ if $m\left((D,E,F), (A,B,C)\right) = M_l.$}
\end{itemize}
Indeed, in these three cases, $D' = D \cup \{x\}$ and $E' =E \setminus \{x\}$, where $x$ belongs respectively
to  $E\cap A$, $E\cap B$ and $E\cap C$ and the cardinalities of these sets are
\begin{itemize}
\item{$|E\cap A| = 1$ if $m\left((D,E,F), (A,B,C)\right) = M_l'$}
\item{$|E \cap B| = b-l+1$ if $m\left((D,E,F), (A,B,C)\right) = M_{l-1}$}
\item{$|E \cap C|=\beta-b+l$ if $m\left((D,E,F), (A,B,C)\right) = M_l.$}
\end{itemize}

This gives 1).

The proof of 2)  is similar and it is left to the reader.
\end{proof}

Recall that $S^{\alpha,\beta,\gamma}=M^{a,b,c}\cap \ker d_{1,2}\cap \ker d_{2,3}$ (James' intersection kernels Theorem: see \cite{James,Sternberg} and \cite{ScarabottiForum} for an elementary proof).

\begin{corollary}\label{c;jelle}
Set $\mathcal{J}_l = \mathcal{T}_l\vert{_{ S^{\alpha,\beta,\gamma}}}$. Then

\[
\Delta_{1,2}\mathcal{J}_l = [al-(l+1)(b-l)]\mathcal{J}_l-(l+1)(\beta-b+l+1)\mathcal{J}_{l+1}.
\]

\end{corollary}
\begin{proof}

From 1 in Lemma \ref{zeil1} (with $l+1$ in place of $l$) we deduce that, for $f \in M^{\alpha, \beta,
\gamma}\cap \ker d_{1,2}$, we have

\[
0 = \mathcal{T}'_{l+1} f + (b-l)\mathcal{T}_lf + (\beta-b+l+1)\mathcal{T}_{l+1} f,
\]

and therefore, by 2

\[
\begin{split}
\Delta_{1,2}\mathcal{T}_l f & = al \mathcal{T}_l f + (l+1)\left[-(b-l)\mathcal{T}_l f
-(\beta-b+l+1)\mathcal{T}_{l+1} f\right]\\
& = \left[al -(l+1)(b-l)\right]\mathcal{T}_l f- (l+1)(\beta-b+l+1)\mathcal{T}_{l+1} f.
\end{split}
\]

Since $S^{\alpha,\beta,\gamma} \subset M^{\alpha,\beta,\gamma}\cap \ker d_{1,2}$,
we have the statement.
\end{proof}

For $a \in \RR$ and  $k \in \NN$, we set

\[
(a)_k = a(a+1)\cdots(a+k-1), \ \ \ \ (a)_0 = 1.
\]

\begin{theorem}\label{t;sabc}
For $\max\{0,b-a, \beta-a, b-\beta\}\leq k \leq \min\{b-\gamma, \alpha-a\}$, the operator

\[
R_k=\sum_{l=k}^{\min\{b-\gamma,\alpha-a\}}\binom{l}{k}\frac{(\beta-b+k+1)_{l-k}}{(a-b+2k+2)_{l-k}}\mathcal{J}_l
\]

intertwines  $S^{\alpha, \beta, \gamma} = M^{\alpha, \beta, \gamma} \cap \ker d_{1,2} \cap \ker d_{2,3}$ with
$(d_{2,1})^k\left[M^{a+k,b-k,c}\cap \ker d_{1,2}\right]$, which is isomorphic to
$\text{\rm Ind}_{S_{a+b}\times S_c}^{S_{a+b+c}}\left[S^{a+k,b-k}\otimes S^{(c)}\right]$.
In particular,

\[
\bigoplus_{k = \max\{0,b-a, \beta-a, b-\beta\}}^{\min\{b-\gamma, \alpha-a\}}R_k S^{\alpha, \beta, \gamma}
\]

is the Gelfand-Tsetlin decomposition of the $S^{\alpha, \beta,
\gamma}$-isotypic component of $M^{a,b,c}$ corresponding to the chain of subgroups $S_{a+b+c}\geq S_{a+b}\times S_c\geq S_a\times S_b\times S_c$.
\end{theorem}
\begin{proof}
We look at the eigenvectors  of $\Delta_{1,2}$ of the form

\[
\sum_l \omega(l)\mathcal{J}_lf
\]

where $f \in S^{\alpha, \beta,\gamma}$ and $\omega(l)$ are coefficients to determine.
In virtue of Corollary \ref{c;jelle}, we have

\[
\Delta_{1,2}\sum_l\omega(l)\mathcal{J}_l =
\sum_l\left\{\omega(l)\left[al-(l+1)(b-l)\right]-\omega(l-1)l(\beta-b+l)\right\}\mathcal{J}_l.
\]

Therefore we must solve the eigenvalue problem
\begin{equation}\label{eautov}
\left\{\begin{array}{l}
\omega(l)\left[al-(l+1)(b-l)\right]- \omega(l-1)l(\beta-b+l)  = \lambda \omega(l)\\
l =  \max\{0,b-a, \beta-a, b-\beta\}, \ldots, \min\{b-\gamma, \alpha-a\}.
\end{array}
\right.
\end{equation}
Since the matrix associated with the system (\ref{eautov}) is upper triangular, the eigenvalues are the
diagonal coefficients:
\[
\left\{\begin{array}{l}
\lambda_k = ak-(k+1)(b-k)\\
k = \max\{0,b-a, \beta-a, b-\beta\}, \ldots, \min\{b-\gamma, \alpha-a\}.
\end{array}
\right.
\]
In order to determine the eigenvectors corresponding to $\lambda_k$, it suffices to set $\omega(l) = 0$ for $l<k$; therefore
(\ref{eautov}) becomes
\[
\left\{\begin{array}{l}
\omega(l)\left[al-(l+1)(b-l)\right]-\omega(l-1)l(\beta-b+l) = \left[ak-(k+1)(b-k)\right]\omega(l)\\
l = k, k+1, \ldots, \min\{b-\gamma, \alpha-a\}
\end{array}
\right.
\]

which is solved by  $\omega(k) = 1$ and, recursively,

\[
\begin{split}
\omega(l)  & = \frac{l(\beta-b+l)}{(l-k)(a-b+l+k+1)}\omega(l-1)\\
& = \cdots = \frac{l(l-1)\cdots (k+1)}{(l-k)!}\cdot\frac{(\beta-b+l)\cdots (\beta-b+k+1)\omega(k)}{(a-b+l+k+1)\cdots (a-b+2k+2)}\\
& = \binom{l}{k}\frac{(\beta-b+k+1)_{l-k}}{(a-b+l+2k+2)_{l-k}}.
\end{split}
\]

We can end the proof by invoking Corollary \ref{charS}.
\end{proof}

\begin{corollary}
$S^{\alpha,\beta,\gamma}$ is contained in $\text{\rm Ind}_{S_{a+b}\times S_c}^{S_{a+b+c}}[S^{a+k,b-k}\otimes S^{(c)}]$
if and only if
\[
\max\{ \beta-a, b-\beta\}\leq k \leq \min\{b-\gamma, \alpha-a\}
\]
and the multiplicity is always equal to 1.
\end{corollary}

\begin{proof}
Follows immediately from the above theorem.
\end{proof}

The operator  $\Delta_{1,2}+ \Delta_{1,3}+ \Delta_{2,3}$ belongs to the center of $\text{Hom}_G(L(X),L(X))$ and therefore
each isotypic component is an eigenspace. Its spectrum has been determined in \cite{ScarabottiSabc}, Theorem 4.3.
This, coupled with  the above theorem, gives

\begin{corollary}
The eigenvalue of the operator $\Delta_{2,3} + \Delta_{1,3}$ restricted to the subspace of $\text{\rm Ind}_{S_{a+b}\times S_c}^{S_{a+b+c}}[S^{a+k,b-k}\otimes S^{(c)}]$ isomorphic to $S^{\alpha,\beta,\gamma}$ is
\[
\frac{1}{2}[\alpha^2+ \beta^2+\gamma^2-2\beta-4\gamma-a^2-b^2-c^2]-ak+(k+1)(b-k)],
\]

for
$\max\{ \beta-a, b-\beta\}\leq k \leq \min\{b-\gamma, \alpha-a\}$.

\end{corollary}

\begin{example}
{\rm
Assume the notations in subsection \ref{decMn211}. Then the semistandard basis of $\text{Hom}_{S_n}(S^{n-1,1},M^{n-2,1,1})$ is given by the operators $\mathcal{J}_0,\mathcal{J}_1$, where

\[
(\mathcal{J}_0 f)(i,j)=f(i) \qquad\quad\text{and}\quad\qquad (\mathcal{J}_1 f)(i,j)=f(j),
\]

for all $f\in S^{n-1,1}$. The Gelfand-Tsetlin decomposition is given by the operators: $R_0=\mathcal{J}_0+\frac{1}{n-1}\mathcal{J}_1$ and $R_1=\mathcal{J}_1$. This agrees with the characterizations in \eqref{Sn11} and \eqref{Sn12}.

}
\end{example}

\begin{example}
{\rm
Now suppose that $(a,b,c)=(n-3,2,1)$ and $(\alpha,\beta)=(n-2,2), \gamma=0$. Identify $\Omega_{n-2,2}$ with the set of all unordered pairs $\{i,j\}$ where $i,j\in I(n), i\neq j$, and $\Omega_{n-3,2,1}$ with the set of all ordered pairs $(\{i,j\},k)$, where $\{i,j\}\in\Omega_{n-2,2}$ and $k\notin \{i,j\}$. In particular, $S^{n-2,2}$ is the space of all functions $f\in M^{n-2,2}$ such that $\sum\limits_{\substack{i=1\\ i\neq j}}^jf(\{i,j\})=0$ for all $j\in I(n)$. Moreover, the semistandard basis of the space $\text{Hom}_{S_n}(S^{n-2,2},M^{n-3,2,1})$ is given by the operators $\mathcal{J}_0,\mathcal{J}_1$ where

\[
(\mathcal{J}_0 f)(\{i,j\},k)=f(\{i,j\}) \qquad\quad\text{and}\quad\qquad (\mathcal{J}_1 f)(\{i,j\},k)=f(\{i,k\})+f(\{j,k\}).
\]

Now the Gelfand-Tsetlin decomposition is given by the operators $R_0=\mathcal{J}_0+\frac{1}{n-3}\mathcal{J}_1$ and $R_1=\mathcal{J}_1$.
}
\end{example}

%%%%%%%%%%%%%%%%%%%%%%%%%%%%%%%%%%%%%%%%%%%%%%%%%%%%%%%%%%%%%%%%%%%%%%%%%%%%
%%%%%%%%%%%%%%%%%%%%%%%%%%%%%%%%%%%%%%%%%%%%%%%%%%%%%%%%%%%%%%%%%%%%%%%%%%%%
\subsection{Concluding remarks}
%%%%%%%%%%%%%%%%%%%%%%%%%%%%%%%%%%%%%%%%%%%%%%%%%%%%%%%%%%%%%%%%%%%%%%%%%%%%
%%%%%%%%%%%%%%%%%%%%%%%%%%%%%%%%%%%%%%%%%%%%%%%%%%%%%%%%%%%%%%%%%%%%%%%%%%%%

The results in this section should be generalized to every isotypic component in any permutation module $M^a$ ($a$ a composition of $n$). But we do not know if the detailed analysis in the special cases presented in this section may be obtained in the general case (cf. the complicated formulas in \cite{ScarabottiSabc}). However, we indicate three lines of research that should be followed to deal with the general case:

\begin{itemize}
\item
The theory of Gelfand-Tsetlin basis, as developed in \cite{Ok-Ver1, Ok-Ver2} for the whole group algebra of the symmetric group (that is, $M^{1^k}$), should be extended to any permutation module $M^a$. In our context, the operators $\sum\limits_{i=1}^{k-1}\Delta_{i,k}$ should play the role of the Young-Jucys-Murphy elements in \cite{Ok-Ver1, Ok-Ver2}.
\item
It should developed a general theory of $H-K$-invariant functions on spaces with multiplicity, along the lines of \cite{DunklSymp} and \cite{ScarabottiRadon}. Indeed, even in the case of the Gelfand pair $(S_{a+b},S_a\times S_b)$, the spherical functions are derived in the more general context of the $S_a\times S_b-S_{a+k}\times S_{b-k}$-invariant functions \cite{CST,CST3}. We mean that the spherical functions in $M^a$ should be studied in the more general context of the $S_a\times S_b$-invariant functions, $b$ another composition of $n$. See again \cite{ScarabottiSabc}.
\item
There should be a clear connections between the spherical functions in $M^a$ and the Clebesh-Gordan coefficients of the unitary group. This has been explored in \cite{Koornwinder} in the case of the Gelfand pair $(S_{a+b},S_a\times S_b)$; see also \cite{Jucys} (which, unfortunately, has not been translated from the Russian) and \cite{JinJalunFan}. The theory of Clebesh-Gordan coefficients of the unitary group has been extensively developed; see \cite{K-V,K-V2}.
\end{itemize}

%%%%%%%%%%%%%%%%%%%%%%%%%%%%%%%%%%%%%%%%%%%%%%%%%%%%%%%%%%%%%%%%%%%%%%%%
%%%%%%%%%%%%%%%%%%%%%%%%%%%%%%%%%%%%%%%%%%%%%%%%%%%%%%%%%%%%%%%%%%%%%%%%
\section{Harmonic analysis of the composition action of a crested product}
%%%%%%%%%%%%%%%%%%%%%%%%%%%%%%%%%%%%%%%%%%%%%%%%%%%%%%%%%%%%%%%%%%%%%%%%
%%%%%%%%%%%%%%%%%%%%%%%%%%%%%%%%%%%%%%%%%%%%%%%%%%%%%%%%%%%%%%%%%%%%%%%%
The results in this section constitute a noncommutative
generalization of the theory crested product of symmetric
association schemes, developed in \cite{Ba-Ca}. We refer also to
\cite{Bailey}, Chapter 10.

%%%%%%%%%%%%%%%%%%%%%%%%%%%%%%%%%%%%%%%%%%%%%%%%%%%%%%%%%%%%%%%%%%%%%%%%
%%%%%%%%%%%%%%%%%%%%%%%%%%%%%%%%%%%%%%%%%%%%%%%%%%%%%%%%%%%%%%%%%%%%%%%%
\subsection{Invariant partitions}
%%%%%%%%%%%%%%%%%%%%%%%%%%%%%%%%%%%%%%%%%%%%%%%%%%%%%%%%%%%%%%%%%%%%%%%%
%%%%%%%%%%%%%%%%%%%%%%%%%%%%%%%%%%%%%%%%%%%%%%%%%%%%%%%%%%%%%%%%%%%%%%%%

In this subsection, we give a noncommutative generalization of Theorem 10.5 in \cite{Bailey}. Moreover, we use the algebra of bi-$K$-invariant function in place of the isomorphic Bose-Mesner algebra.\\

Let $F$ be a finite group and $Y=F/H$ a homogeneous $F$-space. Suppose that $y_0$ is the point stabilized by $H$ and that $Y=\coprod_{j\in\mathcal{J}}\Lambda_j$ is the decomposition of $Y$ into $H$-orbits, with $0\in \mathcal{J}$ and $\Lambda_0=\{y_0\}$. If we set $\tilde{\Lambda}_j=\{(uy,uy_0):u\in F,y\in \Lambda_j\}$ then $Y\times Y=\coprod_{j\in \mathcal{J}}\tilde{\Lambda}_j$ is the decomposition of $Y\times Y$ under the diagonal action of $F$; clearly

\begin{equation}\label{lambdatilde}
\Lambda_j=\{y\in Y:(y,x_0)\in\tilde{\Lambda}_j\}.
\end{equation}

Let $\mathcal{Q}$ be a partition of $Y$, that is $Y=\coprod_{B\in\mathcal{Q}}B$.
Let $\sim_{\mathcal{Q}}$ be the associated equivalence relation, that is
$y\sim_\mathcal{Q}y'$ if and only if $y$ and $y'$ are in the same part of $\mathcal{Q}$.
Define $R_\mathcal{Q}$ by setting

\[
R_{\mathcal{Q}}(y,y')=\left\{\begin{array}{ll}
1&\text{if}\quad y\sim_\mathcal{Q}y'\\
0&\text{otherwise.}
\end{array}\right.
\]

Suppose that $\mathcal{Q}$ is $F$-invariant; this means that if $y\sim_{\mathcal{Q}}y'$ and
$u\in F$ then $uy\sim_{\mathcal{Q}}uy'$. Clearly, $\mathcal{Q}$ is $F$-invariant if and only
$R_\mathcal{Q}$ is constant on the orbits of $F$ on $Y\times Y$; if this is the case, there exists $\mathcal{J}_0\subseteq \mathcal{J}$ such that

\begin{equation}\label{erreq}
R_\mathcal{Q}=\sum_{j\in\mathcal{J}_0}\mathbf{1}_{\tilde{\Lambda}_j}.
\end{equation}

Set $B_0=\{y\in Y:y\sim_{\mathcal{Q}}y_0\}$, that is $B_0$ is the part of $\mathcal{Q}$ containing $y_0$.
Then \eqref{lambdatilde} and \eqref{erreq} ensure us that

\begin{equation}\label{bzero}
\mathbf{1}_{B_0}=\sum_{j\in\mathcal{J}_0}\mathbf{1}_{\Lambda_j}.
\end{equation}

Let $S=\{s\in F:sB_0=B_0\}$ be the stabilizer of $B_0$. Note that $S=\{s\in F:sy_0\sim_{\mathcal{Q}}y_0\}$; in particular, $H\leq S$.

\begin{lemma}\label{equivrel}
Define a relation in $\mathcal{J}$ by setting $i\sim j$ if there exists $s\in S$ such that
$\Lambda_i\cap s\Lambda_j\neq\emptyset$. Then $\sim$ is an equivalence relation. Moreover,

\begin{enumerate}

\item\label{invariance} if $i\in\mathcal{J}$ then the cardinality of the set

\begin{equation}\label{sinS}
\{s\in S:sy'\in\Lambda_i\}
\end{equation}

is the same for all $y'\in\coprod\limits_{\substack{j\in\mathcal{J}\\  i\sim j}}\Lambda_j$;

\item if $[i]$ denotes the equivalence class containing $i$, $\mathcal{J}/\sim$ is the quotient
set and $\Lambda_{[i]}=\coprod_{j\in[i]}\Lambda_j$ then
\begin{equation}\label{rightidealpart}
Y=\coprod_{[i]\in\mathcal{J}/\sim}\Lambda_{[i]}
\end{equation}
is the partition of $Y$ into $S$-orbits. Moreover, $[0]=\mathcal{J}_0$ and $\Lambda_{[0]}=B_0$.
\end{enumerate}

\end{lemma}

\begin{proof}
First of all, we prove \ref{invariance}.
Suppose that $i\sim j$ and $y',y_1'\in\Lambda_j$. Since $\Lambda_j$ is an $H$-orbit,
there exists $h\in H$ such that $hy'_1=y'$. Then the bijection

\[
\begin{array}{ccc}
\{s\in S:sy'\in\Lambda_i\}&\to&\{s_1\in S:s_1y_1'\in\Lambda_i\}\\
s&\longmapsto&s_1=sh
\end{array}
\]

show that the cardinality of \eqref{sinS} does not depend on the choice of $y'\in\Lambda_j$ (in particular, \eqref{sinS} is non empty for any $y'\in\Lambda_j$ with $i\sim j$). Suppose that also $i\sim k$
and take $y''\in\Lambda_k$. There exists $s_1,s_2\in S$ such that $s_1y',s_2y''\in\Lambda_i$,
and therefore if $hs_2y''=s_1y'$, with $h\in H$, we have $y'=s_0y''$, with
$s_0=s_1^{-1}hs_2\in S$.
Then the bijection

\[
\begin{array}{ccc}
\{s\in S:sy'\in\Lambda_i\}&\to&\{s'\in S:s'y''\in\Lambda_i\}\\
s&\longmapsto&s'=ss_0
\end{array}
\]

shows that the cardinality of \eqref{sinS} does not depend on $j$. We have proved
\ref{invariance}.\\

Now we prove that $\sim$ is an equivalence relation. It is reflexive because $1_F\in S$. Taking $s^{-1}$, we get immediately the symmetricity. Suppose that $i\sim j$
and $j\sim k$. Then there exist $y''\in\Lambda_k$ and $s'\in S$ such that
$y'=s'y''\in\Lambda_j$. Using \ref{invariance}, we can say that there exists $s\in S$ such that
$sy'\in\Lambda_i$. Then $ss'y''\in\Lambda_i$ and $i\sim k$. Finally, every $S$-orbit is
$H$-invariant, because $H\leq S$, and therefore it decomposes into a disjoint union of $H$-orbits.
Then the definition of $\sim$ ensures us that $\Lambda_i$ and $\Lambda_j$ are in the same
$S$-orbit if and only if $i\sim j$.

\end{proof}

Now suppose that $Y=\coprod_{C\in\mathcal{T}}C$ is an $H$-invariant partition of $Y$. Then
$\mathbf{1}_C\in L(H\backslash F/H)$ for each $C\in\mathcal{T}$.
We will say that $\mathcal{T}$ is a {\em right ideal partition} of $Y$ if
the vector space $\langle\mathbf{1}_C:C\in\mathcal{T}\rangle$ is a right ideal in the algebra $L(H\backslash F/H)$. This means that $f\in L(H\backslash F/H)$ and $C\in\mathcal{T}$ implies that $\mathbf{1}_C*f\in \langle\mathbf{1}_C:C\in\mathcal{T}\rangle)$.

\begin{lemma}\label{idealpart}
The $H$-partition \eqref{rightidealpart} is a right ideal partition of $Y$.
\end{lemma}

\begin{proof}
Let $\mathcal{J}_0$ and $B_0$ be as in \eqref{bzero}. Denote by $\mathcal{V}$
the subspace of $L(H\backslash F/H)$ spanned by the functions
$\{\mathbf{1}_{\Lambda_{[i]}}:[i]\in\mathcal{J}/\sim\}$.
For every $i\in\mathcal{J}$, denote by $m_i$ the cardinality of the set \eqref{sinS} divided by $\lvert H\rvert$.
By definition of convolution, if $i\sim_\mathcal{Q}j$, $t\in F$ and $ty_0\in \Lambda_j$ then

\[
\begin{split}
\mathbf{1}_{B_0}*\mathbf{1}_{\Lambda_i}(ty_0)=&
\frac{1}{\lvert H\rvert}\sum_{u\in F}\mathbf{1}_{B_0}(uy_0)\mathbf{1}_{\Lambda_i}(u^{-1}ty_0)\\
=&\frac{1}{\lvert H\rvert}\sum_{s\in S}\mathbf{1}_{\Lambda_i}(s^{-1}ty_0)\\
=&\frac{1}{\lvert H\rvert}\lvert\{s\in S:s^{-1}ty_0\in\Lambda_i\}\rvert=m_i,
\end{split}
\]

and therefore

\begin{equation}\label{1B0}
\mathbf{1}_{B_0}*\mathbf{1}_{\Lambda_i}=m_i\mathbf{1}_{\Lambda_{[i]}}\in\mathcal{V}
\end{equation}

for every $i\in\mathcal{J}$, which implies that $\mathbf{1}_{B_0}*f\in\mathcal{V}$ for every $f\in L(H\backslash F/H)$. Another application of \eqref{1B0} yields

\[
\mathbf{1}_{\Lambda_{[i]}}*\mathbf{1}_{\Lambda_j}=\frac{1}{m_i}\mathbf{1}_{B_0}*\mathbf{1}_{\Lambda_i}*\mathbf{1}_{\Lambda_j}\in\mathcal{V},
\]

for all $[i]\in\mathcal{J}/\sim$ and $j\in\mathcal{J}$, and therefore
$\mathcal{V}$ is a right ideal in $L(H\backslash F/H)$.

\end{proof}

The partition \eqref{rightidealpart} will be called the {\em right ideal partition} associated to the invariant
partition $\mathcal{Q}$.

\begin{example}\label{examplewreath}
{\rm
Let $F,H,y_0,J,\Lambda_j$ be as above (without assuming the existence of an invariant partition) and suppose that $F'$ is another finite group, acting transitively on $Y'=F'/H'$, with $H'$ stabilizer of $y'_0$. Suppose also that $Y'=\coprod_{j\in\mathcal{J}'}\Lambda'_j$ is the decomposition of $Y'$ into $H'$-orbits, with $0\in \mathcal{J}'$ and $\Lambda'_0=\{y'_0\}$. The {\em wreath product} of $F'$ by $F$ is the group $F'\wr F=F'^{Y}\times F\equiv\{(f,u):u\in F,f:Y\rightarrow F'\}$ with the multiplication law $(f,u)(f_1,u_1)=(f\cdot uf_1,uu_1)$, where $(f\cdot uf_1)(y)=f(y)f_1(u^{-1}y)$ for every $y\in Y$. The group $F'\wr F$ acts on $Y\times Y'$ by the {\em composition action} \cite{CST2,JK}: $(f,u)(y,y')=(uy,f(uy)y')$, for all $(f,u)\in F'\wr F$ and $(y,y')\in Y\times Y'$. The stabilizer of $(y_0,y'_0)$ is the subgroup $L=\{(f,u)\in F'\wr F:u\in H,f(y_0)\in H'\}$ and

\begin{equation}\label{decorbwr}
Y\times Y'=\left[\coprod_{j\in \mathcal{J}'}\left(\Lambda_0\times \Lambda'_{j}\right)\right]\coprod \left[\coprod_{j\in \mathcal{J}\setminus\{0\}}\left(\Lambda_j\times Y'\right)\right]
\end{equation}

is the decomposition of $Y\times Y'$ into $L$-orbits \cite{CST2}. Now there is a natural invariant partition on $Y\times Y'$:

\begin{equation}\label{ytimesyprime}
Y\times Y'=\coprod_{y\in Y}B_y,
\end{equation}

where $B_y=\{(y,y'):y'\in Y\}$. It is easy to see that
$(f,u)(y_0,y'_0)=(uy_0,f(uy_0)y'_0)\in B_{y_0}\equiv\Lambda_0\times Y'$ if and only if
$u\in H$; in other words, the stabilizer of $B_{y_0}$ is $F'\wr H$.
It follows that, in the present situation,

\begin{equation}\label{rightidealyyprime}
Y\times Y'=\coprod_{j\in \mathcal{J}}\left(\Lambda_j\times Y'\right)
\end{equation}

is the right ideal partition associated to the invariant partition \eqref{ytimesyprime}.
}
\end{example}

\begin{example}\label{simgroup}
{\rm Take $F=S_n$ and $H=S_{n-2}\times S_1\times S_1$, where $S_{n-2}$ acts
on $\{3,4,\dotsc,n\}$. Then the orbits of $F$ on $F/H$ are
$\Lambda_0=\{(1,2)\}$, $\Lambda_1=\{(2,1)\}$, $\Lambda_2=\{(1,j):j\neq 1,2\}$,
$\Lambda_3=\{(2,j):j\neq 1,2\}$, $\Lambda_4=\{(i,1):i\neq 1,2\}$,
$\Lambda_5=\{(i,2):i\neq 1,2\}$, and $\Lambda_6=\{(i,j):i\neq j\;\text{and}\;i,j\neq 1,2\}$.
Take the invariant partition $G/H=\coprod_{i=1}^nB_i$, where $B_i=\{(i,j):i\neq j\}$. Then
the stabilizer of $B_1$ is $S_{n-1}\times S_1$ ($S_{n-1}$ acts on $\{2,3,\dotsc,n\}$)
and the relation $\sim$ is given by
$0\sim 2$, $1\sim 4$ and $3\sim 5\sim 6$.
}
\end{example}

\begin{remark}
{\rm We can define another equivalence relation on $\mathcal{J}$ by setting $i\approx j$ when
there exist $y\in\Lambda_i$ and $y'\in\Lambda_j$ such that $y\sim_\mathcal{Q}y'$. If $\{i\}$
denotes the equivalence class containing $i$ and $\Lambda_{\{i\}}=\coprod_{j\in\{i\}}\Lambda_j$,
then $\Lambda_{0}=B_{\{0\}}$ and

\begin{equation}\label{approx}
Y=\coprod_{\{i\}\in\mathcal{J}/\approx}\Lambda_{\{i\}}
\end{equation}

is the supremum of the partitions $\mathcal{Q}$ and $Y=\coprod_{j\in\mathcal{J}}\Lambda_j$,
that is it is the finest partition that has both $\mathcal{Q}$ and
$Y=\coprod_{j\in\mathcal{J}}\Lambda_j$ as a refinement. Moreover, arguing as in the proofs of
lemmas \ref{equivrel} and \ref{idealpart}, it may be shown that

\[
\mathbf{1}_{\Lambda_i}*\mathbf{1}_{B_0}=M_i\mathbf{1}_{\Lambda_{\{i\}}},
\]

where, for $y'\in\Lambda_{\{i\}}$, $M_i=\lvert\{y\in\Lambda_i:y\sim_{\mathcal{Q}}y'\}\rvert$ (and this cardinality does
not depend on the choice of $y'$). Therefore,
$\{\mathbf{1}_{\Lambda_{\{i\}}}:\{i\}\in\mathcal{J}/\approx\}$ span a {\em left} ideal in
$L(H\backslash F/H)$. In general, the equivalence relations $\sim$ and $\approx$ are different.
For instance, in Example \ref{simgroup} we have $0\approx 2$, $1\approx 3$ and
$4\approx 5\approx 6$. Clearly, the equivalence relations $\sim$ and $\approx$ coincide if
and only if $\mathbf{1}_{B_0}$ is in the center of $L(H\backslash F/H)$. In particular, they coincide in the case of a symmetric association scheme treated in \cite{Bailey,Ba-Ca}. Another case where they coincide is given by Example \ref{examplewreath}.
}
\end{remark}

%%%%%%%%%%%%%%%%%%%%%%%%%%%%%%%%%%%%%%%%%%%%%%%%%%%%%%%%%%%%%%%%%%%%%%%%
%%%%%%%%%%%%%%%%%%%%%%%%%%%%%%%%%%%%%%%%%%%%%%%%%%%%%%%%%%%%%%%%%%%%%%%%
\subsection{Crested products of finite groups}
%%%%%%%%%%%%%%%%%%%%%%%%%%%%%%%%%%%%%%%%%%%%%%%%%%%%%%%%%%%%%%%%%%%%%%%%
%%%%%%%%%%%%%%%%%%%%%%%%%%%%%%%%%%%%%%%%%%%%%%%%%%%%%%%%%%%%%%%%%%%%%%%%

In the present subsection, we introduce the notion of crested product of finite groups \cite{Bailey, Ba-Ca}. It is a generalization of both the notion of wreath product (see example \ref{examplewreath}) and direct product.

We continue to use all the notation in lemma \eqref{idealpart}. Moreover, we assume that $N=\{u\in F:uB=B\;\text{for all}\;B\in\mathcal{Q}\}$, which is a normal subgroup of $F$, is transitive on every $B\in\mathcal{Q}$ (in other words, we assume that $\mathcal{Q}$ is the orbit partition of its stabilizer $N$).
Let $G$ be another group, acting transitively on a set $X/K$, with $K$ stabilizer of $x_0$. Suppose that $X=\coprod_{i\in\mathcal{I}}\Xi_i$ is the decomposition of $X$ into $K$-orbits, with $0\in \mathcal{I}$ and $\Xi_0=\{x_0\}$, and that $X=\coprod_{A\in\mathcal{P}}A$ is a $G$-invariant partition of $X$. Let $F^\text{diag}$ be the group formed by all functions $f:X\to F$ that are constant on the whole $X$ and $N^\mathcal{P}$ the group of all $f:X\to N$ that are constant on each $A\in \mathcal{P}$. Then $F^\text{diag}$ and $N^\mathcal{P}$ are subgroups of $F^X$, $F^\text{diag}\cap N^\mathcal{P}=N^\text{diag}$ and $N^\mathcal{P}$ is normalized by $F^\text{diag}$. Then $F^\text{diag}\cdot N^\mathcal{P}$ is a subgroup of $F^X$. Identify $\{(f,g)\in F\wr G:f(x)=1_F\;\text{for all}\;x\in X\}$ with $G$ and $\{(f,g)\in F\wr G:g=1_G\}$ with $F^X$. Then $F\wr G\cong F^X\ltimes G$.

\begin{definition}
The {\em crested product} of $F$ by $G$ is the subgroup $(F^\text{\rm diag}\cdot N^\mathcal{P})\ltimes G$ of $F\wr G$.
\end{definition}

The crested product $(F^\text{diag}\cdot N^{\mathcal{P}})\ltimes G$ acts on $X\times Y$ as a subgroup of $F\wr G$, via the composition action (Example \ref{examplewreath}). The stabilizer of $(x_0,y_0)$ is the subgroup

\[
R=\{(f,k)\in(F^\text{diag}\cdot N^{\mathcal{P}})\ltimes G:k\in K,f(x_0)\in H\}.
\]

As in \eqref{bzero}, suppose that $A_0=\{x\in X:x\sim_{\mathcal{P}}x_0\}$ and define $\mathcal{I}_0\subseteq\mathcal{I}$ by requiring that

\begin{equation}\label{T-orbits}
\mathbf{1}_{A_0}=\sum_{i\in\mathcal{I}_0}\mathbf{1}_{\Xi_i}.
\end{equation}

Now we can prove the analogous of \eqref{decorbwr} for a crested product, generalizing Theorem 10 in \cite{Ba-Ca}. See also Lemma 4.1 in \cite{CST2}.

\begin{theorem}\label{thdecorbcr}
The decomposition of $X\times Y$ into $R$-orbits is:

\begin{equation}\label{decorbcr}
X\times Y=\left[\coprod_{\substack{i\in \mathcal{I}_0\\j\in \mathcal{J}}}\left(\Xi_i\times \Lambda_{j}\right)\right]\coprod \left[\coprod_{\substack{i\in\mathcal{I}\setminus\mathcal{I}_0\\
[j]\in \mathcal{J}/\sim}}\left(\Xi_i\times \Lambda_{[j]}\right)\right].
\end{equation}

\end{theorem}

\begin{proof}
First of all, we prove that all the subsets in the decomposition \eqref{decorbcr} are $R$-invariant. Note that any $(f,k)\in R$ may be written in the form $(f,k)=(f_1,1_G)(f_2,k)$, with $f_2\in H^\text{diag}$, $f_1\in N^{\mathcal{P}}$ and $f_1(x_0)=1_F$. Moreover, $(f_2,k)\left(\Xi_i\times \Lambda_{j}\right)=\Xi_i\times \Lambda_{j}$ for all $i\in\mathcal{I}$ and $j\in\mathcal{J}$. Then we have to prove that every subset in \eqref{decorbcr} is invariant under $(f_1,1_G)$.

If $i\in\mathcal{I}_0,j\in\mathcal{J},x\in\Xi_i$ and $y\in\Lambda_j$ then $x\sim_{\mathcal{P}}x_0$ and therefore $f_1(x)=f_1(x_0)=1_F$. It follows that

\[
(f_1,1_G)(x,y)=(x,f_1(x)y)=(x,y),
\]

and therefore $\Xi_i\times\Lambda_j$ is $(f_1,1_G)$-invariant.

On the other hand, suppose that
$i\in\mathcal{I}\setminus\mathcal{I}_0,[j]\in\mathcal{J}/\sim$, $x\in\Xi_i$ and that
$y\in\Lambda_j\subseteq\Lambda_{[j]}$. Since $f_1(x)\in N\leq S$, if $f_1(x)y\in\Lambda_k$
then $k\sim j$. Therefore $f_1(x)y\in \Lambda_{[j]}$ and

\[
(f_1,1_G)(x,y)=(x,f_1(x)y)\in\Xi_i\times\Lambda_{[j]}.
\]

We can conclude that $\Xi_i\times\Lambda_{[j]}$ is invariant under $(f_1,1_G)$ too.\\

Now we prove that $R$ is transitive on every subset in the right hand side of \eqref{decorbcr}.
Note that $K$ is transitive on every $\Xi_i$, and therefore we can limit ourselves to consider pairs $(x_1,y_1)$, $(x_2,y_2)$ with $x_1=x_2$.

Suppose that $i\in\mathcal{I}_0,j\in\mathcal{J}$ and that
$(x,y_1),(x,y_2)\in \Xi_i\times\Lambda_j$. Taking $(f,1_G)$ with $f\in H^\text{diag}$ such that
$f(x)y_1=y_2$ we get $(f_1,1_G)(x,y_1)=(x,f_1(x)y_1)=(x,y_2)$, and therefore $R$ is transitive on $\Xi_i\times\Lambda_j$.

On the other hand, suppose that $(x,y_1),(x,y_2)\in \Xi_i\times\Lambda_{[j]}$, with
$i\in\mathcal{I}\setminus\mathcal{I}_0$ and $[j]\in\mathcal{J}/\sim$. Since $S$ is transitive on $\Lambda_{[j]}$, there exists
$s\in S$ such that $y_2=sy_1$. Since $N$ is transitive on $B_0$ and $sy_0\in B_0$, there exists $n\in N$ such that $nsy_0=y_0$, and
therefore $ns\in H$. Since $i\notin \mathcal{I}_0$, we have $x\nsim_{\mathcal{Q}}x_0$ and
therefore we can take $f=f_1\cdot f_2\in N^\mathcal{P}\cdot H^\text{diag}$ such that
$f_1(x)=n^{-1}$ and $f_2(x)=ns$. It follows that

\[
(f,1_G)(x,y_1)=(x,f(x)y_1)=(x,n^{-1}nsy_1)=(x,y_2)
\]

and $R$ is transitive on $\Xi_i\times\Lambda_{[j]}$ too.
\end{proof}

\begin{example}\label{crossnested}{\rm
Let $\mathcal{U}_X$ and $\mathcal{E}_X$ be respectively the {\em universal partition} (the partition with a single part) and the {\em equality partition} (every part is a singleton) of the set $X$; similarly for  $\mathcal{U}_Y$ and $\mathcal{E}_Y$. If we take  $\mathcal{P} =\mathcal{U}_X$ and $\mathcal{Q}=\mathcal{E}_Y$, then the resulting crested product is isomorphic to the direct product $F\times G$.
Now the composition action coincides with the direct product of permutation representations and \eqref{decorbcr} becomes $X\times Y=\coprod_{i\mathcal{I},j\in\mathcal{J}}(\Xi_i\times \Lambda_j)$.
On the other hand,
if we take  $\mathcal{P}=\mathcal{E}_X$ and $\mathcal{Q}=\mathcal{U}_Y$, then the resulting crested product coincides with the whole wreath product $F\wr G$. In this case, we get the usual composition action and \eqref{decorbcr} becomes \eqref{decorbwr} (translated from the $Y\times Y'$ setting to the $X\times Y$ setting.)}
\end{example}

In particular,
Example \ref{crossnested} shows that the crested product generalizes both the direct and the wreath product \cite{Ba-Ca}.

\begin{example}\label{iteratedwreath}
{\rm This is a continuation of Example \ref{examplewreath}. We show how to obtain a crested product in the setting of iterated wreath products. We suppose that $G,K,X,x_0,\mathcal{I}$ and $\Xi_i$ are as above, but we do not assume that there exists an invariant partition on $X$. On the contrary, suppose that $G'$ is another finite group, acting on $X'=G'/K'$, $K'$ the  stabilizer of $x'_0$ and $X'=\coprod_{i\in\mathcal{I}'}\Xi'_i$ the decomposition of $X'$ into $K'$ orbits, with $\Xi'_0=\{x'_0\}$. Then the wreath product $G'\wr G$ acts on $X\times X'$ via the composition action and we have the invariant partition

\begin{equation}\label{xtimesxprime}
X\times X'=\coprod_{x\in X}A_x,
\end{equation}

where $A_x=\{(x,x'):x'\in x\}$ (Example \ref{examplewreath} in the $X\times X'$ setting).
The iterated wreath product $F'\wr F\wr G'\wr G$ is isomorphic to the set of all $(f',f,g',g)$ where $f':X'\times X\times Y\to F'$, $f:X\times X'\to F$, $g':X\to G'$ and $g\in G$, with the multiplication law

\[
(f',f,g',g)(f_1',f_1,g'_1,g_1)=(f'\cdot (f,g',g)f_1',f\cdot(g',g)f_1,g'\cdot gg'_1,gg_1).
\]

(compare with the labeling used in \cite{Grigorchuk}). Then we can
take \eqref{ytimesyprime} and \eqref{xtimesxprime} as invariant
partitions respectively on $Y\times Y'$ and on $X\times X'$. Denote
by $\mathcal{P}$ the partition \eqref{xtimesxprime}. Then the
resulting crested product is $\left[(F'^Y)^\mathcal{P}\cdot(F'\wr
F)^\text{diag}\right]\ltimes (G'\wr G)$. In other words, the
resulting crested product is the subgroup of $F'\wr F\wr G'\wr G$
formed by all $(f',f,g',g)$ such that: $g\in G,g':X\to G',f\in
F,f':X\times Y\to F'$, that is $f$ is constant and $f'$ does not
depend on $x'\in X'$. The stabilizer of $(x_0,x'_0,y_0.y'_0)$ is the
subgroup $R$ of the crested product formed by all $(f',f,g',k)$ such
that $k\in K$, $g'(x_0)\in K'$, $f\in H$ and $f'(x_0,y_0)\in H'$.
Then Lemma \ref{idealpart}, Theorem \ref{thdecorbcr},
\eqref{decorbwr} and \eqref{rightidealyyprime} yield the
decomposition of $X\times X'\times Y\times Y'$ into $R$-orbits:

\[
\begin{split}
X\times X'\times Y\times Y'=&\left[\coprod_{\substack{i'\in \mathcal{I}'\\j'\in \mathcal{J}'}}
\left(\Xi_0\times \Xi'_i\times \Lambda_0\times\Lambda'_j\right)\right]\coprod
\left[\coprod_{\substack{i'\in\mathcal{I}'\\
j\in \mathcal{J}\setminus\{0\}}}\left(\Xi_0\times\Xi'_i \times\Lambda_j\times Y'\right)\right]\\
&\coprod \left[\coprod_{\substack{i\in\mathcal{I}\setminus\{0\}\\
j\in \mathcal{J}}}\left(\Xi_i\times X' \times\Lambda_j\times Y'\right)\right].
\end{split}
\]

}
\end{example}

%%%%%%%%%%%%%%%%%%%%%%%%%%%%%%%%%%%%%%%%%%%%%%%%%%%%%%%%%%%%%%%%%%%%%%%
%%%%%%%%%%%%%%%%%%%%%%%%%%%%%%%%%%%%%%%%%%%%%%%%%%%%%%%%%%%%%%%%%%%%%
\subsection{The permutation representation of the composition action}
%%%%%%%%%%%%%%%%%%%%%%%%%%%%%%%%%%%%%%%%%%%%%%%%%%%%%%%%%%%%%%%%%%%%%%%
%%%%%%%%%%%%%%%%%%%%%%%%%%%%%%%%%%%%%%%%%%%%%%%%%%%%%%%%%%%%%%%%%%%%%

We need some general facts on the decomposition of permutation representations.
Let $G,X,K$ and $X=\coprod_{i\in\mathcal{J}}\Xi_i$ be as in the preceding subsections. Suppose that $L(X)=\bigoplus_{\omega\in\Omega}a_\omega V_\omega$ is the isotypic decomposition of $L(X)$. Then
the sum of the squares of the multiplicities is equal
to the number of orbits of $K$ on $X$, that is

\begin{equation}\label{Wielandt}
\sum_{\omega\in\Omega}(a_\omega)^2=\lvert\mathcal{I} \rvert.
\end{equation}

This is called {\em Wielandt's Lemma} in \cite{CST,CST2,CST3}. See also \cite{Sternberg,Wielandt}.

\begin{lemma}\label{Wieapplic}
For any orthogonal decomposition of $L(X)$ into $G$-invariant subspaces

\begin{equation}\label{isotypic}
L(X)=\bigoplus_{\gamma\in\Gamma}c_\gamma U_\gamma,
\end{equation}

where every block $c_\gamma U_\gamma$ is the orthogonal sum of $c_\gamma$ invariant $G$-isomorphic subspaces, we have

\[
\sum_{\gamma\in\Gamma}(c_\gamma)^2\leq\lvert\mathcal{I}\rvert
\]

with equality if and only if \eqref{isotypic} is the isotypic decomposition.
\end{lemma}

\begin{proof}
Starting from \eqref{isotypic}, we can get the isotypic decomposition in two stages.

\begin{enumerate}

\item First of all, we can decompose every $U_\gamma$ into irreducible representations:
$U_\gamma=m_{\gamma}W_\gamma\oplus m'_{\gamma}W'_\gamma\oplus\dotsb\oplus m''_\gamma W''_\gamma$,
and then replace $c_\gamma U\gamma$ with $c_\gamma m_{\gamma}W_\gamma\oplus c_\gamma m'_{\gamma}W'_\gamma\oplus\dotsb\oplus c_\gamma m''_\gamma W''_\gamma$. This way, we get a decomposition

\begin{equation}\label{decirred}
L(X)=\bigoplus_{\beta\in B}b_\beta W_\beta
\end{equation}

where each $W_\beta$ is irreducible. Clearly, $\sum_{\beta\in B}(b_\beta)^2\geq\sum_{\gamma\in\Gamma}(c_\gamma)^2$, with equality if and only if all the $U_\gamma$'s are irreducible.

\item We can group together the isomorphic representations in \eqref{decirred}: if
$W_\beta,W_{\beta'},\dotsc, W_{\beta''}$ are all the representations isomorphic to $W_\beta$, we can replace $b_\beta W_\beta\oplus b_{\beta'}W_{\beta'}\oplus\dotsb\oplus b_{\beta''}W_{\beta''}$ with $(b_\beta+ b_{\beta'}+\dotsb +b_{\beta''})W_\beta$. This way, we get a decomposition

\[
L(X)=\bigoplus_{\alpha\in A}a_\alpha V_\alpha
\]

that must coincide with the isotypic one. Again, $\sum_{\alpha\in A}(a_\alpha)^2\geq\sum_{\beta\in B}(b_\beta)^2$, with equality if and only if the representations $W_\beta$'s in \eqref{isotypic} are pairwise inequivalent.

\end{enumerate}

By Wielandt's Lemma \eqref{Wielandt}, we have

\[
\sum_{\gamma\in\Gamma}(c_\gamma)^2\leq\sum_{\beta\in B}(b_\beta)^2
\leq\sum_{\alpha\in A}(a_\alpha)^2=\lvert\mathcal{I}\rvert,
\]

with equality if and only if \eqref{isotypic} is the isotypic decomposition.
\end{proof}

In other words, the isotypic decomposition may be characterized, among those of the form \eqref{isotypic}, as the decomposition that maximizes the quantity $\sum_{\gamma\in\Gamma}(c_\gamma)^2$.\\

Now we return to use all the notation in the preceding subsections. Suppose that $L(X)=\bigoplus_{\omega\in\Omega}a_\omega V_\omega$ and
$L(Y)=\bigoplus_{\delta\in\Delta}b_\delta W_\delta$ are the isotypic decompositions
of $L(X)$ and $L(Y)$ into respectively $G$-irreducible and $F$-irreducible representations.
Let $T$ be the stabilizer of $A_0$ in $G$. That is, $T$ plays, in the $X$-setting, the same role
of $S$ in the $Y$-setting. Then $A_0\cong T/K$ and $G/T\cong \mathcal{P}$.
Suppose that $L(A_0)=\bigoplus_{\gamma\in\Gamma}c_\gamma U_\gamma$ is the decomposition of $L(A_0)$
into irreducible $T$-representations. By transitivity of induction, we have

\begin{equation}\label{indlazero}
L(X)=\text{Ind}_T^GL(A_0)=
\bigoplus_{\gamma\in\Gamma}c_\gamma\text{Ind}_T^G(U_\gamma).
\end{equation}

We will need this simple consequence of the theory developed in subsection \ref{trainddec}.

\begin{lemma}\label{pringen}
For each $A\in\mathcal{P}$, choose an element $l_A\in G$ such that $l_A A_0=A$ (and we always take $l_{A_0}=1_G$). Then, for each copy of $U_\gamma$ in $c_\gamma U_\gamma$,
$\text{\rm Ind}_T^G(U_\gamma)$ is the subspace of $L(X)$ spanned by all functions
$l_A\mathcal{G}$, with $\mathcal{G}\in U_\gamma, A\in\mathcal{P}$. Moreover, for every
$g\in G$ we have $gl_A\mathcal{G}=l_{gA}(t\mathcal{G})$, where $t=(l_{gA})^{-1}gl_A\in T$
(and therefore $t\mathcal{G} \in U_\gamma$).
\end{lemma}

\begin{proof}
Clearly, $\{l_A:A\in\mathcal{P}\}$ is a set of representatives for the left cosets of
$T$ in $G$, that is $G=\coprod_{a\in\mathcal{P}}l_A T$. Then, by definition of induced
representation, $\text{Ind}_T^G(U_\gamma)$ is spanned by the functions of the type $l_A\mathcal{G}$.
The rest of the lemma follows from the simple observation that $gl_A A_0=gA=l_{gA}A_0$,
and therefore there exists $t\in T$ such that $gl_A=l_{gA}t$.

\end{proof}

In the $Y$-setting, we identify $L(\mathcal{Q})\equiv L(F/S)$ with the space
of all $\phi\in L(Y)$ that are constant on each
part of $\mathcal{Q}$. Suppose that $W,W'$ are two irreducible $F$-representations contained
in $L(Y)$, with $W\subseteq L(\mathcal{Q})$ and $W'$ orthogonal to $L(\mathcal{Q})$. Then
$\text{Res}^F_N(W)$ is isomorphic the direct sum of $\text{dim}W$ copies of the trivial
representation of $N$, while $\text{Res}^F_N(W')$ must contain some nontrivial
$N$-representation (by hypothesis, $\mathcal{Q}$ is the orbit partition of $N$).
Then $W$ and $W'$ are not equivalent. It follows that there exists a subset
$\Delta_0\subseteq\Delta$ such that:

\begin{equation}\label{deltauno}
L(\mathcal{Q})=\bigoplus_{\delta\in\Delta_0}b_\delta W_\delta
\end{equation}

is the decomposition of $L(\mathcal{Q})$ into irreducible $F$-representations.

We also need to determine the orbits of $S$ on $\mathcal{Q}$.
Suppose that $B\in\mathcal{Q}$ and $B\cap\Lambda_{[i]}\neq\emptyset$. Then if $y\in B$ and $y'\in B\cap\Lambda_{[i]}$ we can take $n\in N\leq S$ such that $ny'=y$ ($N$ is transitive on $B$), and therefore $B\subseteq\Lambda_{[i]}$ (because $\Lambda_{[i]}$ is an $S$-orbit). Setting

\[
\mathcal{B}_{[i]}=\{B:B\subseteq\Lambda_{[i]}\}
\]

(note that $\mathcal{B}_{[i]}$ is a subset of $\mathcal{Q}$), then

\begin{equation}\label{Sorbq}
\mathcal{Q}=\coprod_{[i]\in\mathcal{J}/\sim}\mathcal{B}_{[i]}
\end{equation}

is the decomposition of $\mathcal{Q}$ into $S$-orbits. Note that we have used the fact that $N$ is transitive on each $B\in\mathcal{Q}$; in Example \ref{simgroup}, $N$ is trivial and \eqref{Sorbq} does not hold.

Now we are in position to get the decomposition of $L(X\times Y)$ into irreducible representations of
the crested product $(N^\mathcal{P}\cdot F^\text{diag})\ltimes G$. Compare with Theorem 4.2 in \cite{CST2} and Theorem F in \cite{BPRS}.

\begin{theorem}\label{deccompactionth}

The following

\begin{equation}\label{deccompaction}
\begin{split}
L(X\times Y)=&\left[\bigoplus_{\omega\in\Omega}\bigoplus_{\delta\in\Delta_0}a_\omega b_\delta
\left(V_\omega\otimes W_\delta\right)\right]\\
&\bigoplus\left\{\bigoplus_{\gamma\in\Gamma}\bigoplus_{\delta\in\Delta\setminus\Delta_0}
c_\gamma b_\delta \left[\left(\text{\rm Ind}_T^GU_\gamma\right)\otimes W_\delta
\right]\right\}
\end{split}
\end{equation}

is the decomposition of $L(X\times Y)$ into $(N^\mathcal{P}\cdot F^\text{\rm diag})\ltimes G$-irreducible
representations.
\end{theorem}

\begin{proof}
Suppose that $\mathcal{G}\in L(X)$, $\mathcal{F}\in L(Y)$, $(x,y)\in X\times Y$ and
$(f_1f_2,g)\in(N^\mathcal{P}\cdot F^\text{diag})\ltimes G$, with
$f_2\in F\equiv F^\text{diag},f_1\in N^\mathcal{P}$. Then

\begin{equation}\label{actensorprod}
\begin{split}
[(f_1f_2,g)(\mathcal{G}\otimes\mathcal{F})](x,y)=&(\mathcal{G}\otimes\mathcal{F})
[(f_1f_2,g)^{-1}(x,y)]\\
=&(\mathcal{G}\otimes\mathcal{F})(g^{-1}x,[f_1(x)f_2]^{-1}y)\\
=&\mathcal{G}(g^{-1}x)\cdot\mathcal{F}([f_1(x)f_2]^{-1}y)\\
=&(g\mathcal{G})(x)\cdot [f_1(x)(f_2\mathcal{F})](y).
\end{split}
\end{equation}

The last expression in \eqref{actensorprod} is deceptive: in general, $f_1(x)(f_2\mathcal{F})$ depends on
$x$, and when this is the case, \eqref{actensorprod} is {\em not} a tensor product. But there are two special cases in which it is a tensor product, and we have to examine these cases in order to show that every subspace in the right hand side of \eqref{deccompaction} is $(N^\mathcal{P}\cdot F^\text{diag})\ltimes G$-invariant. \\

If $\omega\in\Omega$,
$\delta\in\Delta_0$, $\mathcal{G}\in V_\omega$ and $\mathcal{F}\in W_\delta$ then
$g\mathcal{G}\in V_\omega$ and $f_1(x)(f_2\mathcal{F})=f_2\mathcal{F}\in W_\delta$,
because $f_1(x)\in N$ and $N$ acts trivially on each $W_\delta\subseteq L(\mathcal{Q})$.
Then \eqref{actensorprod} yields

\begin{equation}\label{actensorprod2}
(f_1f_2,g)(\mathcal{G}\otimes\mathcal{F})=
(g\mathcal{G})\otimes(f_2\mathcal{F})
\end{equation}

and therefore $V_\omega\otimes W_\delta$ is
$(N^\mathcal{P}\cdot F^\text{diag})\ltimes G$-invariant. \\

On the other hand, suppose that $\gamma\in\Gamma$,
$\delta\in\Delta\setminus\Delta_0$, $l_A\mathcal{G}\in\text{Ind}_T^G
U_\gamma$ (cf. Lemma \ref{pringen}) and that $\mathcal{F}\in
W_\delta$. If $x\notin gA$ and $t$ is as in Lemma \ref{pringen},
then $[l_{gA}(t\mathcal{G})](x)=0$, while if $x\in gA$ and $u\in F$
is the constant value of $f_1$ on $gA$ then \eqref{actensorprod}
yields

\[
\begin{split}
\{(f_1f_2,g)[(l_A\mathcal{G})\otimes\mathcal{F}]\}(x,y)&=[l_{gA}(t\mathcal{G})](x)\cdot
[f_1(x)f_2\mathcal{F}](y)\\
&=\{[l_{gA}(t\mathcal{G})]\otimes [uf_2\mathcal{F}]\}(x,y).
\end{split}
\]

We can conclude that

\begin{equation}\label{actensorprod3}
(f_1f_2,g)[(l_A\mathcal{G})\otimes\mathcal{F}]=
[l_{gA}(t\mathcal{G})]\otimes [uf_2\mathcal{F}]\in \text{Ind}_T^G U_\gamma\otimes W_\delta
\end{equation}

and this shows that $\text{Ind}_T^G U_\gamma\otimes W_\delta$ is
$(N^\mathcal{P}\cdot F^\text{diag})\ltimes G$-invariant. \\

It remains to show that the representations in \eqref{deccompaction} are irreducible,
pairwise inequivalent and that their sum is $L(X\times Y)$.

First of all, from

\[
\begin{split}
L(X\times Y)=&L(X)\otimes L(Y)\\
=&\left[\left(\bigoplus_{\omega\in\Omega}a_\omega V_\omega\right)\bigotimes\left(\bigoplus_{\delta\in\Delta_0}b_\delta W_\delta\right)\right]\bigoplus
\left[\left(\bigoplus_{\gamma\in\Gamma}c_\gamma\text{Ind}_T^GU_\gamma\right)\bigotimes\left(\bigoplus_{\delta\in\Delta\setminus\Delta_0}b_\delta W_\delta\right)\right]
\end{split}
\]

we deduce that \eqref{deccompaction}
is a decomposition into mutually orthogonal invariant subspaces (the invariance
has been proved above).
From \eqref{Wielandt}, it follows that

\[
\sum_{\omega\in\Omega}(a_\omega)^2=\lvert\mathcal{I}\rvert,\quad
\sum_{\delta\in\Delta}(b_\delta)^2=\lvert\mathcal{J}\rvert,\quad
\sum_{\gamma\in\Gamma}(c_\gamma)^2=\lvert\mathcal{I}_0\rvert,\quad
\text{and}\quad
\sum_{\delta\in\Delta_0}(b_\delta)^2=\lvert\mathcal{J}/\sim\rvert.
\]

The third equality follows from the fact that the $K$-orbits on $A_0\equiv T/K$ are given by \eqref{T-orbits}. The fourth equality follows from \eqref{deltauno} and \eqref{Sorbq}.
Then

\[
\sum_{\omega\in\Omega}\sum_{\delta\in\Delta_0}(a_\omega)^2(b_\delta)^2+
\sum_{\gamma\in\Gamma}\sum_{\delta\in\Delta\setminus\Delta_0}(c_\gamma)^2(b_\delta)^2=
\lvert\mathcal{I}\rvert\cdot\lvert\mathcal{J}/\sim\rvert+
+\lvert\mathcal{I}_0\rvert(\lvert\mathcal{J}\rvert-\lvert\mathcal{J}/\sim\rvert).
\]

On the other hand, from Theorem \ref{thdecorbcr} it follows that the number of orbits
of the stabilizer $R$ on $X\times Y$ is equal to

\[
\lvert\mathcal{I}_0\rvert\cdot\lvert\mathcal{J}\rvert+
\lvert\mathcal{I}\setminus\mathcal{I}_0\rvert\cdot\lvert\mathcal{J}/\sim\rvert=
\lvert\mathcal{I}\rvert\cdot\lvert\mathcal{J}/\sim\rvert+
+\lvert\mathcal{I}_0\rvert(\lvert\mathcal{J}\rvert-\lvert\mathcal{J}/\sim\rvert).
\]

Then from Lemma \ref{Wieapplic} we deduce that \eqref{deccompaction} is the isotypic
decomposition of $L(X\times Y)$ with respect to the action of the crested product.
\end{proof}

\begin{example}\label{wreathdec}{\rm
Now we apply Theorem \ref{deccompactionth} to the case of a wreath product. We use the notations in Example \ref{examplewreath}. Suppose that $L(Y)=\bigoplus_{\delta\in\Delta}b_\delta W_\delta$ and $L(Y')=\bigoplus_{\delta\in\Delta'}b'_\delta W'_\delta$ are the isotypic decompositions into $F$ and $F'$-irreducible representations. Suppose also that $0\in\Delta,\Delta'$ and that $W_0,W'_0$ are the trivial representations. If we take the equality partition $\mathcal{E}_Y$ on $Y$ and the universal partition $\mathcal{U}_{Y'}$ on $Y'$, then the resulting crested product of $F'$ by $F$ is isomorphic to the ordinary wreath product $F'\wr F$ (cf. Example \ref{crossnested}). Moreover, now \eqref{indlazero} (for $Y$) and \eqref{deltauno} (for $Y'$) become respectively $L(Y)\equiv L(Y)$ and $L(\mathcal{U}_{Y'})=W'_0$.
Then \eqref{deccompaction} yields the isotypic decomposition of $L(Y\times Y')$ as a homogeneous $F'\wr F$-space:

\[
L(Y\times Y')=\left[\bigoplus_{\delta\in\Delta}b_\delta(W_\delta \otimes W'_0)\right]\bigoplus \left[\bigoplus_{\delta\in\Delta'\setminus\{0\}}b'_\delta(L(Y)\otimes W'_\delta)\right].
\]
}
\end{example}

\begin{example}\label{iteratedwreath2}{\rm
In this example, we apply Theorem \ref{deccompactionth} to the setting of Example \ref{iteratedwreath}. We also assume the notation and the results in Example \ref{wreathdec}. Moreover, we suppose that $L(X)=\bigoplus_{\omega\in\Omega}a_\omega V_\omega$ and $L(X')=\bigoplus_{\omega\in\Omega'}a'_\omega V'_\omega$ are the isotypic decompositions into $G$ and $G'$ representations and therefore (by Example \ref{wreathdec})

\[L(X\times X')=\left[\bigoplus_{\omega\in\Omega}a_\omega(V_\omega \otimes V'_0)\right]\bigoplus \left[\bigoplus_{\omega\in\Omega'\setminus\{0\}}a'_\omega(L(X)\otimes V'_\omega)\right]
\]

is the isotypic decomposition of $L(X\times X')$ under the action of $G'\wr G$.
Now we have $A_{x_0}\equiv X'$ and if $\tilde{V}'_\omega\subseteq L(A_{x_0})$ denotes the subspace isomorphic to $V'_\omega\subseteq L(X')$, then $L(A_{x_0})=\bigoplus_{\omega\in\Omega'} a'_\omega \tilde{V}'_\omega$ is the isotypic decomposition of $L(A_{x_0})$ under the action of $G'\wr K$. Then, for $L(X\times X')$ the decomposition \eqref{indlazero} is:

\[
L(X\times X')=\text{Ind}_{G'\wr K}^{G'\wr G}L(A_{x_0})=\bigoplus_{\omega\in\Omega'}a'_\omega\text{Ind}_{G'\wr K}^{G'\wr G}
\tilde{V}'_\omega=
\bigoplus_{\omega\in\Omega'}a'_\omega (L(X)\otimes V'_\omega).
\]

On the other hand, for $L(Y\times Y')$ the decomposition \eqref{deltauno} is equivalent to $L(Y)\cong\bigoplus_{\delta\in\Delta}b_\delta(W_\delta \otimes W'_0)$. Then the isotypic decomposition of $L(X\times X'\times Y\times Y')$ under the action of the crested product $\left[\left((F')^\mathcal{P}\cdot(F'\wr F)^\text{diag}\right)\right]\ltimes (G'\wr G)$ is given by:

\[
\begin{split}
L(X\times X'\times Y\times Y')=&\left[\bigoplus_{\omega\in\Omega}\bigoplus_{\delta\in\Delta}a_\omega b_\delta(V_\omega \otimes V'_0\otimes W_\delta \otimes W'_0)\right]\\
&\bigoplus \left[\bigoplus_{\omega\in\Omega'\setminus\{0\}}\bigoplus_{\delta\in\Delta}a'_\omega b_\delta(L(X)\otimes V'_\omega \otimes W_\delta \otimes W'_0)\right]\\
&\bigoplus \left[\bigoplus_{\omega\in\Omega'}\bigoplus_{\delta\in\Delta'\setminus\{0\}}a'_\omega
b'_\delta(L(X)\otimes V'_\omega\otimes L(Y)\otimes W'_\delta)\right].
\end{split}
\]
}
\end{example}

In the last part of this section,
we give the formulas for the spherical matrix coefficients of the irreducible representations in Theorem \ref{deccompaction} (see subsection \ref{sphefunctsection}). Suppose that

\begin{itemize}

\item  for each $\omega\in\Omega$, $v_1^\omega,v_2^\omega,\dotsc,v_{a_\omega}^\omega$ is an orthonormal basis for the $K$-invariant vectors in $V_\omega$, and  $\phi^\omega_{i,i'},i,i'=1,2,\dotsc,a_\omega$ are the corresponding matrix coefficients: $\phi^\omega_{i,i'}(gx_0)=\langle v_i^\omega,g v_{i'}^\omega \rangle$, for any $g\in G$;

\item for each $\delta\in\Delta$, $w_1^\delta,w_2^\delta,\dotsc,w_{b_\delta}^\delta$ is an orthonormal basis for the $H$-invariant vectors in $W_\delta$, and $\psi^\omega_{j,j'},j,j'=1,2,\dotsc,b_\delta$ are the corresponding matrix coefficients: $\psi^\delta_{j,j'}(uy_0)=\langle w_j^\delta,u w_{j'}^\delta \rangle$, for any $u\in F$;

\item for each $\gamma\in\Gamma$, $u_1^\gamma,u_2^\gamma,\dotsc,u_{c_\gamma}^\gamma$ is an orthonormal basis for the $K$-invariant vectors in $U_\gamma$, and $\theta^\gamma_{h,h'},h,h'=1,2,\dotsc,c_\gamma$ are the corresponding matrix coefficients: $\theta^\gamma_{h,h'}(tx_0)=\langle u_h^\gamma,t u_{h'}^\gamma \rangle$, for any $t\in T$.

\end{itemize}

Note that the $\phi_{i,i'}^\omega$'s, the $\psi_{j,j'}^\delta$'s and the $\theta_{h,h'}^\gamma$'s are seen as functions defined respectively on $X$, $Y$ and $A_0$. Moreover, in the notation of Lemma \ref{pringen}, we set $\tilde{u}_h^\gamma=l_{A_0}u_h^\gamma\in\text{Ind}_T^G U_\gamma$, that is $\tilde{u}_h^\gamma$ is the copy of $u_h^\gamma$ in the subspace $l_{A_0}U_\gamma$ of $\text{Ind}_T^G U_\gamma=\bigoplus_{A\in\mathcal{P}}l_AU_\gamma$. In other words, if we think of $U_\gamma$ as a subspace of $L(A_0)$, then $\tilde{u}_h^\gamma=u_h^\gamma$ on $A_0$ and $\tilde{u}_h^\gamma\equiv 0$ on $X\setminus A_0$. In the same spirit, we set $\tilde{\theta}_{h,h'}^\gamma(x)=\theta_{h,h'}^\gamma(x)$ if $x\in A_0$,
$\tilde{\theta}_{h,h'}^\gamma(x)=0$ if $x\in X\setminus A_0$,

\begin{theorem}

\begin{enumerate}

\item
For $\omega\in\Omega$ and $\delta\in\Delta_0$, the set $v_i^\omega\otimes w_j^\delta$, $i=1,\dotsc,a_\omega$, $j=1,\dotsc,b_\delta$, is an orthonormal basis for the $R$-invariant vectors in the irreducible representation $V_\omega\otimes W_\delta$. Moreover, the corresponding spherical matrix coefficients are: $\phi^\omega_{i,i'}(x)\psi^\delta_{j,j'}(y)$, as functions of $(x,y)\in X\times Y$.

\item
For $\gamma\in\Gamma$ and $\delta\in\Delta\setminus\Delta_0$, the set $\tilde{u}^\gamma_h\otimes w_j^\delta$, $h=1,\dotsc,c_\gamma$, $j=1,\dotsc,b_\delta$, is an orthonormal basis for the $R$-invariant vectors in the irreducible representation $(\text{\rm Ind}_T^G U_\gamma)\otimes W_\delta$. Moreover, the corresponding spherical matrix coefficients are: $\tilde{\theta}^\gamma_{h,h'}(x)\psi_{j,j'}^\delta(y)$.

\end{enumerate}

\end{theorem}

\begin{proof}
\begin{enumerate}
\item If $(f_1f_2,k)\in R$, with $f_1\in N^\mathcal{P}$, $f_2\in H^\text{diag}$ and $k\in K$, then by \eqref{actensorprod2} we have $(f_1f_2,k)(v_i^\omega\otimes w_j^\delta)=(kv_i^\omega)\otimes (f_2w_j^\delta)=v_i^\omega\otimes w_j^\delta$, that is $v_i^\omega\otimes w_j^\delta$ is $R$-invariant. The set of all vectors of this type form an orthonormal basis for the $R$-invariant vectors in $V_\omega\otimes W_\delta$ because their number is equal to the multiplicity of this representation in $L(X\times Y)$. Using again \eqref{actensorprod2}, we can compute the spherical matrix coefficients: if $(x,y)=(f_1f_2,g)(x_0,y_0)\equiv (gx_0,f_1(gx_0)f_2y_0)$, with $(f_1f_2,g)\in (N^{\mathcal{P}} \cdot F^\text{diag})\ltimes G$, then

\[
\begin{split}
\langle v_i^\omega\otimes w_j^\delta, (f_1f_2,g) v_{i'}^\omega\otimes w_{j'}^\delta\rangle=&
\langle v_i^\omega\otimes w_j^\delta, (gv_{i'}^\omega)\otimes (f_2w_{j'}^\delta)\rangle\\
=&\langle v_i^\omega, gv_{i'}^\omega\rangle\langle w_j^\omega, f_2w_{j'}^\delta\rangle\\
=&\phi_{i,i'}^\omega(x)\cdot\psi_{j,j'}^\delta(y).
\end{split}
\]

\item
Now for $(f_1f_2,k)\in R$, \eqref{actensorprod3} yields
$(f_1f_2,k)(\tilde{u}_h^\gamma\otimes w_j^\delta))=(k\tilde{u}^\gamma_h)\otimes(f_2w_j^\delta)=\tilde{u}^\gamma_h\otimes w_j^\delta$, because $k\in T$, $f_1\equiv 1_F$ on $A_0$ and $\tilde{u}^\gamma_h$ is $K$-invariant. Now suppose that $(x,y)=(f_1f_2,g)(x_0,y_0)\equiv (gx_0,f_1(gx_0)f_2y_0)$ and $g\equiv gl_{A_0}=l_{gA_0}t$ with $t\in T$. If $gA_0\neq A_0$ then $l_{gA_0}t\tilde{u}^\gamma_{k'}=0$, while if $gA_0=A_0$ then $l_{gA_0}t=t$.
Therefore, again by \eqref{actensorprod3}, we have ($r$ is the value of $f_1$ on $gA_0$):

\[
\begin{split}
\langle \tilde{u}_k^\gamma\otimes w_j^\delta, (f_1f_2,g) \tilde{u}_{k'}^\gamma\otimes w_{j'}^\delta\rangle=&\langle \tilde{u}_k^\gamma\otimes w_j^\delta, (l_{gA_0}t\tilde{u}_{k'}^\gamma) \otimes[rf_2w_{j'}^\delta]\rangle\\
=&\langle \tilde{u}_k^\gamma,l_{gA_0}t\tilde{u}_{k'}^\gamma\rangle \langle w_j^\delta, rf_2w_{j'}^\delta\rangle\\
=&\tilde{\theta}_{k,k'}(x)\cdot\psi_{j,j'}^\delta(y).
\end{split}
\]

\end{enumerate}

\end{proof}

\begin{remark}{\rm
In \cite{Ba-Ca}, section 9, Bailey and Cameron describe more general notion of crested product involving a set of invariant partitions. It should be interesting to extend their theory to our noncommutative setting. Moreover, in section 11 they also suggest, as an open problem, to develop a theory of {\it generalized crested product} along the lines of the theory of generalized wreath product of groups and association schemes \cite{Ba2,BPRS}. It is also an interesting open problem to develop such theory for spaces with multiplicity. Our Examples \ref{iteratedwreath} and \ref{iteratedwreath2} also deserve to be generalized. The most intriguing aspect of the whole theory is the reciprocity between the decompositions \ref{decorbcr} and \ref{deccompaction}.}
\end{remark}

%%%%%%%%%%%%%%%%%%%%%%%%%%%%%%%%%%%%%%%%%%%%%%%%%%%%%%%%%%%%%%%%%%%
%%%%%%%%%%%%%%%%%%%%%%%%%%%%%%%%%%%%%%%%%%%%%%%%%%%%%%%%%%%%%%%%%%%
\section{Harmonic analysis of exponentiation and wreath product of permutations representations}
%%%%%%%%%%%%%%%%%%%%%%%%%%%%%%%%%%%%%%%%%%%%%%%%%%%%%%%%%%%%%%%%%%%
%%%%%%%%%%%%%%%%%%%%%%%%%%%%%%%%%%%%%%%%%%%%%%%%%%%%%%%%%%%%%%%%%%%

In this section, we want to obtain an explicit decomposition of the
exponentiation action of a wreath product. This is motivated by the
classical Hamming scheme. Actually, we analyze a more  general
notion, suggested by our recent work on finite lamplighter random
walks.

%%%%%%%%%%%%%%%%%%%%%%%%%%%%%%%%%%%%%%%%%%%%%%%%%%%%%%%%%%%%%%%%%%%
%%%%%%%%%%%%%%%%%%%%%%%%%%%%%%%%%%%%%%%%%%%%%%%%%%%%%%%%%%%%%%%%%%%
\subsection{Representation theory of wreath products of finite groups}\label{repwreath}
%%%%%%%%%%%%%%%%%%%%%%%%%%%%%%%%%%%%%%%%%%%%%%%%%%%%%%%%%%%%%%%%%%%
%%%%%%%%%%%%%%%%%%%%%%%%%%%%%%%%%%%%%%%%%%%%%%%%%%%%%%%%%%%%%%%%%%%

Suppose that $G$ and $F$ are finite groups and that $G$ acts transitively on a set $X$.
In this subsection, we give a description of the irreducible representations of the wreath product $F\wr G\equiv F^X\rtimes G$. We refer to \cite{JK,Hu} for complete proofs.\\

Every irreducible representation of the base group
$F^X$ may be written as a tensor product in the form

\[
\bigotimes_{x \in X}\sigma_x
\]

where

\[
\begin{array}{lcl}
X & \to & \widehat{F}\\
x&\mapsto& \sigma_x
\end{array}
\]

is any map from $X$ to $\widehat{F}$, the dual of $F$.
In other words, if $f_0\in F^X$ then

\[
\left(\bigotimes_{x\in X}\sigma_x\right)(f_0,1_G) =\bigotimes_{x\in X}\sigma_x(f_0(x))
\]

and if $\bigotimes\limits_{x\in X}v_x \in \bigotimes\limits_{x\in X}V_{\sigma_x}$, with $V_{\sigma_x}$ the space on which acts
the representation $\sigma_x$, then

\[
\left[\left(\bigotimes\limits_{x\in X}\sigma_x\right)(f_0,1_G)\right]\left(\bigotimes\limits_{x\in X}v_x\right) =
\bigotimes\limits_{x\in X}\sigma_x(f_0(x))v_x.
\]

The group $F\wr G$ acts on $\widehat{F^X}$ by the conjugation action:
the $(f,g)$ conjugate of $\bigotimes\limits_{x\in X}\sigma_x$ is defined by setting

\[
\stackrel{(f,g)\quad\qquad\qquad}{\left(\bigotimes\limits_{x\in X}\sigma_x\right)}\!\!\!\!\!\!\!(f_0,1_G)  =
\left(\bigotimes\limits_{x\in X}\sigma_x\right)[(f,g)^{-1}(f_0,1_G)(f,g)].
\]

Then we have:

\begin{equation}\label{conjsigma}
\stackrel{(f,g)\quad\qquad\qquad}{\left(\bigotimes\limits_{x\in X}\sigma_x\right)}\!\!\!\!\!\! = \bigotimes\limits_{x\in X} {^{f(x)}\!\sigma_{g^{-1}x}}
\sim\bigotimes\limits_{x\in X}\sigma_{g^{-1}x}.
\end{equation}

The {\em in\ae rtia group} $I_{F \wr G}(\sigma)$ of $\sigma = \left(\bigotimes\limits_{x\in X}\sigma_x\right) \in \widehat{F^X}$ is the stabilizer of $\sigma$ with respect to the conjugation action; \eqref{conjsigma} ensures us that

\[
I_{F \wr G}(\sigma) = F \wr T_G(\sigma)\cong F^X \ltimes T_G(\sigma) ,
\]

where $T_G(\sigma) = \{g \in G: \sigma_{gx} \sim \sigma_x \ \forall x \in X\}$.\\

We give two general definitions. Suppose that $H$ is a subgroup of $G$ and that $\theta$ is an $H$ representation. An {\em extension} of $\theta$ to $G$ is a representation $\widetilde{\theta}$ of $G$ such that: $\text{Res}^G_H\widetilde{\theta}=\theta$. Clearly, if $\theta$ is irreducible, $\widetilde{\theta}$ is irreducible too. On the other hand, if $N$ is a normal subgroup of $G$ and $\eta$ is a representation of the quotient group $G/N$, then the {\em inflation} $\overline{\eta}$ of $\eta$ to $G$ is defined by setting $\overline{\eta}(g)=\eta(gN)$, for all $g\in G$. In other words, we compose $\eta$ (which is a homomorphism of $G/N$ into the unitary group of a Hermitian space) with the quotient homomorphism: $G\to G/N$. Clearly, $\eta$ is irreducible if and only if $\overline{\eta}$ is irreducible.\\

In our setting, each $\left(\bigotimes\limits_{x\in X}\sigma_x\right) \in \widehat{F^X}$ has an extension $\widetilde{\sigma}$ to
the whole $I_{F \wr G}(\sigma)$: it is given by setting

\begin{equation}\label{extsigma}
\widetilde{\sigma}(f,g)\left(\bigotimes\limits_{x\in X}v_x\right):=\bigotimes\limits_{x\in X}\sigma_{g^{-1}x}(f(x))v_{g^{-1}x}\equiv\bigotimes\limits_{x\in X}\sigma_x(f(x))v_{g^{-1}x},,
\end{equation}

for all $(f,g)\in F\wr T_G(\sigma)$ and $\bigotimes\limits_{x\in X}v_x\in\bigotimes\limits_{x\in X}V_{\sigma_x}$.

Now let $\Sigma$ be a system of representatives  for the $F\wr G-$conjugacy classes of irreducible
representations of $\widehat{F^X}$. For each $\sigma \in \Sigma$, denote by $\widetilde{\sigma}$ its extension
to $I_{F\wr G}(\sigma)$ as shown in \eqref{extsigma}.
For each $\eta \in \widehat{T}_G(\sigma)$, denote by $\overline{\eta}$ its inflation to $I_{F\wr G}(\sigma)$ (using the homomorphism $I_{F\wr G}(\sigma)\to T_G(\sigma)\cong I_{F\wr G}(\sigma)/F^X$). That is, if $U$ is the representation space of $\eta$, then

\begin{equation}\label{etasegn}
\overline{\eta}(f,g)u=\eta(g)u\qquad\forall (f,g)\in I_{F\wr G}(\sigma),u\in U.
\end{equation}

We are in position to enunciate the main theorem in the representation theory of wreath products \cite{JK,Hu}.

\begin{theorem}\label{irrwreathprod}
The dual of $F\wr G$ is given by:

\[
\widehat{F\wr G} = \{\text{\rm Ind}_{I_{F\wr G}(\sigma)}^G(\widetilde{\sigma}\otimes \overline{\eta}):
\sigma \in \Sigma, \eta \in \widehat{T}_G(\sigma)\},
\]

that is the above is the list of all irreducible representations of $F \wr G$, and for
different values of $\sigma, \eta$ we obtain inequivalent representations.
\end{theorem}

%%%%%%%%%%%%%%%%%%%%%%%%%%%%%%%%%%%%%%%%%%%%%%%%%%%%%%%%%%%%%%%%%%%
%%%%%%%%%%%%%%%%%%%%%%%%%%%%%%%%%%%%%%%%%%%%%%%%%%%%%%%%%%%%%%%%%%%
\subsection{Exponentiations and wreath products}\label{expwr}
%%%%%%%%%%%%%%%%%%%%%%%%%%%%%%%%%%%%%%%%%%%%%%%%%%%%%%%%%%%%%%%%%%%
Let $G,F,X$ be as in the preceding subsection and suppose also that $G$ acts transitively on another set $Z$ and that $F$ acts transitively on $Y$. We form the wreath product $F\wr G$ with respect to the action of $G$ on $X$. The group $F\wr G$ acts on $Y^X$ via the {\em exponentiation} of the action of $F$ on $Y$: if $\varphi\in Y^X$, $(f,g)\in F\wr G$ then $(f,g)\varphi$ is defined by setting $[(f,g)\varphi](x)=f(x)\varphi(g^{-1}x)$, for every $x\in X$.
It also acts on $Z$ via the {\em infation} of the action of $G$ on $Z$: if $(f,g)\in F\wr G$ and $z\in Z$ then $(f,g)z=gz$.
Then we define the {\em wreath product} of the action of $F$ on $Y$ by the actions of $G$ on $X$ and $Z$ as the direct product of the exponentiation and the inflation:

\[
(f,g)(\varphi,z)=((f,g)\varphi,gz)
\]

for all $(f,g)\in F\wr G$ and $(\varphi,z)\in Y^X\times Z$. Note that this action is transitive, simply because $F^X$ is transitive on $Y^X$ and $G$ is transitive on $Z$.
The task of this section is to obtain an explicit decomposition of the permutation representation of $F\wr G$ on $Y^X\times Z$. The study of this permutation representation is motivated by our work on finite lamplighter random walks \cite{ScaTol2} (and for $Z$ trivial, our results apply to the exponentiation action).\\

We will denote by $\lambda$ the permutation representation of $F^X$ on $L(Y^X)$ and by $\widetilde{\lambda}$ its extension to $F\wr G$ (that is, the permutation representation associated to the exponentiation). In particular, for $g\in G$, $\psi\in L(Y^X)$  and $\varphi\in Y^X$, we have $[\widetilde{\lambda}(\mathbf{1}_F,g)\psi](\varphi)=\psi(g^{-1}\varphi)$.
If $\psi_x\in L(Y)$ for all $x\in X$, then  $\bigotimes\limits_{x\in X}\psi_x\in L(Y^X)$ is defined by setting

\[
\biggl(\bigotimes_{x\in X}\psi_x\biggr)(\varphi)= \prod_{x\in X}\psi_x(\varphi(x)),\qquad \forall \varphi\in Y^X.
\]

Moreover, if $\psi\in L(Y)$ and $f\in F$, we set $(f\psi)(y)=\psi(f^{-1}y)$ for all $y\in Y$, that is in this way we denote the permutation representation of $F$ on $Y$.

The following lemma is similar to \eqref{extsigma}, but now $g\in G$. We give the elementary proof for completeness.

\begin{lemma}\label{actionpsi}
If $\bigotimes\limits_{x\in X}\psi_x\in L(Y^X)$, with $\psi_x\in L(Y)$ for all $x\in X$, and $(f,g)\in F\wr G$ then

\[
\widetilde{\lambda}(f,g)\biggl(\bigotimes\limits_{x\in X}\psi_x\biggr)=\bigotimes\limits_{x\in X}f(x)\psi_{g^{-1}x}.
\]

\end{lemma}

\begin{proof}
For any $\varphi\in Y^X$, we have:

\[
\begin{split}
\bigg[[\widetilde{\lambda}(f,g)\biggl(\bigotimes\limits_{x\in X}\psi_x\biggr)\biggr](\varphi)
=&\biggl(\bigotimes\limits_{x\in X}\psi_{x}\biggr)[(f,g)^{-1}\varphi]\\
=&\prod\limits_{x\in X}\psi_{x}[f(gx)^{-1}\varphi(gx)]\\
=&\prod\limits_{x\in X}\psi_{g^{-1}x}[f(x)^{-1}\varphi(x)]\\
=&\prod\limits_{x\in X}\left[f(x)\psi_{g^{-1}x}\right](\varphi(x))\\
=&\biggl[\bigotimes\limits_{x\in X}f(x)\psi_{g^{-1}x}\biggr](\varphi).\\
\end{split}
\]

\end{proof}

Suppose that $\Sigma$ is as in Theorem \ref{irrwreathprod}. Fix a $\sigma\in\Sigma$ and set $I=I_G(\sigma)$. We also suppose that each $\sigma_x, x\in X$, appears in the decomposition of $L(Y)$ into irreducible $F$-representations. Clearly, there exist a partition $X=\coprod\limits_{i=1}^n\Omega_i$ of $X$ and
$\sigma_1,\sigma_2,\dotsc,\sigma_n$ irreducible, pairwise inequivalent $F$-representations  such that: $\sigma_x=\sigma_i$ for all $x\in \Omega_i$, $i=1,\dotsc,n$. Moreover, $I=\{g\in G:g\Omega_i=\Omega_i, i=1,2,\dotsc,n\}$. Let $m_i$ be the multiplicity of $\sigma_i$ in the permutation representation of $F$ on $Y$. If $x\in\Omega_i$ and $V_{\sigma_x}\equiv V_i$ is the representation space of $\sigma_x\equiv \sigma_i$, we fix an orthonormal basis $T_{x,1},T_{x,2},\dotsc,T_{x,m_i}$ in $\text{Hom}_F(V_{\sigma_x},L(Y))$ (orthonormal with respect to the Hilbert-Schmidt scalar product); we suppose that, for every $1\leq h\leq m_i$, the operator $T_{x,h}$ is the same for all $x\in\Omega_i$. Denote by $J$ the space of all functions $j:X\to\mathbb{N}$ such that $j(x)\in\{1,2,\dotsc,m_i\}$ for all $x\in\Omega_i$, $i=1,2,\dotsc,n$. The group $I$ acts on $J$ in a natural way: if $g\in I$ and $j\in J$, then $gj$ is defined by setting: $(gj)(x)=j(g^{-1}x)$ (recall that $I$ stabilizes eve!
 ry $\Omega_i$). For any $j\in J$, set $T_j=\bigotimes\limits_{x\in X}T_{x,j(x)}$. Then the set $\{T_j:j\in J\}$ is an orthonormal basis for $\text{Hom}_{F^X}\left(\bigotimes\limits_{x\in X}V_{\sigma_x},L(Y^X)\right)$.
Note also that

\begin{equation}\label{Tx}
T_{gx,j(gx)}=T_{x,j(gx)}\qquad\text{if}\qquad g\in I
\end{equation}

because $T_{x,h}$ does not depend on $x\in\Omega_i$.

\begin{lemma}
For $g\in I$ and $T\in\text{\rm Hom}_{F^X}\left(\bigotimes\limits_{x\in X}V_{\sigma_x},L(Y^X)\right)$, define a linear operator $\pi(g)T$ by setting

\[
\pi(g)T=\widetilde{\lambda}(\mathbf{1}_F,g)T\widetilde{\sigma}(\mathbf{1}_F,g^{-1}).
\]

Then $\pi$ is a representation of $I$ on $\text{\rm Hom}_{F^X}\left(\bigotimes\limits_{x\in X}V_{\sigma_x},L(Y^X)\right)$.
\end{lemma}

\begin{proof}

We show that $\widetilde{\lambda}(\mathbf{1}_F,g)T\widetilde{\sigma}(\mathbf{1}_F,g)^{-1}\in\text{Hom}_{F^X}\left(\bigotimes\limits_{x\in X}V_{\sigma_x},L(Y^X)\right)$ for every $T\in \text{Hom}_{F^X}\left(\bigotimes\limits_{x\in X}V_{\sigma_x},L(Y^X)\right)$ and $g\in I$.
Indeed,

\[
\begin{split}
\left[\widetilde{\lambda}(\mathbf{1}_F,g)T\widetilde{\sigma}(\mathbf{1}_F,g^{-1})\right]\widetilde{\sigma}(f,1_G)=&\widetilde{\lambda}(\mathbf{1}_F,g)T\widetilde{\sigma}(g^{-1}f,1_G)\widetilde{\sigma}(\mathbf{1}_F,g^{-1})\\
=&\widetilde{\lambda}(\mathbf{1}_F,g)\widetilde{\lambda}(g^{-1}f,1_G)T\widetilde{\sigma}(\mathbf{1}_F,g^{-1})\\
=&\widetilde{\lambda}(f,1_G)\left[\widetilde{\lambda}(\mathbf{1}_F,g)T\widetilde{\sigma}(\mathbf{1}_F,g^{-1})\right].
\end{split}
\]

It is clear that $\pi$ is a representation of $I$.
\end{proof}

\begin{lemma}\label{lambdaTsigma}
For $(f,g)\in F\wr I$ and $j\in J$, we have:

\[
\widetilde{\lambda}(f,g)T_j=T_{gj}\widetilde{\sigma}(f,g).
\]

\end{lemma}

\begin{proof}

From Lemma \ref{actionpsi} and \eqref{Tx}, it follows that

\[
\begin{split}
\left[\widetilde{\lambda}(f,g)T_j\right]\biggl(\bigotimes_{x\in X}v_x\biggr)=&
\widetilde{\lambda}(f,g)\biggl(\bigotimes_{x\in X}T_{x,j(x)}v_x\biggr)\\
=&\bigotimes_{x\in X}f(x)T_{x,j(g^{-1}x)}v_{g^{-1}x}.\\
\end{split}
\]

On the other hand, from \eqref{extsigma} and the fact that $T_{x,j(x)}\in\text{Hom}_F(V_{\sigma_x},L(Y))$, we get

\[
\begin{split}
\left[T_{gj}\widetilde{\sigma}(f,g)\right]\biggl(\bigotimes_{x\in X}v_x\biggr)=&
\biggl(\bigotimes_{x\in X}T_{x,(gj)(x)}\biggr)\biggl(\bigotimes_{x\in X}\sigma_x(f(x))v_{g^{-1}x}\biggr)\\
=&\bigotimes_{x\in X}T_{x,(gj)(x)}\sigma_x(f(x))v_{g^{-1}x}\\
=&\bigotimes_{x\in X}f(x)T_{x,j(g^{-1}x)}v_{g^{-1}x}.
\end{split}
\]
\end{proof}

\begin{corollary}\label{corpi}
The representation $\pi$ is equivalent to the permutation representation of $I$ on $J$.
\end{corollary}

\begin{proof}
As a particular case of Lemma \ref{lambdaTsigma}, we get the identity
$\pi(g)T_j\equiv\widetilde{\lambda}(\mathbf{1}_F,g)T_j\widetilde{\sigma}(\mathbf{1}_F,g^{-1})=T_{gj}$.
But $\{T_j:j\in J\}$ is an orthonormal basis for $\text{Hom}_{F^X}\left(\bigotimes\limits_{x\in X}\sigma_x,L(Y^X)\right)$, and therefore the map

\[
\begin{array}{rcl}
L(J)&\longrightarrow& \text{Hom}_{F^X}\left(\bigotimes\limits_{x\in X}\sigma_x,L(Y^X)\right)\\
\delta_j\quad&\longmapsto&\quad T_j
\end{array}
\]

is an isomorphism of $I$-representations.

\end{proof}

\begin{theorem}\label{decFwrI}
Suppose that $\mathcal{T}\in\text{\rm Hom}_I(U,L(Z\times J))$.
For
\[
v:=\bigotimes\limits_{x\in X}v_x\in \bigotimes\limits_{x\in X}V_{\sigma_x}\qquad \text{and}\qquad u\in U,
\]

define $\widehat{\mathcal{T}}(v\otimes u)\in L(Y^X\times Z)$ by setting:

\begin{equation}\label{That}
\left[\widehat{\mathcal{T}}\left(v\otimes u\right)\right](\varphi,z)=\left[\sum_{j\in J}(\mathcal{T}u)(z,j)\cdot T_jv\right](\varphi)\equiv \sum_{j\in J}(\mathcal{T}u)(z,j)\cdot (T_jv)(\varphi),
\end{equation}

for all $(\varphi,z)\in Y^X\times Z$.
Then $\widehat{\mathcal{T}}\in \text{\rm Hom}_{F\wr I}\left(\Bigl(\bigotimes\limits_{x\in X}V_{\sigma_x}\Bigr)\bigotimes U,L(Y^X\times Z)\right)$ and the map

\[
\begin{array}{rcl}
\text{\rm Hom}_I(U,L(Z\times J))&\longrightarrow&\text{\rm Hom}_{F\wr I}\left(\Bigl(\bigotimes\limits_{x\in X}V_{\sigma_x}\Bigr)\bigotimes U,L(Y^X\times Z)\right)
\\
\mathcal{T}\qquad&\longmapsto&\qquad\widehat{\mathcal{T}}
\end{array}
\]

is a linear isometric isomorphism.

\end{theorem}
\begin{proof}
First we show that
$\widehat{\mathcal{T}}\in \text{\rm Hom}_{F\wr I}\left(\Bigl(\bigotimes\limits_{x\in X}V_{\sigma_x}\Bigr)\bigotimes U,L(Y^X\times Z)\right)$. If $(f,g)\in F\wr I$ then

\begin{multline*}
\left\{\widehat{\mathcal{T}}\left[(\widetilde{\sigma}(f,g)\otimes\overline{\eta}(f,g))\left(v\otimes u\right)\right]\right\}(\varphi,z)
=\left\{\widehat{\mathcal{T}}\left[\widetilde{\sigma}(f,g)v\otimes \eta(g)u\right]\right\}(\varphi,z)\quad(\text{by}\;\eqref{etasegn})\\
=\left\{\sum_{j\in J}[(\mathcal{T}\eta(g)u)(z,j)]\cdot T_j\widetilde{\sigma}(f,g)v\right\}(\varphi)\quad(\text{by}\;\eqref{That})\\
=\left\{\widetilde{\lambda}(f,g)\left[\sum_{j\in J}[\mathcal{T}(u)(g^{-1}z,g^{-1}j)]\cdot T_{g^{-1}j}v\right]\right\}(\varphi)\\
(\text{because}\;\mathcal{T}\in\text{\rm Hom}_I(U,L(Z\times J)\;\text{and}\;\widetilde{\lambda}(f,g)T_{g^{-1}j}=T_j\widetilde{\sigma}(f,g)\;\text{by Lemma}\; \ref{lambdaTsigma})\\
=\left\{\sum_{j\in J}[\mathcal{T}(u)(g^{-1}z,j)]\cdot T_jv\right\}((f,g)^{-1}\varphi)\quad (\text{replacing}\;g^{-1}j\;\text{by}\;j)\\
=\left[\widehat{\mathcal{T}}(v\otimes u)\right]\left((f,g)^{-1}\varphi,(f,g)^{-1}z\right)\quad(\text{again by}\;\eqref{That}).
\end{multline*}
\qquad\\
\qquad\\
This proves that $\widehat{\mathcal{T}}$ commutes with $F\wr I$.\\

Now we define the inverse correspondence $\mathcal{T}\longrightarrow\widehat{\mathcal{T}}$.
Suppose that

\[
\mathcal{T}\in\text{\rm Hom}_{F\wr I}\left(\biggl(\bigotimes\limits_{x\in X}V_{\sigma_x}\biggr)\bigotimes U,L(Y^X\times Z)\right).
\]

For all fixed $u\in U$ and $z\in Z$, define a map

\[
\mathcal{T}_{u,z}^\sharp:\bigotimes_{x\in X}V_{\sigma_x}\longrightarrow L(Y^X)
\]

by setting

\begin{equation}\label{Tdiesis}
\left[\mathcal{T}_{u,z}^\sharp v\right](\varphi)=
\left\{\mathcal{T}\left[v\otimes u\right]\right\}(\varphi,z)
\end{equation}

for all $\bigotimes\limits_{x\in X}v_x\in \bigotimes\limits_{x\in X}V_{\sigma_x}$ and $\varphi\in Y^X$.
Then $\mathcal{T}_{u,z}^\sharp\in\text{\rm Hom}_{F^X}\left(\bigotimes\limits_{x\in X}V_{\sigma_x},L(Y^X)\right)$. Indeed, if $f\in F^X$ then $\overline{\eta}(f,1_G)u=u$ and therefore

\[
\begin{split}
\left[\mathcal{T}_{u,z}^\sharp\widetilde{\sigma}(f,1_G)v\right](\varphi)=&
\bigr\{\mathcal{T}\left[\left(\widetilde{\sigma}(f,1_G)v\right)\otimes u\right]\bigl\}(\varphi,z)\\
=&\left\{\mathcal{T}(\widetilde{\sigma}(f,1_G))\otimes\widetilde{\eta}(f,1_G))\left[v\otimes u\right]\right\}(\varphi,z)\\
=&\left\{\mathcal{T}\left[v\otimes u\right]\right\}((f,1_G)^{-1}\varphi,z)\\
=&\left[\mathcal{T}_{u,z}^\sharp v\right]((f,1_G)^{-1}\varphi).\\
\end{split}
\]

It follows that there exists $\alpha_{u,z}\in L(J)$ such that

\begin{equation}\label{alpha}
\mathcal{T}_{u,z}^\sharp=\sum_{j\in J}\alpha_{u,z}(j)T_j.
\end{equation}

Thus we can define a linear map $\widetilde{\mathcal{T}}:U\rightarrow L(Z\times J)$ by setting, for $u\in U$ and $z,j\in Z\times J$,

\begin{equation}\label{Ttilde}
[\widetilde{\mathcal{T}}u](z,j)=\alpha_{u,z}(j),
\end{equation}

where $\alpha_{u,z}$ is given by \eqref{alpha}. In other words, from \eqref{Tdiesis}, \eqref{alpha} and \eqref{Ttilde} it follows that

\begin{equation}\label{TTtilde}
\left\{\mathcal{T}\left[v\otimes u\right]\right\}(\varphi,z)=\left\{\sum_{j\in J}[\widetilde{\mathcal{T}}u](z,j)\cdot T_jv\right\}(\varphi).
\end{equation}

Now we prove that $\widetilde{\mathcal{T}}\in\text{\rm Hom}_I(U,L(Y^X,J))$: for any $g\in I$ we have

\begin{multline*}
\left\{\sum_{j\in J}[(\widetilde{\mathcal{T}}\eta(g)u)(z,j)]\cdot T_jv\right\}(\varphi)
=\left\{\mathcal{T}\left[\widetilde{\sigma}(\mathbf{1}_F,g)\widetilde{\sigma}(\mathbf{1}_F,g^{-1})v\otimes \eta(g)u\right]\right\}(\varphi,z)\\
=\left\{\mathcal{T}\left[\widetilde{\sigma}(\mathbf{1}_F,g^{-1})v\otimes u\right]\right\}\left((\mathbf{1}_F,g^{-1})\varphi,(\mathbf{1}_F,g)^{-1}z\right)\\
=\left\{\sum_{j\in J}[\widetilde{\mathcal{T}}u](g^{-1}z,j)\cdot T_j\widetilde{\sigma}(\mathbf{1}_F,g^{-1})v\right\}\left((\mathbf{1}_F,g^{-1})\varphi\right)\\
=\left\{\sum_{j\in J}[\widetilde{\mathcal{T}}u](g^{-1}z,j)\cdot T_{gj}v\right\}(\varphi)\qquad\;(\text{by Lemma}\;\ref{lambdaTsigma})\\
=\left\{\sum_{j\in J}[\widetilde{\mathcal{T}}u](g^{-1}z,g^{-1}j)\cdot T_jv\right\}(\varphi).\\
\end{multline*}

and therefore $\widetilde{\mathcal{T}}$ commutes with $I$.
From \eqref{That} and \eqref{TTtilde}, it is also clear that
$\widehat{\widetilde{\mathcal{T}}}=\mathcal{T}$ and that
$\widetilde{\widehat{\mathcal{T}}}=\mathcal{T}$.

Finally, we prove that $\mathcal{T}\mapsto \widehat{\mathcal{T}}$ is an isometry. Suppose that by $\mathcal{B}_1$ (resp. $\mathcal{B}_2$) is an orthonormal basis in $\bigotimes_{x\in X}V_{\sigma_x}$ (resp. $U$). If $\mathcal{T}_1,\mathcal{T}_2\in\text{Hom}_I(U,L(J\times Z))$, then

\[
\begin{split}
\langle \widehat{\mathcal{T}_1},\widehat{\mathcal{T}_2} \rangle_{HS}=&
\sum_{v\in\mathcal{B}_1}\sum_{u\in\mathcal{B}_2}\sum_{(\varphi,z)\in Y^X\times Z}
[\widehat{\mathcal{T}_1}(v\otimes u)](\varphi,z)\cdot \overline{[\widehat{\mathcal{T}_2}(v\otimes u)](\varphi,z)}\\
=&\sum_{v,u}\sum_{(\varphi,z)}\left[\sum_{j\in J}(\mathcal{T}_1u)(z,j)\cdot(T_jv)(\varphi)\right]\cdot\left[\sum_{i\in J}\overline{(\mathcal{T}_2u)(z,i)\cdot(T_iv)(\varphi)}\right]\\
=&\sum_{u,j,i,z}(\mathcal{T}_1u)(z,j)\cdot\overline{(\mathcal{T}_2u)(z,i)}\langle T_j,T_i\rangle_{HS}\\
&\sum_{u,j,z}(\mathcal{T}_1u)(z,j)\cdot\overline{(\mathcal{T}_2u)(z,j)}\\
&=\langle\mathcal{T}_1,\mathcal{T}_2\rangle_{HS}.
\end{split}
\]

\end{proof}

\begin{remark}
{\rm
Using the same techniques in the proof of Theorem \ref{decFwrI}, one can prove that the map $\mathcal{T}^\sharp: u\otimes\delta_z\mapsto\mathcal{T}_{u,z}^\sharp$ belongs to
$\text{\rm Hom}_I(U\otimes L(Z),L(J))$ and that

\[
\begin{array}{rcl}
\text{\rm Hom}_{F\wr I}\left(\Bigl(\bigotimes\limits_{x\in X}V_{\sigma_x}\Bigr)\bigotimes U,L(Y^X\times Z)\right)&\longrightarrow&
\text{\rm Hom}_I(U\otimes L(Z),L(J))
\\
\mathcal{T}\qquad&\longmapsto&\qquad\mathcal{T}^\sharp
\end{array}
\]

is a linear isomorphism. Recall also Corollary \ref{corpi}.

}
\end{remark}

Frobenius reciprocity, as stated in Proposition \ref{Frobenius}, yields an explicit isometric isomorphism

\[
\begin{array}{rcl}
\text{\rm Hom}_{F\wr I}\left(\Bigl(\bigotimes\limits_{x\in X}V_{\sigma_x}\Bigr)\bigotimes U,L(Y^X\times Z)\right)&\longrightarrow&\text{\rm Hom}_{F\wr G}\left(\text{\rm Ind}_{F\wr I}^{F\wr G}\Bigl[\Bigl(\bigotimes\limits_{x\in X}V_{\sigma_x}\Bigr)\bigotimes U\Bigr],L(Y^X\times Z)\right).
\\
\mathcal{T}\qquad&\longmapsto&\qquad\stackrel{\diamond}{\mathcal{T}}
\end{array}
\]

Then by combining Theorem \ref{decFwrI} and the isomorphism $\mathcal{T}\mapsto \stackrel{\diamond}{\mathcal{T}}$, we get the main result of this section.

\begin{theorem}\label{decFwrG}
The map

\[
\begin{array}{rcl}
\text{\rm Hom}_I(U,L(Z\times J))&\longrightarrow&\text{\rm Hom}_{F\wr G}\left(\text{\rm Ind}_{F\wr I}^{F\wr G}\Bigl[\Bigl(\bigotimes\limits_{x\in X}V_{\sigma_x}\Bigr)\bigotimes U\Bigr],L(Y^X\times Z)\right)
\\
\mathcal{T}\qquad&\longmapsto&\qquad\stackrel{\diamond}{\widehat{\mathcal{T}}}
\end{array}
\]
is a linear isometric isomorphism.

\end{theorem}

The following Corollary is an immediate consequence.

\begin{corollary}\label{cormolt}
The multiplicity of the irreducible representation $\text{\rm Ind}_{F\wr I}^{F\wr G}(\widetilde{\sigma}\otimes \overline{\eta})$ in the decomposition of $L(Y^X\times Z)$ into irreducible $F\wr G$-representations is equal to the multiplicity of $\eta$ in the decomposition of $L(J\times Z)$ into irreducible $I$-representations.
\end{corollary}

\begin{remark}{\rm
If $Y\equiv F$ and $Z\equiv G$, both with the left regular
representation, then Theorem \ref{decFwrG} leads to a decomposition
of the left regular representation of $F\wr G$ and therefore
Corollary \ref{cormolt} yields a formula for the dimension of the
irreducible representation $\text{\rm Ind}_{F\wr I}^{F\wr
G}(\widetilde{\sigma}\otimes \overline{\eta})$: such a dimension is
equal to the multiplicity of $\eta$ into the decomposition of
$L(J\times G)$ into irreducible $I$-representations. We want to show
that this fact agree with usual formula for the dimension of an
induced representation. The action of $I$ on $J\times G$ has $\lvert
J\rvert\cdot\frac{\lvert G\rvert}{\lvert I\rvert}$ orbits, and each
orbit is equivalent to the left action of $I$ on itself: the
stabilizer of any point $(j,g)\in J\times G$ is the trivial
subgroup. Then the above computed multiplicity is equal to
$\dim\eta\cdot\lvert J\rvert\cdot\frac{\lvert G\rvert}{\lvert
I\rvert}$. But the dimension of $\sigma$ is just $\lvert J\rvert$
(recall that $J$ is obtained by decomposing the left regular
representation of $F^X$). Then the formula for the dimension of an
induced representation tells us that in fact $\dim \text{\rm
Ind}_{F\wr I}^{F\wr G}(\widetilde{\sigma}\otimes
\overline{\eta})=\dim\eta\cdot\dim\sigma\cdot\frac{\lvert
G\rvert}{\lvert I\rvert}$. }
\end{remark}

%%%%%%%%%%%%%%%%%%%%%%%%%%%%%%%%%%%%%%%%%%%%%%%%%%%%%%%%%%%%%%%%%%%%%%%%%%
%%%%%%%%%%%%%%%%%%%%%%%%%%%%%%%%%%%%%%%%%%%%%%%%%%%%%%%%%%%%%%%%%%%%%%%%%
\subsection{The case $G=C_2$ and $Z$ trivial}
%%%%%%%%%%%%%%%%%%%%%%%%%%%%%%%%%%%%%%%%%%%%%%%%%%%%%%%%%%%%%%%%%%%%%%%%%%
%%%%%%%%%%%%%%%%%%%%%%%%%%%%%%%%%%%%%%%%%%%%%%%%%%%%%%%%%%%%%%%%%%%%%%%%%

As an example, in this section we examine the case $G=C_2\equiv X$ and $Z$ trivial. We set $C_2=\{1,-1\}$ and $U=$ the representation space of the nontrivial (the alternating) representation of $C_2$.  Suppose that $L(Y)=\bigotimes\limits_{\sigma\in R}m_\sigma V_\sigma$ is the isotypic decomposition of $L(Y)$ into irreducible $F$-representations.
More precisely, we suppose that $T^\sigma_1,T^\sigma_2,\dotsc,T^\sigma_{m_\sigma}$ is an orthonormal basis for $\text{Hom}_F(V_\sigma,L(Y))$.
Then

\begin{equation}\label{LYdec}
L(Y\times Y)=\bigoplus\limits_{\sigma,\sigma'\in R}m_{\sigma}m_{\sigma'} (V_\sigma\otimes V_{\sigma'})\equiv\bigoplus\limits_{\sigma,\sigma'\in R}\bigoplus_{j=1}^{m_\sigma}\bigoplus_{j'=1}^{m_{\sigma'}}\left(T^\sigma_jV_\sigma\oplus T^{\sigma'}_{j'}V_{\sigma'}\right)
\end{equation}

is a decomposition into irreducible $F\times F$-representations. We want to show how to obtain a decomposition into irreducible $F\wr C_2$-representations.\\

When $\sigma\neq \sigma'$, the in\ae rtia group of $V_\sigma\otimes
V_{\sigma'}$ is trivial. Moreover, just from the definition of
induced representation (or from Proposition \eqref{Frobenius}) it
follows that

\[
\text{Ind}_{F\times F}^{F\wr C_2}\left(V_\sigma\otimes
V_{\sigma'}\right)=(V_\sigma\otimes V_{\sigma'})\oplus( V_{\sigma'}\otimes V_\sigma).
\]

The multiplicity of this irreducible $F\wr C_2$-representation into $L(Y\times Y)$ is clearly $m_\sigma m_{\sigma'}$. With more cumbersome but more precise notation, we can define

\[
W_{\sigma,\sigma'}^{j,j'}=\left(T^\sigma_jV_\sigma\otimes T^{\sigma'}_{j'}V_{\sigma'}\right)\oplus\left(T^{\sigma'}_{j'}V_{\sigma'}\otimes T^\sigma_jV_\sigma\right)
\]

which are mutually orthogonal subspaces of $L(Y\times Y)$ isomorphic to $V_\sigma\otimes V_{\sigma'}$.\\

When $\sigma=\sigma'$, the in\ae rtia group of $V_\sigma\otimes V_{\sigma'}$ is $F\wr C_2$. Then there is no induction and we only need to apply Theorem \ref{decFwrI}. We have $J=\{(i,j):1\leq i,j\leq m_\sigma\}$ and the orbits of $C_2$ on $J$ are $\{(i,j),(j,i)\},1\leq i\neq j\leq m_\sigma$ and $\{(i,i)\},i=1,\dotsc,m_\sigma$. This means that $L(J)$ contains $\frac{m_\sigma(m_\sigma+1)}{2}$ times the trivial representation of $C_2$ and $\frac{m_\sigma(m_\sigma-1)}{2}$ times the nontrivial representation $U$.
For $i,j=1,2,\dotsc,m_\sigma$, we can define the subspaces

\[
W^{i,j}_{\sigma,+}=\left\langle \left(T_i^\sigma v_1\otimes T_j^\sigma v_2\right)+
\left(T_j^\sigma v_1\otimes T_i^\sigma v_2\right)
:v_1,v_2\in V_\sigma\right\rangle
\]

that corresponds to the choice $\eta$= the trivial representation in \eqref{That}; that is each $W^{i,j}_{\sigma,+}$ is isomorphic to $V_\sigma\otimes V_\sigma$ as a $F\wr C_2$ representation. Analogously, for $i\neq j$ we can define

\[
W^{i,j}_{\sigma,-}=\left\langle \left(T_i^\sigma v_1\otimes T_j^\sigma v_2\right)-
\left(T_j^\sigma v_1\otimes T_i^\sigma v_2\right)
:v_1,v_2\in V_\sigma\right\rangle
\]

that corresponds to the choice $\eta$= the nontrivial representation in \eqref{That}; that is each $W^{i,j}_{\sigma,-}$ is isomorphic to $(V_\sigma\otimes V_\sigma)\otimes U$ as a $F\wr C_2$ representation. Then the decomposition of $L(Y\times Y)$ into irreducible $F\wr C_2$-representations is precisely:

\[
L(Y\times Y)=\left(\bigoplus_{\substack{\sigma,\sigma'\in R\\\sigma\neq\sigma'}}\bigoplus_{i=1}^{m_\sigma}\bigoplus_{j=1}^{m_{\sigma'}}W^{j,j'}_{\sigma,\sigma'}\right)\bigoplus\left\{\bigoplus_{\sigma\in R}\left[\left(\bigoplus_{i,j=1}^{m_\sigma}W^{i,j}_{\sigma,+}\right)\bigoplus
\left(\bigoplus_{\substack{i,j=1\\ i\neq j}}^{m_\sigma}W^{i,j}_{\sigma,-}\right)\right]\right\}.
\]

In more simple terms, every representation $(V_\sigma\otimes V_{\sigma'})\oplus( V_{\sigma'}\otimes V_\sigma)$, $\sigma\neq\sigma'$, appears with multiplicity $m_\sigma m_{\sigma'}$, every representation $V_\sigma\otimes V_\sigma$ appears with multiplicity $\frac{m_\sigma(m_\sigma+1)}{2}$ and every representation $(V_\sigma\otimes V_\sigma)\otimes U$ appears with multiplicity $\frac{m_\sigma(m_\sigma-1)}{2}$ .

%%%%%%%%%%%%%%%%%%%%%%%%%%%%%%%%%%%%%%%%%%%%%%%%%%%%%%%%%%%%%%%%%%%%%%%%%%
%%%%%%%%%%%%%%%%%%%%%%%%%%%%%%%%%%%%%%%%%%%%%%%%%%%%%%%%%%%%%%%%%%%%%%%%%
\subsection{The case in which $L(Y)$ is multiplicity free}
%%%%%%%%%%%%%%%%%%%%%%%%%%%%%%%%%%%%%%%%%%%%%%%%%%%%%%%%%%%%%%%%%%%%%%%%%%
%%%%%%%%%%%%%%%%%%%%%%%%%%%%%%%%%%%%%%%%%%%%%%%%%%%%%%%%%%%%%%%%%%%%%%%%%

In this subsection, we examine the particular case in which $L(Y)$
decomposes without multiplicity. Suppose that
$L(Y)=\bigoplus\limits_{h=0}^n V_h$ is the decomposition of $L(Y)$
into inequivalent irreducible representations. We identify each
$V_h$ with a subspace of $L(Y)$; if $v\in V_h$ then $v$ is a
function defined on $Y$ and $v(y)$ is the value of $v$ on $y\in Y$.
Moreover, we denote by $\sigma_h$ the representation of $F$ on
$V_h$; this means that for any $f\in F$ we have a unitary operator
$\sigma_h(f):V_h\rightarrow V_h$ such that
$v(f^{-1}f_0)=[\sigma_h(f)v](f_0)$ for all $v\in V_h,f_0\in F$.

Let $H$ be the set of all functions $h:X\rightarrow\{0,1,\dotsc,n\}$. If $\bigotimes\limits_{x\in X}v_x\in L(Y^X)$ and  $v_x\in V_{h(x)}$ for all $x\in X$, with $h\in H$, we say that $\bigotimes\limits_{x\in X}v_x$ is a vector of {\em type} $h$. Set $V_h=\bigotimes_{x\in X}V_{h(x)}$. Then $V_h$ is the set of all vectors of type $h$ and

\[
L(Y^X)=\bigoplus\limits_{h\in H} V_h
\]

is the decomposition of $L(Y^X)$ into irreducible $F^X$-representations. If $u\in L(Z)$ and $g\in G$ then we denote by $gu$ its $g$-translate, that is $(gu)(z)=u(g^{-1}z)$ for all $z\in Z$. The following Lemma is a variation of Lemma \ref{actionpsi}, but it is more specific to the present situation.

\begin{lemma}\label{actionVhLZ}
If $(f,g)\in F\wr G$, $h\in H$, $\bigotimes\limits_{x\in X}v_x\in V_h$ and $u\in L(Z)$ then

\[
(f,g)\left[\left(\bigotimes_{x\in X}v_x\right)\bigotimes u\right]=\left\{\bigotimes_{x\in X}[\sigma_{h(g^{-1}x)}(f(x))v_{g^{-1}x}]\right\}\bigotimes gu.
\]

\end{lemma}

\begin{proof}
If $(\varphi,z)\in Y^X\times Z$ then

\[
\begin{split}
\left\{(f,g)\left[\left(\bigotimes_{x\in X}v_x\right)\bigotimes u\right]\right\}(\varphi,z)=&
\left[\left(\bigotimes_{x\in X}v_x\right)\bigotimes u\right]((f,g)^{-1}\varphi,g^{-1}z)=\\
=&\left\{\prod_{x\in X}v_x[f(gx)^{-1}\varphi(gx)]\right\}\cdot (gu)(z)=\\
=&\left\{\prod_{x\in X}[\sigma_{h(g^{-1}x)}(f(x))v_{g^{-1}x}](\varphi(x))\right\}\cdot (gu)(z)=\\
=&\left(\left\{\bigotimes_{x\in X}[\sigma_{h(g^{-1}x)}(f(x))v_{g^{-1}x}]\right\}\bigotimes gu\right)(\varphi,z).
\end{split}
\]
\end{proof}

As in subsection \ref{expwr}, fix a $\sigma\in\Sigma$. Now it means
that we have fixed an $h\in H$. Moreover, now $I=I_G(\sigma)$ is
just the stabilizer of $h$ ($G$ acts on $H$ in the obvious way).
Suppose that $\eta$ is an irreducible $I$-representation contained
in the permutation representation on $Z$ and that $U_1\oplus
U_2\oplus\dotsb\oplus U_m$ is an orthogonal decomposition of the
$\eta$-isotypic component in $L(Z)$. We think of each $U_i$ as a
subspace of $L(Z)$. Suppose that $S$ is a system of representatives
for the right cosets of $I$ in $G$.

\begin{theorem}
For $i=1,2,\dotsc,m$, set

\[
W_i=\bigoplus_{s\in S}\left(V_{sh}\otimes sU_i\right).
\]

Then each $W_i$ is $F\wr G$-invariant, irreducible and equivalent to $\text{\rm Ind}_{F\wr I}^{F\wr G}(\widetilde{\sigma}\otimes\overline{\eta})$. Moreover,

\[
W_1\oplus W_2\oplus\dotsb\oplus W_m
\]

is an orthogonal decomposition of the $\text{\rm Ind}_{F\wr I}^{F\wr G}(\widetilde{\sigma}\otimes\overline{\eta})$-isotypic component of $L(Y^X\times Z)$.
\end{theorem}

\begin{proof}
Now $J$ is trivial and \eqref{That} becomes $\widehat{\mathcal{T}}(v\otimes u)=v\otimes \mathcal{T}u$.
Therefore Theorem \ref{decFwrI} just tells us that

\[
\bigoplus_{i=1}^m\left[\biggl(\bigotimes_{x\in X}V_{h(x)}\biggr)\bigotimes U_i\right]
\]

is an orthogonal decomposition of the $\widetilde{\sigma}\otimes\overline{\eta}$-isotypic component in $L(Y^X\times Z)$. Note that $\{(\mathbf{1}_F,s):s\in S\}$ is a system of representatives for the right cosets of $F\wr I$ in $F\wr G$. Then an application of Proposition \ref{Frobenius} by mean of the identity in Lemma \ref{actionVhLZ} ends the proof: it suffices to notice that
if $s\in S$, $\bigotimes\limits_{x\in X}v_x\in V_h$ and $u\in U_i$ then

\[
(\mathbf{1}_F,s)\left[\biggl(\bigotimes\limits_{x\in X}v_x\biggr)\bigotimes u\right]=\biggl(\bigotimes\limits_{x\in X}v_{s^{-1}x}\biggr)\bigotimes su
\]

is an element of $V_{sh}\otimes sU_i$.

\end{proof}

\begin{corollary}
The multiplicity of $\text{\rm Ind}_{Fwr I}^{F\wr G}\left(\widetilde{\sigma}\otimes\overline{\eta}\right)$ in $L(Y^X\times Z)$ is equal to the multiplicity of $\eta$ in $L(Z)$.
\end{corollary}

\begin{remark}
{\rm
The representation theoretic results in \cite{AM,Mi,ScaTol2,Schoolfield} (and the exponentiation of a finite Gelfand pair is studied in \cite{CST2}) are all particular cases of this example.
}
\end{remark}

\end{document}